\documentclass[10pt,a4paper]{article}
\usepackage[pdftex]{graphicx} %\usepackage[dvips]{graphicx} \usepackage{psfrag} % incompatibile con dvipdfm
\usepackage{amsmath, amsfonts, amssymb, amsbsy}
\usepackage{latexsym} \usepackage{stmaryrd} %\usepackage{MnSymbol} 

\usepackage[labelfont=bf]{caption}
\usepackage{hyperref}
\hypersetup{
    colorlinks,%
    citecolor=black,%
    filecolor=black,%
    linkcolor=black,%
    urlcolor=black
}

\begin{document}

%%%%%%%%%%%%%%%%%%%%%%%%%%%%%%%%%%%%%%%%%%%% MACROS
%%%%%%%%%%%%%%%%%%%%%%%%%%%%%%%%%%%%%%%%%%%%%% NUMERAZIONE
\newtheorem{theorem}{Theorem}[section]
\newtheorem{proposition}[theorem]{Proposition}
\newtheorem{corollary}[theorem]{Corollary}
\newtheorem{lemma}[theorem]{Lemma}
\newtheorem{definition}[theorem]{Definition}
\newtheorem{remark}[theorem]{Remark}
\newtheorem{example}[theorem]{Example}

%%%%%%%%%%%%%%%%%%%%%%%%%%%%%%%%%%%%%%%%%%%%%%% MACRO TESTO
\newcommand{\proof}{\noindent \textbf{Proof. }}
\newcommand{\qed}{ \hfill {\vrule width 6pt height 6pt depth 0pt} \medskip }

\newcommand{\separe}{\medskip \centerline{\tt -------------------------------------------- } \medskip}
\renewcommand{\separe}{}

\newcommand{\doit}[1]{\hfill {\small \tt [#1]}}
\newcommand{\note}[1]{\noindent \fbox{\parbox{\textwidth}{\tt #1}} \medskip}

%%%%%%%%%%%%%%%%%%%%%%%%%%%%%%%%%%%%%%%%%%%%%%%% MACRO MATEMATICA
\newcommand{\eps}{\varepsilon} \renewcommand{\det}{\mathrm{det}} \newcommand{\argmin}{ \mathrm{argmin} \,}
\newcommand{\Om}{\Omega} \def\interior{\mathaccent'27} 
\newcommand{\weakto}{ \rightharpoonup}  \newcommand{\weakstarto}{\stackrel{*}{\rightharpoonup}}
\newcommand{\R}{\mathbb{R}}

\newcommand{\stress}{\boldsymbol{\sigma}} \newcommand{\strain}{\boldsymbol{\epsilon}} 
\newcommand{\neu}{ \partial_{\mbox{\it \tiny N}} } \newcommand{\dir}{\partial_{\mbox{\it \tiny D}} }
\newcommand{\I}{ {\mbox{\tiny \rm I}} } \newcommand{\II}{ {\mbox{\tiny \rm II}} }
\newcommand{\jump}[1]{\llbracket #1 \rrbracket}
\def\Xint#1{\mathchoice
 {\XXint\displaystyle\textstyle{#1}}%
 {\XXint\textstyle\scriptstyle{#1}}%
 {\XXint\scriptstyle\scriptscriptstyle{#1}}%
 {\XXint\scriptscriptstyle\scriptscriptstyle{#1}}%
 \!\int}
\def\XXint#1#2#3{{\setbox0=\hbox{$#1{#2#3}{\int}$}
 \vcenter{\hbox{$#2#3$}}\kern-.5\wd0}}
 \def\ddashint{\Xint=} %%%%%%%%%%%%%%%%%%%%%% questi sono i comandi 
 \def\dashint{\Xint-}

 % use symbolic links

\newcommand{\F}{\mathcal{F}}\newcommand{\E}{\mathcal{E}}

%%%%%%%%%%%%%%%%%%%%%%%%%%%%%%%%%%%%%%%%%%%%% ARTICLE 

% !TeX root = article.tex

%%%%%%%%%%%%%%%%%%%%%%%%%%%%%%%%%%%%%%%%%%%%%%%%%%%% COVER 

\thispagestyle{empty} %\setcounter{page}{0}

\phantom{1}

\vspace{0.5cm}
{\LARGE
{\bf Weak solutions for gradient flows }

\medskip
{\bf under monotonicity constraints}
}
\vspace{36pt}

\begin{small}
{\bf M.~Negri}

%\medskip
{Department of Mathematics -  University of Pavia} 

{Via A.~Ferrata 1 - 27100 Pavia - Italy}

{matteo.negri@unipv.it} %\quad {\tt http://www-dimat.unipv.it/~negri/}

\medskip
{and} 
\medskip

{Institute for Applied Mathematics and Information Technologies - CNR}

{Via A.~Ferrata 5 - 27100 Pavia - Italy}

\vspace{24pt}
{\bf M.~Kimura}

%\medskip
{Faculty of Mathematics and Physics - Kanazawa University}

{Kanazawa 920-1192 - Japan}

{mkimura@se.kanazawa-u.ac.jp}

\vspace{36pt}
\noindent 
{\bf Abstract.} We consider the gradient flow of a quadratic non-autonomous energy under monotonicity constraint in time and natural regularity assumptions. We provide first a notion of weak solution, inspired by the theory of curves of maximal slope, and then existence (employing time-discrete schemes with different ``implementations'' of the constraint), uniqueness, power and energy identity, comparison principle and continuous dependence. As a byproduct, we show that the energy identity gives a selection criterion for the (non-unique) evolutions obtained by other notions of solutions. We finally show that, for autonomous energies, the solutions obtained with the monotonicity constraint actually coincide with those obtained with a fixed obstacle, given by the initial datum. 
%acle constraint, given by the initial condition. 


\end{small}

 %\thispagestyle{empty}
% !TeX root = article.tex
% ============================================================================================  INTRO
% motivation 
% our pb 
% alternative notion of solutions 
% phase-field fracture?
% our results 

\section{Introduction} 

% ------------------------------------------------------------ MOTIVATION 
Parabolic evolution equations with monotonicity constraints naturally arise in several mathematical models; for instance, the behaviour of materials undergoing inelastic processes (like fracture, damage, plasticity etc.) requires monotonicity constraints, due to the irreversibility of such phenomena. Few specific results have been recently obtained for applications in mechanics, see e.g. \cite{BabadjianMillot_AHHPANL14,   BonettiFreddiSegatti_CMT17, Negri_ACV19}; for abstract evolutions let us mention \cite{AkagiKimura_JDE19}, dealing with strong solutions, and the recent \cite{AkagiEfendiev_EJAM}, dealing with the an autonomous Allen-Cahn equation. In the context of \cite{AkagiKimura_JDE19}, and inspired by the applications 
%In this context, our work aims at both facing basic theoretical questions and providing applicable results.
%for this class of evolutions. in between theory and applications for this class of evolutions, our work faces basic theoretical questions and provides applicable result. 
%
% like well-posedness, continuous dependece) %looking both toward theoretical understanding and ``practical'' 
%
% and facing this class of evolutions is not straightforward as it may seem at first glance; for instance, the available notions of solution are not equivalent or not always applicable.
%
%time-depending and time-independent data make a major technical difference both in the proofs and (even) in the set and notion of solutions.
%
%
%------------------------------------------------------------- OUR PB
%\medskip
%Aiming at and having in mind applications of the phase-field approach, %, our intention is to study parabolic evolution equations with monotonicity constraints. 
%To this end, 
we consider in particular a prototype energy of the form 
$$ 
\F (t,u) = \E (u) - \langle f(t) , u \rangle = \tfrac12 a (u,u) - \langle f(t) , u \rangle 
$$ 
where $a(\cdot, \cdot)$ is a coercive, continuous bi-linear form in $H^1_0$ while $f$ belongs to $L^2(0,T; L^2)$. 
 We write $\langle \cdot , \cdot \rangle$ for the $L^2$-scalar product and $(\cdot , \cdot)$ for the duality between $H^1_0$ and $H^{-1}$, we will also employ the operator $A:H^1_0 \to H^{-1}$ defined by $ a ( u, v) = - ( A u ,v )$. We consider weak solutions $u \in L^\infty (0,T; H^1_0) \cap H^1 ( 0, T ; L^2)$ with $u(0) = u_0 \in H^1_0$.  
Before switching to the mathematical content let us make a comment on time-depending data in applications: in the context of phase-field models for fracture stored energies often take the form 
$$
     \F_\eps ( t , v , w) = \tfrac12 \int_\Omega \eps | \nabla v |^2 \, + \eps^{-1} |v |^2 \, dx + \int_\Omega \mathcal{W} (  t, v , w) \, dx 
$$
where $v \in H^1(\Omega , [0,1])$ is the phase-field variable, $w$ is the displacement field, while $\mathcal{W}  ( t, v, w)$ is the elastic (phase-field) energy. The first integral plays the role of the energy $\E$, while the second (non-linear) term corresponds, roughly speaking, to the (linear) term $ \langle - f(t) , v \rangle$; indeed, among the many (see for instance \cite{Wu_JMPS17} and the references therein) a possible, simple choice is $\mathcal{W} ( t, v, w) = (v-1) W ( Dw(t) )  $ where $W$ denotes linear elastic energy density, whose regularity in time is, in general, not better that $L^\infty(0,T; L^p)$ for some $p < 1$, see e.g.~\cite{KneesRossiZanini_M3AS13, AN}.  In this specific application, the fact that $p<2$ is balanced by the fact that $v \in L^\infty$ which is not the case in our setting, however it is important to note that differentiability of $W$ in time is out of reach. We anticipate that the time regularity of data and solutions will play a crucial role also in the analysis.  %In our work we face basic theoretical questions and provide results of interest in the applications.
%for this class of evolutions. in between theory and applications for this class of evolutions, our work faces basic theoretical questions and provides applicable result. 
%
% like well-posedness, continuous dependece) %looking both toward theoretical understanding and ``practical'' 
%
% and facing this class of evolutions is not straightforward as it may seem at first glance; for instance, the available notions of solution are not equivalent or not always applicable.
%

\medskip

%In our work aims at both facing basic theoretical questions and providing applicable results.
%for this class of evolutions. in between theory and applications for this class of evolutions, our work faces basic theoretical questions and provides applicable result. 
%
% like well-posedness, continuous dependece) %looking both toward theoretical understanding and ``practical'' 
%
% and facing this class of evolutions is not straightforward as it may seem at first glance; for instance, the available notions of solution are not equivalent or not always applicable.
%

Our very first target is a suitable notion of solution and, equivalently, an effective way of writing the unilateral (constrained) gradient flow. This basic question is delicate, in particular as far as well-posedness, since different notions may provide different solutions; let us briefly list the main options available in the literature (more details are in \S 1), highlighting the possible issues and the main differences.

%\medskip

%------------------------------------------------------------- ALTERNATIVE NOTION OF SOLUTIONS 
%---------- parabolic variational inequalities [cfr. Brezis]
A convenient framework to handle unilateral (monotonicity) constraints is given by parabolic variational inequalities. In our setting, it is natural to search for solutions $u$ such that for a.e.~$t \in (0,T)$ it holds 
\begin{equation} \label{e.0}
	(A u (t) , z - u(t) ) + \langle f(t) , z - u(t) \rangle \le \langle \dot{u} (t) , z - u(t) \rangle 
\end{equation}
for every $z$ in the convex cone $K(t) = \{ z \in H^1_0: z \ge u(t) \}$. This is an elliptic-parabolic problem with a  time depending constraint, cf. \cite{KuboShirakawaYamazaki_JMAA12, FukaoKenmochi_AMSA13}; however here $K(t)$ is not a datum, because it depends on the solution $u$ itself, and this changes significantly the problem. Indeed, existence is easily proved (see Proposition \ref{p.pde}) but the set of solutions turns out to be but far too large for uniqueness (see the counter-example in \S \ref{non-uni}). 

% ---------  doubly non-linear ?? ... Goro
Another possible way of writing the evolution is to employ a sort of ``doubly non-linear inclusion''; in our setting a feasible formulation could be 
\begin{equation} \label{e.1}
%	\begin{cases}
		\dot{u}(t) + \partial I_+ (\dot{u}(t)) - A u(t) - f(t) \ni 0 % \\
%\partial \Psi (\dot{u}(t)) \ni A u(t) + f(t) & \text{in $L^2$ for a.e.~$t \in (0,T)$,} \\
%		u (0) = u_0 ,
%	\end{cases}
\end{equation}
where $I_+$ is the indicator function of the set $\{ v \ge 0\}$ and $\partial I_+$ is its $L^2$-subdifferential. However, this inclusion to hold, $A u(t) + f(t)$ should be in $L^2$, which is not true in general (cf. \S\ref{PDE}). As a matter of fact, this approach is suitable under more restrictive conditions on $u_0$ and $f$, which ensure $ u(t) \in H^2$ and thus $A u(t) + f(t)$ in $L^2$, see e.g.~\cite{AkagiEfendiev_EJAM, AkagiKimura_JDE19}.

%----------  gradient flows
Another natural approach is to consider the $L^2$-projection of the gradient of the energy $\F$. For sake of simplicity, note that if $u \in H^2$ then $-d \F ( t , u )  [z] = \langle A u + f(t) , z \rangle$, as a consequence, the positive part $[ A u + f(t) ]_+$ is the $L^2$-projection of $Au + f(t) = -\nabla_{\!L^2} \F ( t, u)$ on the cone of positive functions. Therefore, it makes sense to search for solutions $u$ such that 
\begin{equation} \label{e.2}
%	\begin{cases}
		\dot{u} (t) = [ A u (t) + f (t) ]_+  \quad \text{in $L^2$ for a.e.~$t \in (0,T)$.} %\\
%		u (0) = u_0 ,
%	\end{cases}
\end{equation}
Technically, if $u (t) \not\in H^2$ then $A u (t) + f (t)$ is a locally finite Radon measure and $[A u (t) + f (t)]_+$ is its positive part, in the sense of Hahn decomposition. Once again, existence of solutions is true (see Proposition \ref{p.pde}) but uniqueness is not (see the counter-example in \S \ref{non-uni}).

% On the other hand, doubly non-linear evolutions ... [cfr. Colli-Visintin, Goro, Gianazza-Savare']

% ---------  curve of maximal unilateral slope 
Finally, let us introduce our notion of solution, which provides existence and uniqueness. We employ the theory and the language of curves of maximal slope \cite{AmbrosioGigliSavare05}, starting, for sake of clarity, with $f \in AC(0,T;L^2)$.  In this case (see Theorem \ref{t.exist}) there exists a unique $u$ such that the energy $t \mapsto \F ( t , u(t))$ is absolutely continuous in $(0,T)$ and for a.e.~$t \in (0,T)$ the following power balance holds  
\begin{equation} \label{e.3}
	\dot{\F} ( t, u(t) ) %=  d \F ( t , u(u) ) [ \dot{u} (t) ] + \partial_t \F ( t ,u(t)) 
\le - \tfrac12 | \dot{u} (t) |^2_{L^2} \,  - \tfrac12 | \partial \F |^2_{L^2_+} ( t, u (t) ) -  \langle \dot{f} (t) , u(t) \rangle ,
\end{equation}
where $| \dot{u} |_{L^2_+}$ and $| \partial \F |_{L^2_+} (t, u)$ denote respectively a singular (unilateral) norm and the unilateral slope respectively given by 
\begin{equation*} %\label{e.slope}
	| v |_{L^2_+} = \begin{cases}
		\| v \|_{L^2} & \text{if $v\ge0$,} \\
		+ \infty & \text{otherwise},
	\end{cases}
	\qquad 
 | \partial \F |_{L^2_+} (t, u) = %\limsup_{\substack{v \to u \\ v \text{\tiny $\ge$} u}}  
\limsup_{v \to u} \frac{[\F(t, v) - \F(t, u)]_-}{ | v - u |_{L^2_+}} = \| [ A u (t) + f (t) ]_+ \|_{L^2} .
%	\quad	\text{ with $v \ge u$,}
\end{equation*}
Following \cite{AmbrosioGigliSavare05} we will say that $u$, satisfying \eqref{e.3}, is a curve of maximal unilateral slope. When $f \in L^2(0,T;L^2)$ the power balance inequality \eqref{e.3} does not make sense since the time derivative of $f$ is not available. However, there exists a unique $u$ such that 
\begin{equation} \label{e.4}
      \dot{\E} ( u (t) ) 
%			& \le & - \tfrac12 | \dot{u} (t) |_{L^2_+}^2 - \tfrac12 | \partial \F |^2_{L^2_+} (t, u (t))\, + \langle f (t) , \dot{u} (t) \rangle  \\
			\le  - \tfrac12 | \dot{u} (t) |_{L^2_+}^2 - \tfrac12 \| [ A u(t) + f(t)]_+ \|^2_{L^2_+} , + \langle f (t) , \dot{u} (t) \rangle .
\end{equation}
Actually, in Definition \ref{d.curve} we will employ an equivalent time-integral formulation which is more convenient in the proofs and which is strictly related to the energy identity 
\begin{equation} \label{e.4bis}
       \E ( u (t) )  = \E ( u_0 )  -  \, \tfrac12 \int_{0}^{t} | \dot{u} (s) |_{L^2_+}^2 + \| [ A u (s) + f(s) ]_+ \|^2_{L^2} \, ds + \int_{0}^{t} \langle f (s) , \dot{u} (s) \rangle \, ds  \quad \text{for every $t \in (0,T)$. } \\
\end{equation} 
%More precisely, \eqref{e.4} and \eqref{e.4bis} are equivalent. 
At this point, it is important to remark that the unique solution of \eqref{e.4} is also a solution to \eqref{e.0} and \eqref{e.2}; in other terms, the energy balance turns out to select a unique solution of the parabolic variational inequality \eqref{e.0} and of the unilateral gradient flow \eqref{e.2}. Moreover, if the solution is sufficiently regular then it solves also \eqref{e.1} and \eqref{e.3}.

% ------------------------------------------------------------ PHASE-FIELD FRACTURE?

% ------------------------------------------------------------ OUR RESULTS 
\medskip
% --------------  well posedness

% -------------- discrete schemes 
Now, let us describe the structure and the content of the article. Sections and results are organized according to the time regularity of the datum $f$, which plays an important role both in the analysis and in the applications. First of all we consider the most general case, i.e.~$f \in L^2(0,T;L^2)$, which occupies most of the paper. We prove existence and uniqueness of a solution in the sense of \eqref{e.4}. 
Existence is obtained by time discretization, employing three different incremental problems of interest in the applications \cite{TakaishiKimura_K09, AlmiBelzNegri_M2AN19, GerasimovDeLorenzis_18}. Let $t_{n,k} = k \tau_n$ is a uniform discretization of the interval $[0,T]$ with $\tau_n = T/ n$. In the first scheme, given $u_{n,k}$ at time $t_{n,k}$, let the configuration $u_{n,k+1} $ at time $t_{n,k+1}$ be simply given by 
\begin{equation*} %\label{e.scheme1}
	u_{n,k+1} \in \argmin \big\{ \F ( t_{n,k+1} , u  ) + \tfrac{1}{2\tau_n} \, | u - u_{n,k} |^2_{L^2_+} :  u \in H^1_0 \big\}  .
\end{equation*}
In the second we employ instead an a posteriori truncation, i.e., given $u_{n,k}$ we define $u_{n,k+1} $ by
\begin{gather*} %\label{e.scheme2}
	\begin{cases} \tilde u_{n,k+1}  \in \argmin \big\{ \F ( t_{n,k+1} , u  ) + \tfrac{1}{2\tau_n} \, \| u - u_{n,k} \|^2_{L^2} :  u \in H^1_0 \big\}  \\
                                    u_{n,k+1} = \max \big\{  \tilde u_{n,k+1} , u_{n,k}  \big\} . 
            \end{cases}
\end{gather*}
The fact that the first minimization is unconstrained makes this scheme very convenient in the numerical implementation \cite{TakaishiKimura_K09,AlmiBelzNegri_M2AN19}, on the other hand the analysis is slightly more involved. Last, we consider a penalty method, i.e., given $u_{n,k}$ we get $u_{n,k+1}$ by solving  
\begin{gather*} %\label{e.scheme3}
	u_{n,k+1} \in \argmin \big\{ \F ( t_{n,k+1} , u  ) + \tfrac{1}{2\tau_n} \, | u - u_{n,k} |^2_{L^2_{\tau_n}} :  u \in H^1_0 \big\}  .
\end{gather*}
where 
$$
	| v |^2_{L^2_{\tau_n}} = \int_\Omega \psi_{\tau_n} (v) \, dx \quad \text{and} \quad  \psi_{\tau_n} ( v ) = 
		\begin{cases}
		v^2 & \text{if $v \ge 0$} \\
		\alpha_n v^2 & \text{if $v < 0$}
		\end{cases} 
	\quad \text{for $\alpha_n \to +\infty$.}
$$
Each scheme defines a sequence of discrete solutions $u_n$ (depending on $\tau_n$) which  enjoys suitable compactness properties and which converges (weakly and up to subsequences) to a solution of \eqref{e.4}; a posteriori we actually prove that the whole sequence converges strongly.

Note that for $f \in L^2(0,T;L^2)$ the time regularity of solutions is rather low, since in general $u \in H^1 (0,T;L^2)$. As a consequence, uniqueness does not follow from classical tools, we use instead a contradiction argument of \cite{Gigli_CVPDE10} based on energy balance and convexity. For the same reason, the energy identity does not follow by the chain rule, which would require at least $u \in H^1_{loc} (0,T; H^1)$, rather, it is proved employing a measure theory argument, see also \cite{DalMFrancfToad05, Negri_ACV19}. 

In the second case we consider $f \in AC (0,T;L^2)$. This is obviously contained in the previous one, however, from a theoretical point of view, it is interesting to know that in this case solutions are of class $H^1_{loc} (0,T; H^1_0)$; as a consequence, a better representation holds and few issues, due to the lack of time regularity, are avoided.

Lastly, when $f$ is independent of time, besides recovering the classical results of \cite{GianazzaSavare_RANSXLMM94}, we prove a (rather surprising) property: the unique solution of \eqref{e.4} turns out to coincide with the unique solution of the unconstrained $L^2$-gradient flow for the functional $\tilde{\F} (u) = \F (u) + I_+ (u -u_0)$, in other terms, the monotonicity constraint can be replaced by a fixed obstacle, given by the initial datum $u_0$; however, this property does not hold when $f$ depends on time (see the counter-example in Remark \ref{r.auto}).

To complete our analysis, we prove a comparison principle and a (non-quantitative) continuous dependence property for solutions of \eqref{e.4}; moreover, for the interested reader, we provide in the appendix  further properties, representations and remarks on the unilateral slope. Finally, we remark that several results can be generalized, for instance to (non-quadratic) convex or $\lambda$-convex energies (see \S \ref{s.gen}). 

% -------------- autonomous functionals and dependence on time

% ===========================================================================================  ABSTRACT 

%Moreover, we provide a comparative study, highlighting equivalences, issues, and above all differences between different notion of solutions. 
%Beside the theoretical nature of this article, the results, in our intention, will serve as a reference for the choice of the proper technical setting in future applications.

% ===========================================================================================  EXTRA

%On the technical side, the analysis of this peculiar class of PDEs could rely on several technical approaches: parabolic (or elliptic-parabolic) variational inequalities, doubly non-linear evolution equations and gradient flows. These theories, dating back to the 70's, are classical and can be adapted to monotonicity constraints, however, on the base of a case study under development, they may provide different sets of solutions, for instance when monotonicity constraints are applied to non-autonomous PDEs. 

\tableofcontents

%\bigskip
%{\note{Possible interesting problems: \break 
%%1) $f (t,u)$ and semi-implicit scheme, see alternate minimization \break 
%%- existence for p-laplacian energy, see [KRZ], or $BV$, see [ST] and [SK] \break
%%*) check if and how the re-formulation with fixed obstacle of Akagi-Efendiev fits into our framework \break
%*) continuous dependence of solutions in quantitative form \break
%**) Neumann problem 
%}

%\note{Changed time to more natural $t_{n,k+1}$ in the discrete scheme and made notation consistent throughout paper.}

%\pagebreak 

%%%%%%%%%%%%%%%%%%%% SECTIONS
% !TeX root = article.tex

\section{Setting and statement of the main results\label{state}}

Let us consider an open, bounded domain $\Omega \subset \R^d$. Throughout the paper we will employ the short-hand notation $L^2$ for $L^2(\Omega)$ and similarly for other functional spaces. We will use the notation $\langle \cdot , \cdot \rangle$ for the scalar product in $L^2$, while $( \cdot, \cdot)$ will denote the duality between $H^{-1}$ and $H^1_0$. Consider a coercive, continuous and symmetric bi-linear form $a(\cdot, \cdot)$ in $H^1_0 \times H^1_0$ given by 
$$
	a ( u ,v) = \int_\Omega \nabla u (x)  \cdot B(x) \,\nabla v(x) + b(x) \, u(x) v(x) \, dx 
$$
and the corresponding operator $A: H^1_0 \to H^{-1}$ given by $(A u , v) = - a ( u , v)$. Accordingly, we introduce the stored energy 
\begin{equation} \label{e.en-stored}
   \E ( u ) = \tfrac12 a (u , u) = \tfrac12 \int_\Omega \nabla u \cdot B \,\nabla u + b \, u^2 \, dx .
\end{equation}
Clearly $ a^{1/2} (u,u)$ is the energy-norm which is equivalent to the standard norm in $H^1_0$.

Finally, let us introduce the following convenient notation
$$
	| w |_{L^2_+} = \begin{cases}
		\| w \|_{L^2} & \text{if $w\ge0$,} \\
		+ \infty & \text{otherwise}.
	\end{cases}
$$
Accordingly, we will say that $v \to u$ in $L^2_+$ when $| v - u |_{L^2_+} \to 0 $, i.e.~when $v \ge u$ and $v \to u$ in $L^2$.

\subsection[The case $f \in L^2(0,T;L^2)$]{The case \boldmath{$f \in L^2(0,T;L^2)$}}

Let $[0,T]$ be a time interval and let $ f \in L^2 ( 0,T ; L^2)$. Let us choose a representative of $f$ (defined for every $t \in [0,T]$) and consider the free energy $\F : [0,T] \times L^2 \to \R \cup \{ +\infty \}$ given by 
\begin{equation}\label{e.en}
	\F ( t, u ) = \begin{cases}
	{\displaystyle \E(u) - \langle f(t) , u \rangle } &  \text{for } u \in H^1_0 , \\
	+ \infty & \text{otherwise}
	\end{cases}
\end{equation}
(in the sequel we will see that the evolution is independent of the choice of the representative).
%Often, we will write $\F ( t ,u ) = \E (u) - \langle f(t) , u \rangle$.
Clearly, for $u,v \in H^1_0$ we have
$$
	d \F ( t, u) [v] = a ( u , v ) - \langle f(t) , v \rangle = - ( A u + f(t) , v ) .
$$
In particular, for $u \in H^2$ we have $Au \in L^2$ and thus $-d \F ( t , u )  [z] = \langle A u + f(t) , z \rangle$, then we can write $ - \nabla_{\!L^2} \F ( t, u) = A u + f(t)$ and thus $[ A u + f(t) ]_+$ (the positive part) turns out to be the $L^2$-projection of $- \nabla_{\!L^2} \F ( t, u)$ on the cone of positive functions, which is indeed the set of admissible variations. A qualitatively similar property holds, in a suitable sense, even if $u \in H^1_0 \setminus H^2$, see \S\,\ref{en-sl}. Inspired by the theory of curves of maximal slope \cite{ AmbrosioGigliSavare05} we provide the following definition (further connections will be given in the sequel).

%\begin{example} The simplest possible example is given by 
%$$
%	 \F (t, u) = \tfrac12 \int_\Omega | \nabla u |^2 \, dx - \int_\Omega f (t) \, u \, dx  .
%$$
%In this case $A u = - \Delta u$. 
%\end{example}

%\note{%Possibly $f \in L^q$ for some $q$ by Sobolev embedding. \hfill \break 
%In general $f(t)$ and $g ( \nabla u)$ and Dirichlet or Neumann b.c. \hfill \break 
%Possibly we could assume that $\F$ is $\lambda$-convex, coercive and lsc. \hfill \break
%%Show that $\lim_{t \to +\infty} u(t) \in \argmin \{ \F(u) : u \ge u_0 \}$
%Speed of convergence at $\infty$ ...~with different algorithms
%}

\begin{definition} \label{d.curve} An evolution $u \in L^\infty (0,T; H^1_0) \cap H^1 ( 0, T ; L^2)$ is a unilateral gradient flow for the energy $\F$ if $A u(t) + f(t)$ is a Radon measure for a.e.~$t \in (0,T)$ and if for every $0 \le t^* \le T$
%%$\F ( u ( \cdot))$ is monotone non-increasing in $[0,T)$ and if 
\begin{equation}\label{e.evol-ineq}
      \E ( u (t^*) ) \le \E ( u (0) )  -  \, \tfrac12 \int_{0}^{t^*} | \dot{u} (t) |_{L^2_+}^2 + \| [ A u (t) + f(t) ]_+ \|^2_{L^2} \, dt + \int_0^{t^*} \langle f (t) , \dot{u} (t) \rangle \, dt.
\end{equation}
The fact that $u$ is monotone in time, i.e.~$\dot{u} \ge 0$, is implicitely written in \eqref{e.evol-ineq}.
\end{definition}

%\note{ EQUIVALENTLY ... if the energy $t \mapsto \E ( u(t))$ is absolutely continuous in $[0,T]$ and for a.e.~$t \in (0,T)$ it holds %if $t \mapsto \E(u(t))$ is absolutely continuous and  
%$$
%    \dot{\E} (u(t)) \le - | \dot{u} (t) |_{L^2_+} \, | \partial \F |_{L^2_+} ( u(t)) + \langle f(t) , \dot{u}(t) \rangle . 
%$$
%}

The next theorem contains the main result: existence, uniqueness, and energy identity; it will be will be proven in \S\,\ref{sec.3}, employing several different time-discrete schemes. 

\begin{theorem}\label{t.exist} Given $u_0 \in H^1_0$ there exists a unique unilateral gradient flow $u$ for $\F$ with  $u(0) =u_0$. Moreover, for every $0 \le t_1 \le t_2 \le T$ the following energy identities hold:
\begin{align}\label{e.evol-eq}
      \E ( u (t_2) ) & = \E ( u (t_1) )  -  \, \tfrac12 \int_{t_1}^{t_2} | \dot{u} (t) |_{L^2_+}^2 + \| [ A u (t) + f(t) ]_+ \|^2_{L^2} \, dt + \int_{t_1}^{t_2} \langle f (t) , \dot{u} (t) \rangle \, dt \\
	& = \E ( u (t_1) ) -  \int_{t_1}^{t_2} | \dot{u} (t) |_{L^2_+} \, \| [ A u (t) + f(t) ]_+ \|_{L^2} \, dt + \int_{t_1}^{t_2} \langle f (t) , \dot{u} (t) \rangle \, dt . \label{e.evol-eqbis}
\end{align}
Note that \eqref{e.evol-eq} is independent of the choice of the representative of the datum $f$.
% since all terms depending on $f$ are integrated over $[t_1, t_2]$.
\end{theorem}

From \eqref{e.evol-eqbis} it follows that the energy $t \mapsto \E ( u(t))$ is absolutely continuous in $[0,T]$ and for a.e.~$t \in (0,T)$ it holds %if $t \mapsto \E(u(t))$ is absolutely continuous and  
$$
    \dot{\E} (u(t)) \le - | \dot{u} (t) |_{L^2_+} \, \| [ A u (t) + f(t) ]_+ \|_{L^2} + \langle f(t) , \dot{u}(t) \rangle .
$$
%(which is indeed equivalent to characterization of unilateral gradient flows.)

%\begin{remark} \label{r.physics} 

%\end{remark} 

Remember that $[A u + f (t)]_+$ plays the role of the projection of $- \nabla_{\!L^2} \F ( t, u)$ on the unilateral cone of positive functions; it is thus natural that unilateral gradient flows solve also the parabolic problem \eqref{e.pde} and the elliptic-parabolic variational inequality \eqref{e.pvi} below. 
% their solutions however are not unique 
Actually, in \S\,\ref{PDE}  we will see that solutions of \eqref{e.pde} or \eqref{e.pvi} are not unique. Therefore, \eqref{e.pde} or \eqref{e.pvi}, by themselves, are not characterizations of unilateral gradient flows, in the sense of Definition \ref{d.curve}. Lack of uniqueness is essentially due to the constraint, indeed, by classical results (see e.g.~\cite{Brezis_73}) the solution of the uncostrained gradient flow for $\F$ would be unique. In other terms, not all the solutions of \eqref{e.pde} or \eqref{e.pvi} satisfy the energy identity \eqref{e.evol-eq}, which instead selects a unique solution. However, if $u_0 \in H^2$ and if $f$ is suitably controlled then \eqref{e.pde} has a unique strong solution, see \cite{AkagiKimura_JDE19}.

\begin{proposition} \label{p.pde} 
Let $u$ be the unilateral gradient flow provided by Theorem \ref{t.exist}, then $u$ solves the parabolic partial differential equation
\begin{equation}\label{e.pde}
	\begin{cases}
		\dot{u} (t) = [ A u (t) + f (t) ]_+  & \text{in $L^2$ for a.e.~$t \in (0,T)$,} \\
		u (0) = u_0 ,
	\end{cases}
\end{equation}
where $A u (t) + f (t)$ is a (locally finite) Radon measure and $[A u (t) + f (t)]_+$ is its positive part. In particular, $\| \dot{u} \|_{L^2} = | \dot{u} |_{L^2_+} = \| [ A u (t) + f (t)]_+ \|_{L^2}$ a.e.~in $(0,T)$.
Moreover, if $u$ solves the parabolic problem \eqref{e.pde} then it solves also the elliptic-parabolic variational inequality
\begin{equation} \label{e.pvi}
	(A u (t) , \phi ) + \langle f(t) , \phi \rangle \le \langle \dot{u} (t) , \phi \rangle 
	\qquad \text{for every $\phi \in H^1_0$ with $\phi \ge 0$.}
\end{equation}
Clearly, writing $\phi = z - u(t)$ the previous inequality reads 
\begin{equation} \label{e.pvi2}
	\int_\Omega \dot{u} (t) (u(t) - z) \, dx + \int_\Omega \nabla u(t) \cdot B \, \nabla ( u(t) - z ) + b \, u (t)  ( u(t) -z ) \, dx \le  \int_\Omega f(t) ( u(t) - z )  \, dx , 
\end{equation}
for every $z$ in the convex cone $K(t) = \{ z \in H^1_0: z \ge u(t) \}$.
\end{proposition} 

%\begin{remark} \label{remo} 

We remark that in \eqref{e.pvi2} the set $K(t)$ is unknown, since it depends on $u$; this is a major difference comparing with elliptic-parabolic problems with time depending constraints, see e.g.~the recent  \cite{KuboShirakawaYamazaki_JMAA12, FukaoKenmochi_AMSA13} and the references therein.

From the ``physical point of view'' solutions in the sense of Definition \ref{d.curve} could be equivalently characterized by
the parabolic problem \eqref{e.pde} together with the energy identity
\begin{equation} \label{e.phys}
	     \E ( u (t^*) ) = \E ( u_0 )  -  \, \int_{0}^{t^*} \mathcal{D} ( \dot{u} (t) ) \, dt + \int_{0}^{t^*} \mathcal{P}_{ext} ( t, \dot{u} (t) ) \, dt ,
\end{equation}
where 
$$
	\mathcal{D} ( \dot{u} (t) ) = \| \dot{u} \|^2_{L^2_+} \quad \text{and} \quad 
	\mathcal{P}_{ext} ( t , \dot{u} (t) ) = \langle f(t) ,  \dot{u} \rangle 
$$
denote respectively the dissipation and the power of external forces.

Finally, notice that in general $A u(t) + f(t) \in H^{-1} \setminus L^2$; an explicit example is given in \S\ref{PDE}. In particular, if $\Psi : L^2 \to [0,+\infty]$ is given by 
$\Psi(v) = \tfrac12 \| v \|^2_{L^2} + I_+ (v)$ (where $I_+$ is the indicator function of the set $\{ v \ge 0\}$) we cannot re-write \eqref{e.pde} in the form 
\begin{equation}\label{e.pde-akagi}
	\begin{cases}
		\partial \Psi (\dot{u}(t)) \ni A u(t) + f(t) & \text{in $L^2$ for a.e.~$t \in (0,T)$,} \\
		u (0) = u_0 ,
	\end{cases}
\end{equation}
because $\partial \Psi \subset L^2$ and thus $A u(t) + f(t)$ should be in $L^2$, which is not always the case. The latter equation, in the form 
$$
	\dot{u}(t) + \partial I_+ (\dot{u}(t)) - A u(t) - f(t) \ni 0 ,
$$
is adopted e.g.~in \cite{AkagiEfendiev_EJAM, AkagiKimura_JDE19} under stronger regularity on the data, in order to have $A u(t) + f(t)$ in $L^2$, see e.g.~\cite[Theorem 2.6]{AkagiKimura_JDE19}. 
%\end{remark}

Unilateral gradient flows, in the sense of Definition \ref{d.curve}, enjoy comparison principle and continuous dependence; on the contrary, by lack of uniqueness, solutions of \eqref{e.pde} or \eqref{e.pvi} do not satisfy them.

\begin{proposition} \label{p.comppri} If $u$ and $v$ are curves of maximal unilateral slope for $\F$ with initial values $u_0 \le v_0$ then $u \le v$ in $[0,T]$.
\end{proposition}

\begin{proposition} \label{p.contdep} Let $f^m \to f$ in $L^2(0,T; L^2)$ and $u_0^m \to u_0$ in $H^1_0$; let  $u^m$ and $u$ be the corresponding unilateral gradient flows. Then $u^m \weakto u$ in $H^1 (0,T;L^2)$ and $u^m (t) \to u(t)$ in $H^1_0$ for every $t \in [0,T]$.
\end{proposition} 

\subsection[The case $f \in AC(0,T;L^2)$]{The case \boldmath{$f \in AC(0,T;L^2)$}}

If $f \in AC ( 0, T; L^2)$ the results of Theorem \ref{t.exist} and Proposition \ref{p.pde} can be improved.
To this end, let us introduce the unilateral $L^2$-slope, defined as follows.

\begin{definition}[\bf Unilateral slope] For $u \in H^1_0$ define 
\begin{equation}\label{e.slope}
 | \partial \F |_{L^2_+} (t, u) = %\limsup_{\substack{v \to u \\ v \text{\tiny $\ge$} u}}  
\limsup_{v \to u} \frac{[\F(t, v) - \F(t, u)]_-}{ | v - u |_{L^2_+}}  ,
%	\quad	\text{ with $v \ge u$,}
\end{equation}
where $[\,\cdot\,]_-$ denotes negative part and $v \to u$ in $L^2_+$. Set $ | \partial \F |_{L^2_+} (t, u) = +\infty$ if $u \not\in H^1_0$.
\end{definition}
In \S\,\ref{en-sl} we will see that
$$
	| \partial \F |_{L^2_+} (t, u) = \sup \, \big\{ \!- d\F(t, u) [z] :  z \in H^1_0 , \, | z |_{L^2_+} \le 1 \big\}    
	=  \| [ A u + f(t) ]_+ \|_{L^2} .
$$
For equivalent ways of writing the slope, with a ``singular metric'' and with a ``unilateral subdifferential'', see instead Appendix \ref{app.A} and \ref{app.B}. Actually, in the study of unilateral gradient flows we will sistematically employ the unilater slope, also in the case $f \in L^2(0,T;L^2)$ since it is technically very convenient. In particular, if $f \in AC ( 0, T; L^2)$  we have the following result. 

\begin{proposition} \label{ACsol}  If $f \in AC ( 0, T; L^2)$ then the energy $t \mapsto \F ( t , u(t))$ is absolutely continuous in $(0,T)$ and 
$u$ (the unique solution in the sense of Definition \ref{d.curve}) is also characterized by %the inequality
\begin{equation} \label{e.maxsl}
	\dot{\F} ( t, u(t) ) %=  d \F ( t , u(u) ) [ \dot{u} (t) ] + \partial_t \F ( t ,u(t)) 
\le - \tfrac12 | \dot{u} (t) |^2_{L^2_+} \,  - \tfrac12 | \partial \F |^2_{L^2_+} ( t, u (t) ) -  \langle \dot{f} (t) , u(t) \rangle .
\end{equation}
Moreover, following \cite{GianazzaSavare_RANSXLMM94}, $u \in W^{1,\infty}_{loc} ( 0 , T ; L^2 ) \cap W^{1,2}_{loc} (0,T ; H^1_0)$ is also the unique solution of 
\begin{equation}\label{e.pdeH1}
	\begin{cases}
		\partial \Phi ( \dot{u} (t) ) \ni A u (t) + f (t)  & \text{ in $H^{-1}$  for a.e.~$t \in (0,T)$,} \\
		u (0) = u_0 ,
	\end{cases}
\end{equation}
where $\Phi: H^1_0 \to [0,+\infty]$ is defined by  $\Phi (u) = \tfrac12 \| u \|_{L^2}^2 + I_{ \{u \,\ge\, 0\} } $ while $\partial \Phi (u) \subset H^{-1}$ denotes its subdifferential ($I$ is the indicator function) .
\end{proposition}

%\begin{remark}
Inequality \eqref{e.maxsl} provides, in the non-autonomous case, a notion of curve of maximal (unilateral) slope in the spirit of \cite[Definition 1.3.2]{AmbrosioGigliSavare05}. %Note that the previous inequality does not make sense if $f \in L^2(0,T;L^2)$. 
Moreover, as a consequence of Proposition \ref{ACsol}, for every $0 \le t_1 \le t_2 \le T$ the energy identities read
\begin{align} 
       \F  ( t_2 , u  (t_2) )  & = \F ( t_1 , u (t_1)  ) \, -  \, \tfrac12 \int_{t_1}^{t_2} | \dot{u} (t) |_{L^2_+}^2 + | \partial \F |_{L^2_+}^2 ( t, u (t) ) \, dt - \int_{t_1}^{t_2} \langle \dot{f} (t) , u(t) \rangle \, dt \label{e.maxslen}  \\
	& = \F ( t_1 , u (t_1)  ) \, -  \int_{t_1}^{t_2} | \dot{u} (t) |_{L^2_+} \,  | \partial \F |_{L^2_+} ( t, u (t) ) \, dt - \int_{t_1}^{t_2} \langle \dot{f} (t) , u(t) \rangle \, dt . \nonumber
\end{align}
Finally, note that the functional $\Phi$ is indeed the restriction to $H^1_0$ of the functional $\Psi$ appearing in \eqref{e.pde-akagi}. % Remark \ref{remo}.
%\end{remark}

%\begin{remark} 

%\end{remark}

\subsection[A characterization when $f$ is independent of time]{A characterization when \boldmath{$f$} is independent of time}

The case in which $f$ is independent of time has been already treated in the literature, see e.g.~\cite{GianazzaSavare_RANSXLMM94}; however, in this case we show that the monotonicity constraint on the speed $\dot{u}$ can be replaced by a fixed obstacle, as a consequence we provide a further characterization of solutions, in the spirit of the recent \cite[Remark 5.3]{AkagiEfendiev_EJAM} with a different proof.

\begin{proposition}\label{p.AE} Let $f \in L^2$ and $\F (u) = \tfrac12 a ( u , u) - \langle f , u \rangle$ (both independent of time). Given $u_0 \in H^1_0$, let $u \in H^1 (0,T;L^2) \cap L^\infty (0,T;H^1_0)$ be the unilateral gradient flow for $\F$ with initial datum $u_0$. 

Then, $u$ turns out to be the (unconstrained) $L^2$-gradient flow for the functional $\tilde{\F} (u) = \F (u) + I_+ (u -u_0)$. Moreover, $u$ is also the unique solution of the following parabolic obstacle problem: 
\begin{equation} \label{e.AE}
	\begin{cases} 
		\dot{u}(t) - A u(t) - f \ge 0 & \text{in $H^{-1}$ for a.e.~$t \in (0,T)$} \\
		( u(t) - u_0 ,  \dot{u}(t) - A u(t) - f ) = 0 & \text{for a.e.~$t \in (0,T)$} \\
		u(0) = u_0 , \quad u(t) \ge u_0 & \text{for a.e.~$t \in (0,T)$.} 
	\end{cases}
\end{equation}
\end{proposition}
%Note that in \eqref{e.AE} the obstacle is the initial condition.
As we will see in \S\ref{s.auto} the above characterization does not hold when the force $f$ depends on time.

\bigskip
\subsection{Generalizations}  \label{s.gen}

To conclude this section, let us mention that several of the above results can be extended to more general functionals, with few modifications in the definitions and in proofs. The choice of quadratic functionals is motivated by sake of simplicity and by the fact that quadratic, or separately quadratic, functionals are mostly used in applications, since they allow for easy numerical implementations. 

For instance, let $p \in (1,+\infty)$ such that $W^{1,p} \subset L^2$ (by Sobolev embedding), let $f \in L^2(0,T; L^2)$ and consider $w: \R \to [0,+\infty)$ to be $\lambda$-convex, i.e.~(see \cite[Definition 2.4.1]{AmbrosioGigliSavare05}),
$$
	w( s z_1 + (1-s) z_0 ) \le s w(z_1) + (1-s) w(z_0) - \tfrac12 \lambda s (1-s) ( z_1 -z_0)^2 ,
$$
for some $\lambda<0$ and for every $z_0, z_1 \in \R$ and $s\in (0,1)$. The double-well potential $w(z) = z^2 ( z-1)^2$, appearing in the Allen-Cahn equation, is a prototype $\lambda$-convex functions for $\lambda \le \min_z w''(z)$. Under these assumptions, we can define the stored energy $\E : W^{1,p} \to \mathbb{R}$ and the free energy $\F : [0,T] \times W^{1,p} \to \mathbb{R}$ 
$$
	\E (u) = \int_\Omega | \nabla u|^p + w(u) \, dx \,, \qquad 
	\F (t ,u) = \E (u) - \langle f(t) , u \rangle .
$$
In this case, adapting the arguments of the following sections, it is not difficult to see that the unilateral slope is still well defined and weakly lower semi-continuous in $W^{1,p}$. Thus, we can still show that for $u_0 \in W^{1,p}_0$ there exists an evolution $u \in H^1(0,T; L^2) \cap L^\infty (0,T;W^{1,p})$ which satisfies the energy identity \eqref{e.evol-eq}. However, in this weak setting, uniqueness is still open since the unilateral slope is not convex and thus the arguments of \S\ref{s.uni} do not apply.

\section{Energy and unilateral slope \label{en-sl}}

 If $u \in H^1_0$ the differential of $\F ( t, \cdot)$ restricted to $H^1_0$ is 
\begin{equation}\label{e.dF}
	d \F (t, u) [z] = a ( u ,z) - \langle f(t) , z \rangle   = - ( A u + f (t), z) \quad \text{for $z \in H^1_0$.}
\end{equation}

\begin{lemma} \label{l.slope-diff} If $u \in H^1_0$ then 
\begin{equation}\label{e.slope-diff}
| \partial \F |_{L^2_+} (t, u)  =  \sup  \big\{  \!- d\F(t, u) [z] :  z \in H^1_0 , \, | z |_{L^2_+} \le 1 \big\}   .
\end{equation}
\end{lemma}

\proof Denote $S = \sup \{  - d\F(t, u) [z] :  z \in H^1_0 , \, | z |_{L^2_+} \le 1 \} $. Since $z=0$ is an admissible variation, it is clear that $S \ge 0$.

Given $z \neq 0$ as in \eqref{e.slope-diff} let $v_s = u + s z$ for $s \ge 0$. Then, being $[r]_- \ge -r$ for $r \in \mathbb{R}$ 
$$
	| \partial \F |_{L^2_+} (t, u) \ge \limsup_{s \to 0} \frac{[\F(t, v_s) - \F(t, u)]_-}{ \| v_s - u \|_{L^2}}   \ge \limsup_{s \to 0}   \frac{\F(t, u) - \F(t, u+sz )}{s} =  - d \F (t, u) [z] .
$$
Taking the supremum on the right hand side we get $| \partial \F |_{L^2_+} (t, u) \ge S$.

Let us prove that $| \partial \F |_{L^2_+} (t, u) \le S$.  If $| \partial \F |_{L^2_+} (t, u)=0$ there is nothing  to prove. Otherwise, let $v_n \to u$ with $v_n \ge u$ s.t. $0 < | \partial \F |_{L^2_+} (t, u) = \lim_{n \to +\infty} [\F(t, v_n) - \F(t, u)]_- / \| v_n - u \|_{L^2}  $. Hence, $\F(t, v_n) < \F (t, u)$ for $n \gg 1$. % and $v_n \in H^1_0$ . 
By convexity,  $\F (t, v_n) \ge \F(t, u) + d \F (t, u) [v_n- u]$ and then, for $n \gg 1$, 
$$
	\frac{[\F(t, v_n) - \F(t, u)]_-}{\| v_n - u \|_{L^2}}  = \frac{\F(t, u) - \F(t, v_n)}{\| v_n - u \|_{L^2}} \le  \frac{- d \F (t, u) [v_n- u]}{\| v_n - u \|_{L^2}} = - d \F (t, u) [\xi_n] \le S ,
$$
where $\xi_n = (v_n - u) / \| v_n -u \|_{L^2}$ belongs to $H^1_0$, $\xi_n \ge 0$ and $\| \xi_n \|_{L^2} \le 1$. \qed

A direct consequence of the previous Lemma is the following useful result. 

\begin{corollary}\label{c.conv-meas}  Given $t \in [0,T]$ the map $ u \mapsto | \partial \F |_{L^2_+} (t, u) $ is convex in $H^1_0$. Moreover, if $u \in L^\infty(0,T; H^1_0)$ the map $t \mapsto | \partial \F |_{L^2_+} ( t , u(t))$ is measurable.
\end{corollary}

\proof Let $u, v \in H^1_0$ and $\lambda \in (0,1)$. Being $d \F (t,u) [z] = a ( u , z) - \langle f(t) , z \rangle $, 
it turns out that 
$$d \F (t, \lambda u + (1-\lambda) v) [z] = \lambda d \F (t, u) [z] + (1-\lambda) d \F (t,v) [z] . $$ 
By Lemma \ref{l.slope-diff} we have
\begin{align*}
	| \partial \F |_{L^2_+} (t, \lambda u + (1-\lambda) v) 
	& = \sup  \left\{  - d\F(t, \lambda u + (1-\lambda) v) [z] :  z \in H^1_0 , \, z \ge 0 , \, \| z \|_{L^2} \le 1 \right\}   
\end{align*}
with
\begin{align*}
	- d \F (t, \lambda u + (1-\lambda) v) [z] 
		& = - \lambda d \F (t, u) [z] - (1-\lambda) d \F (t, v) [z]  \\
		& \le \lambda | \partial \F|_{L^2_+} (t,u)  + ( 1 - \lambda) | \partial \F |_{L^2_+} (t,  v) .
\end{align*}
It follows that $| \partial \F |_{L^2_+} (t, \lambda u + (1-\lambda) v)  \le \lambda | \partial \F |_{L^2_+} (t, u)  + ( 1 - \lambda) | \partial \F |_{L^2_+} (t, v) $. 

\separe

Given $z \in H^1_0$ with $z \ge 0$ and $\| z \|_{L^2} \le 1$ the map 
$$t \mapsto - d\F ( t , u (t) ) [z] = - \int_\Omega \nabla u (t) \cdot B \, \nabla z + b \,u(t) z \, dx + \int_\Omega f (t) \, z \, dx   $$
is measurable. Taking a dense countable subset $\{ z_n \}$  of $\{ z \in H^1_0$ : $z \ge 0$ and $\| z \|_{L^2} \le 1\}$ yields
$$   | \partial \F |_{L^2_+} ( t , u(t) ) = \sup_{ n \in \mathbb{N}} \, \{  - d\F ( t , u (t) ) [z_n] \} , $$ 
where the supremum is pointwise in $(0,T)$. Measurability follows.  \qed 

\begin{corollary} \label{c.lsc-slope} %Let $f_n \in L^2(0,T; L^2)$ and let $\F_n$ be the corresponding energy. 
Given $t \in [0,T]$ let $f_n (t) \weakto f(t)$ in $L^2$ and consider the energies 
%$\F_n ( t , \cdot) : H^1_0 \to \mathbb{R}$ given by 
$$ \F_n (t , u)=  \tfrac12 a ( u , u) - \langle f_n (t) , u \rangle , \qquad \F (t , u)=  \tfrac12 a ( u , u) - \langle f (t) , u \rangle . $$
If $u_n \weakto u$ in $H^1_0$ then 
\begin{equation} \label{e.lscFn}
	\F  ( t, u )  \le \liminf_{n \to +\infty} \F_n ( t , u_n) \,, \qquad 
	| \partial \F |_{L^2_+} ( t, u )  \le \liminf_{n \to +\infty}  | \partial \F_n |_{L^2_+} ( t , u_n) \,.
\end{equation}
\end{corollary}

\proof The weak lower semicontinuity of the energy is obvious. Let $z \in H^1_0$ with $z \ge 0$ and $\| z \|_{L^2} \le 1$. If $f_n (t) \weakto f(t)$ and $u_n \weakto u$ in $H^1_0$ then we get
$$
     d\F ( t , u ) [z] = a ( u ,z ) - \langle f(t) , z \rangle  = \lim_{n \to +\infty} a ( u ,z ) - \langle f_n(t) , z \rangle 
    =  \lim_{n \to +\infty}  d\F_n ( t , u_n ) [z] .
$$
By Lemma \ref{l.slope-diff} we deduce that 
$$  - d \F ( t , u ) [z]  = \liminf_{n \to +\infty} -  d\F_n ( t , u_n ) [z]  \le \liminf_{n \to +\infty} | \partial \F_n |_{L^2_+} (t, u_n) , $$
from which \eqref{e.lscFn} follows by taking the supremum with respect to $z$.  \qed 

\separe

\begin{remark} \label{r.en} %$\F ( t, \cdot)$ is strictly convex, $H^1_0$-coercive, bounded from below and (weakly) $L^2$-lsc
The energy $\F ( t , \cdot)$ is T-monotone, i.e.~$( \xi_u  - \xi_v , [u- v]_+ ) \ge 0$ for $\xi_u \in \partial \F ( t, u )$ and $\xi_v \in \partial \F ( t, v )$, where $\partial \F (t,\cdot) \subset H^{-1}$ denotes  the subdifferential in $H^1_0$. For the details, see \cite[Lemma 2.1]{ShirakawaKimura_NA05}.
\end{remark}

\separe

To conclude, we provide in Corollary \ref{c.B1} an $L^2$ ``representation'' of the slope, which is fundamental to connect the unilateral gradient flow and the parabolic problem \eqref{e.pde}; its proof is a direct consequence of the next abstract lemmas on Radon measures.
%Stated in \cite{GianazzaSavare_RANSXLMM94}, 

\begin{lemma} \label{l.measure} Let $\zeta \in H^{-1}$. %and $ \tilde\xi \in H^1$ with $ \tilde\xi \ge 0$. 
If 
$$
	 \sup \big\{ ( \zeta,  \xi )  : \xi \in H^1_0  , \, \xi \ge 0, \,  \| \xi \|_{L^2} \le 1 \big\} < +\infty
$$
then $\zeta$ is a (locally finite) Radon measure whose positive part $\zeta_+$ belongs to $L^2$. Moreover
$$
	\sup \big\{  (\zeta , \xi )  : \xi \in H^1_0  , \, \xi \ge 0, \,  \| \xi \|_{L^2} \le 1 \big\} = \| \zeta_+ \|_{L^2} .
$$
\end{lemma}

For a proof see \cite{Negri_ACV19} or \cite{CrismaleLazzaroni_CVPDE16}.

\begin{lemma} \label{l.ineqL2}  Let $\zeta$ be as in the previous lemma
% \in H^{-1} \cap \mathcal{M}_{loc}$ such that $\zeta_+ \in L^2$. Let 
and $z \in L^2$ with $z \ge 0$ and $z \ge \zeta$ in $H^{-1}$. Then $z \ge \zeta_+$ in $L^2$.
\end{lemma}

\proof By definition $( z , \phi ) \ge ( \zeta , \phi)$ for every $\phi \in C^\infty_0$ with $\phi \ge 0$. Denote $\Omega_+$ the support of  $\zeta_+$ and let $\Omega_- = \Omega \setminus \Omega_+$. To prove the lemma it is enough to  show that $\langle z , \phi \rangle \ge \langle \zeta_+  , \phi  \rangle$  for every $\phi \in L^2$ with $\phi \ge 0$ and $\phi=0$ in $\Omega_-$. Since $\Omega_+$ is a Borel set there exists an increasing sequence $K_n$ of compact sets with $K_n \subset \subset \Omega_+$ such that $\phi \chi_n \to \phi$ in $L^2$ (where $\chi_n$ is the characteristic function of $K_n$). Let $\rho_k$ denote a smooth convolution kernel. Since $K_n \subset\subset \Omega_+$ it follows that $\phi \chi_n * \rho_k \in C^\infty_0 (\Omega_+)$ (for $k$ sufficiently large). By a diagonal argument there exists $\phi_n \in C^\infty_0$such that: $\phi_n \ge 0$, $\phi_n=0$ in $\Omega_-$ and $\phi_n \to \phi$ in $L^2$. Then 
$$
	\langle z , \phi_n \rangle = ( z , \phi_n ) \ge ( \zeta , \phi_n) =  ( \zeta_+ , \phi_n) = \langle \zeta_+  , \phi_n  \rangle .
$$
Passing to the limit concludes the proof. \qed

Invoking Lemma \ref{l.measure} together with Lemma \ref{l.slope-diff} and \eqref{e.dF} we get this  Corollary.

\begin{corollary} \label{c.B1} The following conditions are equivalent:
\begin{itemize}
\item[i)] $| \partial \F |_{L^2_+} (t, u) < +\infty$ ,
\item [ii)] 
$- d \F ( t, u) = A u + f (t) $ is a (locally finite) Radon measure with $[ Au + f (t)]_+ \in L^2$.
\end{itemize}
In this case $| \partial \F |_{L^2_+}  (t, u) = \| [ Au + f (t)]_+ \|_{L^2}$. \end{corollary}

% !TeX root = article.tex
\renewcommand{\thefootnote}{\fnsymbol{footnote}}

\section{Solutions for \boldmath{$f \in L^2(0,T;L^2)$}\label{sec.3}}

%\section{Existence by implicit Euler schemes}

In the following subsections we prove existence (and approximation) of unilateral gradient flows, in the sense of Definition \ref{d.curve}, by means of three discrete schemes, which take into account the monotonicity constraint in different ways. We remark that all these ways of representing monotonicity are currently employed in applications to phase-field fracture. We provide complete proofs, however, those parts which are very similar are not repeated. 

\subsection{Constrained incremental problem \label{3.1}}

Let $\tau_n = T /n$ and $t_{n,k} = k \tau_n$ for $k = 0, ..., n$. First of all, for every $k=0,...,n-1$ define
$$
	f_{n,k+1} = \dashint_{t_{n,k}}^{t_{n,k+1}} \hspace{-5pt} f(t) \, dt .
$$
Let $f_n \in L^2(0,T; L^2)$ given by $f_n (t) = f_{n,k+1}$ for every $t \in (t_{n,k} , t_{n,k+1}]$. Note that $f_n \to f$ in $L^2(0,T; L^2)$ and that $f_n(t) \to f(t)$ in $L^2$ for a.e.~$t\in (0,T)$.\footnote{It is enough to show that: $f_n \to f$ a.e. in $(0,T)$, $\| f_n \|_{L^2(L^2)} \le \| f \|_{L^2(L^2)}$. Then $f_n \weakto g$ in $L^2(L^2)$ and $g= f$. Get convergence of norms in $L^2(L^2)$.} Denote by $\F_n$ the corresponding energy, i.e.~$\F_n ( t , u ) = \tfrac12 a (u,u) - \langle f_n (t) , u \rangle$.  Note that $\F_n ( t , \cdot ) = \F_n ( t_{n,k+1} , \cdot)$ for every $t \in (t_{n,k} , t_{n,k+1}]$.

\separe

Define $u_{n,0} = u_0$ at time $t_{n,0}$, and then, given $u_{n,k}$ at time $t_{n,k}$, let the configuration at time $t_{n,k+1}$ be given by 
\begin{equation} \label{e.scheme1}
	u_{n,k+1} \in \argmin \big\{ \F_n ( t_{n,k+1} , u  ) + \tfrac{1}{2\tau_n} \, | u - u_{n,k} |^2_{L^2_+} :  u \in H^1_0 \big\}  .
\end{equation}
Note that a unique minimizer exists by standard arguments and that $u_{n,k+1} \ge u_{n,k}$.
%An implicit constrained Euler scheme. 

Define $u_n : [0,T] \to L^2$ and $u^\sharp_n : [0,T] \to L^2$ respectively as the piecewise affine interpolation and the piecewise constant backward (left-continuous) interpolation of the values $u_{n,k}$ in the points $t_{n,k}$. %Let us define also $f^\flat_n$ as the piecewise constant forward (right-continuous) interpolation of $f$ in the points $t_{n,k}$ and finally $t_n^\flat$ as the piecewise constant forward interpolation of $t$ in the points $t_{n,k}$. 
In this section we will prove the following poposition.

\begin{proposition} \label{p.mm} Upon extracting a subsequence (not relabelled) $u_n \weakto u$ in $H^1( 0,T; L^2)$ where $u$ is a unilateral gradient flow in the sense of Definition \ref{d.curve}. 
\end{proposition}

\begin{remark}\label{r.uniqueness} After \S\,\ref{s.uni} we will see that actually the whole sequence $u_n$ converges weakly to $u$ in $H^1(0,T;L^2)$ and that both $u_n^\sharp \to u$ and $u_n \to u$ (strongly) in $H^1_0$ pointwise in $[0,T]$. %Moreover we will prove that \eqref{e.stasera} is actually an identity. 
\end{remark}

\begin{lemma} \label{l.EL} For every $t \in (t_{n,k} , t_{n,k+1})$ it holds
\begin{equation} \label{e.due}
	 | \dot{u}_{n} (t) |^2_{L^2_+} = \| \dot{u}_{n} (t) \|^2_{L^2} =  - d \F_n ( t_{n,k+1} , u_{n,k+1})  )  [ \dot{u}_{n} (t) ]  = | \partial \F_n |_{L^2_+}^2 ( t , u^\sharp_{n} (t) )  \,.
\end{equation}
\end{lemma} 

\proof Write for simplicity $u_{n,k+1} = u_n ( t_{n,k+1})$ and $\dot{u}_{n,k+1} = (u_{n,k+1} - u_{n,k} )/ \tau_n$  instead of $\dot{u}_{n} (t)$. 
By definition $u_{n,k+1}$ is the solution (in $H^1_0$) of the variational problem
\begin{equation*} %\label{e.dvA}
	d \F_n ( t_{n,k+1}, u_{n,k+1} )  [ v - u_{n,k+1} ] + \langle \dot{u}_{n,k+1} , v - u_{n,k+1} \rangle \ge 0  \quad \text{for every $v \in H^1_0$ with $v \ge u_{n,k}$.}
\end{equation*}
Choosing $v = 2 u_{n,k+1} - u_{n,k}$ and $v = u_{n,k}$ yields $(v - u_{n,k+1} )=  \tau_n  \dot{u}_{n,k+1}$ and $(v - u_{n,k+1} )= - \tau_n  \dot{u}_{n,k+1}$, respectively; hence 
\begin{equation} \label{e.dvAA}
  d \F_n ( t_{n,k+1}, u_{n,k+1})  [ \tau_n \dot{u}_{n,k+1}]  + \tau_n \|  \dot{u}_{n,k+1} \|^2_{L^2} = 0 ,
\end{equation}
which gives the second equality in \eqref{e.due}. Moreover, the variational inequality implies that 
\begin{equation} \label{e.dPDE}
    - d \F_n ( t_{n,k+1}, u_{n,k+1} )  [ z ] \le  \langle \dot{u}_{n,k+1} , z \rangle  \quad \text{for every $z \in H^1_0$ with $z \ge 0$.}
\end{equation}
Then, by Lemma \ref{l.slope-diff} 
\begin{align}
	| \partial \F_n |_{L^2_+} ( t_{n,k+1}, u_{n,k+1}) 
	& = \sup \big\{ \!- d \F_n  ( t_{n,k+1}, u_{n,k+1} ) [z] : z \in H^1_0 ,\, z \ge 0 ,\, \| z \|_{L^2} \le 1 \big\} \nonumber \\
	& \le \sup \big\{  \langle \dot{u}_{n,k+1} , z  \rangle_{L^2} : z \in H^1_0 ,\, z \ge 0 ,\, \| z \|_{L^2} \le 1  \big\}  % \\ &
	    =  \| \dot{u}_{n,k+1} \|_{L^2} .  \label{oggi} 
\end{align}
Assume that $\dot{u}_{n,k+1} \neq 0$, otherwise \eqref{e.due} is trivial; equation \eqref{e.dvAA} provides $$d \F_n ( t_{n,k+1}, u_{n,k+1})  [ \dot{u}_{n,k+1} /  \| \dot{u}_{n,k+1} \|_{L^2} ]  + \|  \dot{u}_{n,k+1} \|_{L^2} = 0 . $$ Hence
$| \partial \F_n  |_{L^2_+} (t_{n,k}, u_{n,k+1})  =  \| \dot{u}_{n,k+1} \|_{L^2} $ and the last equality in \eqref{e.due} is proved. \qed

\begin{lemma} \label{p.enest} For every $0 \le k \le n-1$, the following energy estimate holds
\begin{align}
	\E ( u_n ( t_{n,k+1} ) )   & = \E ( u_n ( t_{n,k} ) ) -  \int_{t_{n,k}}^{t_{n,k+1}} \tfrac12 | \dot{u}_{n} (t) |^2_{L^2_+}  + \tfrac12 | \partial \F_n |^2_{L^2_+}  ( t , u^\sharp_{n} (t) )  \, dt \ + \nonumber  \\ & 
+ \int_{t_{n,k}}^{t_{n,k+1}} \langle f_n (t) , \dot{u}_n (t) \rangle \, dt %- \tau_n \int_{t_{n,k}}^{t_{n,k+1}}  \| \dot{u}_n  (t) \|^2_{a} \, dt  
\,. \label{e.enest-2}
\end{align}
\end{lemma}

\proof Write, $u_{n,k+1} = u_n ( t_{n,k+1})$ etc. Using \eqref{e.due} and being $\F_n (t_{n,k+1}, \cdot)$ convex we get 
\begin{align*}
	\F_n ( t_{n,k+1} , u_{n,k} )  & \ge  \F_n (  t_{n,k+1} , u_{n,k+1} ) + d \F_n ( t_{n,k+1} ,  u_{n,k+1} ) [ u_{n,k} - u_{n,k+1}]  \phantom{1_{L^2_+}} \\ % + \tfrac12 \|   \nabla u_{n,k}  - \nabla  u_{n,k+1} \|^2_{L^2} \\
	& = \F_n ( t_{n,k+1} ,  u_{n,k+1} ) - \tau_n \,d \F_n ( t_{n,k+1} ,  u_{n,k+1} ) [ \dot{u}_{n}] \\ % + \tfrac12 \|  u_{n,k} - u_{n,k+1} \|^2_{a} \\
	& = \F_n ( t_{n,k+1} ,  u_{n,k+1} ) + \tau_n  \big(  \tfrac12 | \dot{u}_{n} |^2_{L^2_+} +  \tfrac12 | \partial \F_n |^2_{L^2_+} ( t_{n,k+1}  , u_{n,k+1})  \big) .% + 	\tfrac12 \|  u_{n,k} - u_{n,k+1} \|^2_{a} 
\end{align*}
Writing
\begin{align*}
	  \F_n ( t_{n,k+1} , u_{n,k} )  = \E (  u_{n,k} ) - \langle f_{n,k} \,,  u_{n, k} \rangle \,,
	\qquad
    \F_n ( t_{n,k+1} , u_{n,k+1} )  = \E (  u_{n,k+1} ) - \langle f_{n,k+1} \,,  u_{n, k+1} \rangle \,,
\end{align*}
we get 
\begin{align*}
  	\E (  u_{n,k}  ) - \langle f_{n,k} \,,  u_{n, k} \rangle \ge & \ \E (  u_{n,k+1} ) - \langle f_{n,k} \,,  u_{n, k+1} \rangle \\  & + \tau_n  \big(  \tfrac12 | \dot{u}_{n} |^2_{L^2_+} +  \tfrac12 | \partial \F_n |^2_{L^2_+} ( t_{n,k+1}  , u_{n,k+1})  \big) . % + \tfrac12 \tau_n^2  \|  \dot{u}_{n} \|^2_{a} \,.
\end{align*}
Using the interpolant $u_n$, $u_n^\sharp$, and $f_n$ we get  \eqref{e.enest-2}. \qed

\noindent {\bf Proof of Proposition \ref{p.mm}.} Using \eqref{e.enest-2} for $0 \le k \le n-1$ and \eqref{e.due} we get 
\begin{align*}
	\E (  u_n ( T ) ) & \le \E ( u_0 ) - \int_{0}^{T} \|\dot{u}_{n} (t) \|^2_{L^2}  + \int_{0}^{T} \langle f_n (t) ,\dot{u}_n (t) \rangle \, dt .  \label{e.enest-2}
\end{align*}
Hence
$ %\begin{align*}
       \|\dot{u}_n \|^2_{L^2(0,T; L^2)} % = \int_0^T \|\dot{u}_n (t) \|^2_{L^2} \, dt  
%                & \le \tfrac12 \int_\Omega | \nabla u_0 |^2 \, dx   + \int_0^T \langle f^\flat _n (t) ,\dot{u}_n (t) \rangle \, dt  \\
	      \le C + \| f_n \|_{L^2(0,T; L^2)} \, \|\dot{u}_n \|_{L^2(0,T; L^2)} 
$. %\end{align*}
Being $f_n$ bounded in $L^2(0,T; L^2)$ the sequence $u_n$ turns out to be bounded in $H^1 (0,T; L^2)$ and thus, up to subsequences (not relabelled), $u_n \weakto u$ in $H^1 (0,T ; L^2)$. We will identify $u$ with its absolutely continuous representantive, i.e. 
$$
	u(t) = u_0 + \int_0^t \dot{u} (s) \, ds .
$$
By weak convergence it is easy to check that the limit $u$ is monotone in time. For every  $0 \le k \le n-1$  from \eqref{e.enest-2} we get 
\begin{align*}
	\E ( u_n ( t_{n,k+1} ) ) & \le \E (  u_0 ) + \int_{0}^{t_{n,k+1}} \langle f_n (t) ,\dot{u}_n (t) \rangle \, dt .  
\end{align*}
By coercivity of the stored energy, we deduce that $u_n$ is bounded in $L^\infty (0,T ; H^1_0)$ and thus $u_n^\sharp$ is bounded in $L^\infty (0,T ; H^1_0)$ as well. It follows that $u^\sharp_n (t) \weakto u(t)$ in $H^1_0$ for every $t \in [0,T]$; indeed, given $t \in (0,T)$, let $k_n$ s.t.~$k_n < t \le k_n+1$, then $u_n^\sharp ( t) = u_n ( t_{n,k_n+1} )$ converges weakly to $u(t)$ in $L^2$, because $t_{n,k_n+1} \to t$ and $u_n \weakto u$ in $H^1 (0,T ; L^2)$, and $u_n^\sharp ( t)$ is bounded in $H^1_0$.

%By  \eqref{e.scheme1} $\F_n ( u_n  (t_{n,k}) ) \le \F_n ( u_0)$ for every $k$. Hence, by coercivity of $\F_n$, $u_n \in L^\infty ( 0,T; H^1_0)$where $C$ is a lower bound for $\F_n$ in $H^1_0$  (cf.~Proposition \ref{p.en}).
% and $u_n \weakstarto u$ in $L^\infty(0,T ;H^1_0)$.

Given $t^* \in (0,T]$ let $k_n$ such that $t_{n,k_n} <  t^* \le t_{n,k_n+1}$. Using \eqref{e.enest-2} we get 
\begin{equation} \label{e.123}
	\E (  u_n ( t_{n,k_n+1} ) ) + \int_{0}^{t_{n,k_n+1}} \hspace{-5pt} \tfrac12 |\dot{u}_{n} (t) |^2_{L^2_+}  +  \tfrac12 | \partial \F_n |^2_{L^2_+}  ( t , u^\sharp_{n} (t) ) - \langle f_n(t) , \dot{u}_n (t)  \rangle  \, dt  \le  \E ( u_0 ) \,.
\end{equation}
Clearly $t_{n,k_n+1} \to t^*$ and  hence the sequence $u_{n,k_n+1} = u_n^\sharp ( t^* )$ converges weakly to $u (t^*)$ in $H^1_0$. Hence
$$
       \E ( u (t^*) )  \le \liminf_{n \to +\infty} \E ( u_n (  t_{n,k_n+1} ) ) .
$$ 
By weak convergence in $H^1(0,T; L^2)$ we get 
\begin{equation} \label{e.alprim}
	\int_0^{t^*} \|\dot{u} (t) \|^2_{L^2} \, dt \le \liminf_{n \to +\infty} \int_0^{t^*} \|\dot{u}_n (t) \|^2_{L^2} \, dt \le  \liminf_{n \to +\infty} \int_0^{t_{n,k_n+1}} \|\dot{u}_n (t) \|^2_{L^2} \, dt .
\end{equation}
All the above $\| \cdot \|_{L^2}$ can be replaced with $| \cdot |_{L^2_+}$ since $\dot{u}_n$ and $\dot{u}$ are non-negative. Since $f_n(t) \to f(t)$ in $L^2$ and $u^\sharp_n(t) \weakto u(t)$ in $H^1$ for a.e.~$t\in (0,T)$, by \eqref{e.lscFn} we get 
$$
             | \partial \F |_{L^2_+} ( t, u (t) )  \le \liminf_{n \to +\infty}  | \partial \F_n |_{L^2_+} ( t , u_n^\sharp (t)) .
$$
Them, by Fatou's Lemma 
\begin{equation} \label{e.alsecond}
     \int_0^{t^*}  | \partial \F |^2_{L^2_+} ( t, u (t) ) \, dt \le \liminf_{n \to +\infty} \int_0^{t^*} | \partial\F_n |^2_{L^2_+} ( t , u_n^\sharp (t)) \, dt \le \liminf_{n \to +\infty} \int_{0}^{t_{n,k_n+1}}   | \partial\F_n |^2_{L^2_+}  ( t , u^\sharp_{n} (t) ) \, dt .
\end{equation}
For last term in the left hand side of \eqref{e.123} we easily have
\begin{equation} \label{e.alters}
   \lim_{n \to +\infty}  \int_0^{t^*} \langle f_n(t) , \dot{u}_n (t) \rangle \, dt =   \int_0^{t^*} \langle f(t) , \dot{u} (t) \rangle \, dt .
\end{equation}
Taking the liminf in \eqref{e.123} gives the thesis. \qed

\subsection{Unconstrained incremental problem with a posteriori truncation \label{s.3-2}}

In this \S\,we consider an alternative discrete scheme, numerically more convenient than \eqref{e.scheme1}, employed in \cite{TakaishiKimura_K09}. Let $\tau_n$ and $t_{n,k}$ be as in the previous subsection. Define $u_{n,0} = u_0$ at time $t_{n,0}$ and then by induction let 
\begin{gather} \label{e.scheme2}
	\begin{cases} \tilde u_{n,k+1}  \in \argmin \big\{ \F_n ( t_{n,k+1} , u  ) + \tfrac{1}{2\tau_n} \, \| u - u_{n,k} \|^2_{L^2} :  u \in H^1_0 \big\}  \\
                                    u_{n,k+1} = \max \big\{  \tilde u_{n,k+1} , u_{n,k}  \big\} . 
            \end{cases}
\end{gather}
Note that the first minimization is unconstrained, the constraint is taken into account a posteriori, simply by truncation. In this way $u_{n,k+1} - u_{n,k} = [ \tilde u_{n,k+1} - u_{n,k}]_+$. 
  
As in the previous subsection we define $u_n : [0,T] \to L^2$ as the piecewise affine interpolation and $u^\sharp_n$ as the piecewise constant backward (left-continuous) interpolation of $u_{n,k}$ in the points $t_{n,k}$. Moreover, we define $\tilde {u}_n : [0,T] \to L^2$ as the piecewise constant backwards (left-continuous) interpolation of $\tilde{u}_{n,k}$, again in the points $t_{n,k}$. 

\begin{proposition} \label{p.mmbis} Upon extracting a subsequence (not relabelled) $u_n \weakto u$ in $H^1( 0,T; L^2)$ where $u$ is a unilateral gradient flow in the sense of Definition \ref{d.curve}. 
\end{proposition}

\begin{lemma} \label{l.ELbis} For every index $0 \le k \le n-1$ 
\begin{align} 
	 | \dot{u}_{n} (t) |_{L^2_+}  = \| \dot{u}_{n} (t) \|_{L^2}  & =  | \partial \F_n |_{L^2_+} ( t , \tilde{u}_{n} (t) )  \quad \text{for $t \in (t_{n,k} , t_{n,k+1})$} .%= - \tau d \F_n ( u_{n} ( t_{n,k+1})  )  [ \dot{u}_{n} (t) ] 
         \label{e.tre} 
\end{align}
Moreover, 
\begin{align}
         d \F_n ( t_{n,k+1} , u_{n,k+1}) [u_{n,k+1} - u_{n,k}] & = d \F_n ( t_{n,k+1}, \tilde u_{n,k+1}) [u_{n,k+1} - u_{n,k}] \label{e.tremezzo} \\ 
     		& = - | \partial \F_n |_{L^2_+} ( t_{n,k+1}, \tilde{u}_{n,k+1} ) \| u_{n,k+1} - u_{n,k}  \|_{L^2}  \label{e.quattro}.
\end{align}
\end{lemma} 

\proof By minimality
\begin{equation} \label{e.dFutilde}
   d\F_n (t_{n,k+1}, \tilde u_{n,k+1}) [z] + \tfrac{1}{\tau_n} \langle   \tilde u_{n,k+1} - u_{n,k} , z \rangle = 0 \quad \text{for every $z \in H^1_0$.} 
\end{equation}
Hence, by \eqref{e.slope-diff} 
\begin{align*}
     | \partial \F_n |_{L^2_+} (t_{n,k+1}, \tilde u_{n,k+1})  
                & \ =  \sup  \left\{  - d\F_n( t_{n,k+1}, \tilde u_{n,k+1}) [z] :  z \in H^1_0 , \, z \ge 0 , \, \| z \|_{L^2} \le 1 \right\}    \\
                & \ =  \max  \left\{ \tfrac{1}{\tau_n} \langle   \tilde u_{n,k+1} - u_{n,k} , z \rangle :  z \in L^2 , \, z \ge 0 , \, \| z \|_{L^2} \le 1 \right\}   \\
                & \ =  \tfrac{1}{\tau_n} \langle  \tilde u_{n,k+1} - u_{n,k} , \frac{ [ \tilde u_{n,k+1} - u_{n,k}]_+ }{\| [ \tilde u_{n,k+1} - u_{n,k}]_+ \|_{L^2} } \rangle  \\
	                & \ = \tfrac{1}{\tau_n} \| [ \tilde u_{n,k+1} - u_{n,k}]_+ \|_{L^2} = \tfrac{1}{\tau}  | u_{n,k+1} - u_{n,k} |_{L^2_+} ,
\end{align*}
which gives \eqref{e.tre}. In particular, remembering that $u_{n,k+1} - u_{n,k} = [ \tilde u_{n,k+1} - u_{n,k}]_+$ we get
\begin{equation} \label{e.aus}
     - d\F_n( t_{n,k+1} , \tilde u_{n,k+1}) [ u_{n,k+1} - u_{n,k} ] =  | \partial \F_n |_{L^2_+} ( t_{n,k+1} , \tilde u_{n,k+1})  \| u_{n,k+1} - u_{n,k} \|_{L^2} ,
\end{equation}
which gives \eqref{e.quattro}.

\separe

It remains to prove \eqref{e.tremezzo}. Define $\Omega_- = \{  \tilde u_{n,k+1} \le u_{n,k}  \}$ and $\Omega_+ = \{  \tilde u_{n,k+1} > u_{n,k}  \}$.
We claim that
\begin{equation} \label{e.dFuz}
     d\F_n ( t_{n,k+1}, u_{n,k+1}) [z] + \tfrac{1}{\tau_n} \langle u_{n,k+1} - u_{n,k} , z \rangle = 0 \quad \text{for every $z \in H^1_0$ with $z=0$ in $\Omega_-$.} 
\end{equation}
Since $z=0$ in $\Omega_-$ 
 $$\int_{\Omega_-} \nabla u_{n,k+1} \cdot B \, \nabla z + b \, u_{n,k+1} z - f_{n,k+1} z + \tfrac{1}{\tau_n} (u_{n,k+1} - u_{n,k}) z \, dx = 0 . $$  
As $\tilde u_{n,k+1} = u_{n,k+1} $ in $\Omega_+$ and  $z=0$ in $\Omega_-$ 
\begin{align*}
  \int_{\Omega_+} \nabla u_{n,k+1} \cdot B \, \nabla z & + b \, u_{n,k+1} z - f_{n,k+1} z + \tfrac{1}{\tau_n} (u_{n,k+1} - u_{n,k}) z \, dx  \\ & =   \int_{\Omega_+} \nabla \tilde u_{n,k+1} \cdot B \, \nabla z + b \, \tilde{u}_{n,k+1} z - f_{n,k+1} z + \tfrac{1}{\tau_n} ( \tilde u_{n,k+1} - u_{n,k}) z \, dx  \\
 & =  \int_{\Omega\phantom{+}} \nabla \tilde u_{n,k+1} \cdot B \, \nabla z + b \, \tilde{u}_{n,k+1} z - f_{n,k+1} z + \tfrac{1}{\tau_n} ( \tilde u_{n,k+1} - u_{n,k}) z \, dx  = 0 ,
\end{align*}
where the last equality follows by minimality. Joining the integrals on $\Omega^\pm$ proves \eqref{e.dFuz}.

\separe

Using \eqref{e.dFutilde} and  \eqref{e.dFuz} with $z = u_{n,k+1} - u_{n,k}$ (note that $z=0$ in $\Omega_-$)
\begin{align*}
	& d \F_n ( t_{n,k+1} , \tilde u_{n,k+1}) [u_{n,k+1} - u_{n,k}] + \tfrac{1}{\tau_n} \langle \tilde u_{n,k+1} - u_{n,k} , u_{n,k+1} - u_{n,k} \rangle = 0 \\
           &   d \F_n ( t_{n,k+1} , u_{n,k+1}) [u_{n,k+1} - u_{n,k}] + \tfrac{1}{\tau_n} \langle u_{n,k+1} - u_{n,k} , u_{n,k+1} - u_{n,k} \rangle = 0 \,.
\end{align*}
Let us see that $\langle \tilde u_{n,k+1} - u_{n,k} , u_{n,k+1} - u_{n,k} \rangle  = \langle  u_{n,k+1} - u_{n,k} , u_{n,k+1} - u_{n,k} \rangle $, indeed their difference reads $\langle \tilde u_{n,k+1} - u_{n,k+1} , u_{n,k+1} - u_{n,k} \rangle$ which vanishes because $\tilde u_{n,k+1} - u_{n,k+1} = 0$ in $\Omega_+$ and $u_{n,k+1} - u_{n,k} = 0$ in $\Omega_-$. Hence   $d \F_n (t_{n,k+1}, \tilde u_{n,k+1}) [u_{n,k+1} - u_{n,k}] = d \F_n ( t_{n,k+1}, u_{n,k+1}) [u_{n,k+1} - u_{n,k}] $. \qed

\separe

Next lemma follows closely the corresponding one in the previous subsection. 

\begin{lemma} \label{p.enestbis} For every $0 \le k \le n-1$, the following energy estimate holds
%\begin{align} 
% 	\F_n (t_{n,k+1} , u_{n} (t_{n,k+1}) ) & \le \F_n ( t_{n,k+1} , u_{n} (t_{n,k+1}) )  -  \tfrac12 \int_{t_{n,k+1}}^{t_{n,k+1}} \| \dot{u}_{n} (t) \|^2_{L^2}  +  | \partial_+ \F_n |^2  ( t_{n}^\flat(t)  , \tilde{u}_{n} (t) )  \, dt \ + \nonumber  \\ & - \int_{t_{n,k+1}}^{t_{n,k+1}} \langle f'(t) , u^\sharp_n(t) \rangle \, dt \label{e.enestbis}  
%\end{align}
\begin{align}
	\E (  u_n ( t_{n,k+1} ) ) & \le \, \E ( u_n ( t_{n,k} ) ) -  \int_{t_{n,k}}^{t_{n,k+1}} \tfrac12  | \dot{u}_{n} (t) |^2_{L^2_+} + \tfrac12  | \partial \F_n |^2_{L^2_+}  ( t , \tilde{u}_{n} (t) )  \, dt \ + \nonumber  \\ & 
+ \int_{t_{n,k+1}}^{t_{n,k+1}} \langle f_n (t) , \dot{u}_n (t) \rangle \, dt .  \label{e.enest-2bis}
\end{align}
\end{lemma}

\proof Write, $u_{n,k+1} = u_n ( t_{n,k+1})$ etc. By convexity and by \eqref{e.quattro} and \eqref{e.tre}
\begin{align*}
	\F_n ( t_{n,k+1} , u_{n,k} )  & \ge \F_n (  t_{n,k+1} , u_{n,k+1} ) + d \F_n ( t_{n,k+1} ,  u_{n,k+1} ) [ u_{n,k} - u_{n,k+1}] \phantom{1_{L^2_+}} \\
	& \ge \F_n ( t_{n,k+1} ,  u_{n,k+1} ) + \tau_n | \partial \F_n |_{L^2_+}  ( t_{n,k+1} ,  \tilde{u}_{n,k+1} ) \| \dot{u}_{n} \|_{L^2} \\
	& \ge \F_n ( t_{n,k+1} ,  u_{n,k+1} ) + \tau_n  \big(  \tfrac12 \| \dot{u}_{n}  \|^2_{L^2} +  \tfrac12 | \partial \F_n |^2_{L^2_+} ( t_{n,k+1}  , \tilde{u}_{n,k+1})  \big) .
\end{align*}
%Hence
%$$
%     \F_n ( t_{n,k+1} , u_{n,k} )   \ge \F_n ( t_{n,k+1} , u_{n, k+1} )  + \tfrac12 \int_{t_{n,k+1}}^{t_{n,k+1}} \| \dot{u}_{n} (t) \|^2_{L^2}  +  | \partial_+ \F_n |^2  ( t_{n}^\flat(t)  , \tilde{u}_{n} (t) )  \, dt
%$$
Following line by line the proof of Lemma \ref{p.enest} provides \eqref{e.enest-2bis}.  \qed

Note that formally (since the sequences do not coincide) the only difference between \eqref{e.enest-2} and \eqref{e.enest-2bis} is the slope: in the former it is evaluated in $(t , u^\sharp_n (t))$ while in the second in $(t , \tilde{u}_n (t))$.

\bigskip
\noindent {\bf Proof of Proposition \ref{p.mmbis}.} Following line by line the first step in the proof of Proposition \ref{p.mm} we get that the sequence $u_n$ is bounded in $L^\infty(0,T; H^1_0)$ and in $H^1(0,T; L^2)$. Hence $u_n \weakto u$ in $H^1( 0,T; L^2)$, upon extracting a (non-relabelled) a subsequence. 

We claim that $\tilde{u}_n (t) \weakto u(t)$ in $H^1_0$ for a.e.~$t\in (0,T)$. Fix $t \in (0,T)$ s.t. $f_{n} (t) \to f(t)$ in $L^2$. For every $n \in \mathbb{N}$ let $k_n$ (depending on $t$) s.t. $t_{n,k_n} < t \le t_{n,k_n+1}$. Note that $f_{n,k+1} = f_n(t)  \to f(t)$ in $L^2$. By minimality we can write that
$$
	\F_n ( t_{n,k+1} ,  \tilde{u}_{n,k+1})  + \tfrac1{2\tau_n} \| \tilde{u}_{n,k+1} - u_{n,k} \|^2_{L^2} \le \F_n ( t_{n,k+1} ,  u_{n,k}) \le c \, \| u_{n,k} \|_{H^1_0}^2 + \| f_{n,k+1} \|_{L^2}  \| u_{n,k} \|_{L^2} \,,
$$
where we used the continuity of the bi-linear form. Since $u_n$ is bounded in $L^\infty(0,T;H^1_0)$ the right hand side is bounded uniformly w.r.t.~$n$ and $k$; thus there exists $C>0$ s.t.
\begin{gather*}
	\| \tilde{u}_{n,k+1} - u_{n,k} \|^2_{L^2} \le C \tau_n , \\
	c' \| \tilde{u}_{n,k+1} \|^2_{H^1_0} - \| f_{n,k+1} \|_{L^2}  \| \tilde{u}_{n,k+1} \|_{H^1_0} \le \F_n ( t_{n,k+1} ,  \tilde{u}_{n,k+1}) \le C ,
\end{gather*}
where in the second line we used coercivity. As $f_{n,k+1} \to f(t)$ in $L^2$, simple algebraic estimate yields $\| \tilde{u}_{n,k+1} \|_{H^1_0} \le C'$.

Since $u_{n,k} = u_n ( t_{n,k}) \to u(t)$ in $L^2$ it follows that $\tilde{u}_n (t) = \tilde{u}_{n,k+1}  \to u(t)$ in $L^2$. Being $\tilde{u}_{n,k+1}$ bounded uniformly in $H^1_0$ we get $\tilde{u}_n (t) \weakto u(t)$ in $H^1_0$.

To conclude the proof it is enough to argue as in the proof of Proposition \ref{p.mm} \qed

%\begin{remark} It is important to note that the convergence of the scheme ... follows easily from the convergence of the scheme ... because our evolutions are of class $H^1 (0,T;L^2)$. In this way the only difference is estimate ...  \end{remark}

\subsection{Unconstrained incremental problem with penalty \label{s.pen}}

Let $\alpha: (0,+\infty) \to [1, +\infty)$ be monotone non-increasing with $\lim_{\tau \to 0^+} \alpha(\tau) = +\infty$. For $\tau >0$ and $v \in L^2 (\Omega)$ let us denote 
$$
	\psi_\tau ( v ) = 
		\begin{cases}
		v^2 & \text{if $v \ge 0$} \\
		\alpha(\tau)v^2 & \text{if $v < 0$}
		\end{cases}
	, \qquad 
	\Psi_\tau (v) = \int_\Omega \psi_\tau (v) \, dx ,
	\qquad | v |_{L^2_\tau} = \Psi^{1/2}_\tau (v) \,.
$$
Clearly, when $\tau$ is small $\alpha (\tau)$ is large and thus $|v  |_{L^2_\tau}$ penalizes $v_-$. Note that $\Psi_\tau$ can be equivalently seen as the Yosida regularization of the indicator function of the set $\{ v \in L^2 : v \ge 0\}$.

Moreover, we have $ | v |_{L^2_\tau} \ge \| v \|_{L^2} \ge | v |_{L^2_\tau} / \alpha (\tau)$, hence $v_n \to v$ in $L^2_\tau$ actually means that $v_n \to v$ in $L^2$. Clearly
$$
	d \Psi_{\tau} ( v ) [ z ] = \int_\Omega \psi'_\tau (v) z \, dx = 2 \int_{ \{ v \,\ge\, 0 \} } \! \! vz  \, dx  + 2 \, \alpha(\tau) \int_{ \{ v \,< \, 0 \} } \! \! vz \, dx .
$$
Before proceeding, let us prove this lemma. 

\begin{lemma} \label{l.Psi}
For every $\tau >0$ and $v \in L^2$ it holds
\begin{equation} \label{e.wlaf}
	\sup \big\{ d \Psi_{{\tau}} ( v ) [ z ] : z \in H^1_0 ,\, | z |_{L^2_{\tau}} \le 1  \big\} 
		= \max  \big\{ d \Psi_{{\tau}} ( v ) [ z ] : z \in L^2 ,\, | z |_{L^2_{\tau}} \le 1  \big\} 
 = 2 | v |_{L^2_{\tau}} .
\end{equation} 
\end{lemma}

\proof We introduce the set 
$\hat{V} = \{ z \in L^2 : z=0 \text{ if } v=0 , \,z\ge 0 \text{ if } v>0, \, z<0 \text{ if } v<0  \}$ and we will prove that 
\begin{align*} 
	\sup \big\{ d \Psi_{{\tau}} ( v ) [ z ] : z \in H^1_0 ,\, | z |_{L^2_{\tau}} \le 1  \big\}  
		& = \, \sup \big\{ d \Psi_{{\tau}} ( v ) [ z ] : z \in L^2 ,\, | z |_{L^2_{\tau}} \le 1  \big\}   \\
		& = \, \sup \big\{ d \Psi_{{\tau}} ( v ) [ z ] : z \in L^2 ,\, | z |_{L^2_{\tau}} \le 1 ,\, 
						   z \in \hat{V} \big\}   \\
		& = 2 | v |_{L^2_{\tau}} .
\end{align*} 
It is enough to consider $| v |_{L^2_{\tau}}\neq 0$, otherwise there is nothig to prove. 
The first identity follows by density of  $H^1_0$ in $L^2$ and by continuity of $ | \cdot |_{L^2_\tau}$ and $d \Psi_\tau ( v) [ \cdot]$. Let us check the second identity. Clearly 
$$
 \sup \big\{ d \Psi_{{\tau}} ( v ) [ z ] : z \in L^2 ,\, | z |_{L^2_{\tau}} \le 1  \big\} \ge \sup \big\{ d \Psi_{{\tau}} ( v ) [ z ] : z \in L^2 ,\, | z |_{L^2_{\tau}} \le 1 ,\, 
						   z \in \hat{V} \big\} .
$$
In order to prove the opposite inequality, given $z \in L^2$
%Note that in \eqref{e.wlaf} it is not restrictive to consider that $\mathrm{sign} (z) = \mathrm{sign} (v)$, indeed, consider the function
let $\tilde{z} \in L^2$ be defined by 
$$
	\tilde{z} = \begin{cases} 
		z &  \text{if $v > 0$ and $z \ge 0$} \\ 
		0 &  \text{if $v > 0$ and $z < 0$} \\
		0 &  \text{if $v=0$} \\
		z &  \text{if $v < 0$ and $z<0$} \\
		0 &  \text{if $v < 0$ and $z \ge 0$.}
\end{cases}
$$ 
Then, $\tilde{z}\in \hat{V}$,  $| \tilde{z} |_{L^2_{\tau}} \le | z |_{L^2_{\tau}} \le 1 $ and $ d \Psi_{\tau} ( v ) [ z ] \le d \Psi_{\tau} ( v ) [ \tilde{z} ]$ because 
$$
	\int_{ \{ v \,>\, 0 \} } \! \! vz  \, dx \le \int_{ \{ v \,>\, 0 \} } \! \! v \tilde{z}  \, dx 
	\quad \text{and} \quad 
     \int_{ \{ v \,< \, 0 \} } \! \! vz \, dx \le \int_{ \{ v \,< \, 0 \} } \! \! v \tilde{z} \, dx .
$$
To prove the last equality, given $\tau$ and $v$, let us introduce the Hilbert space $L^2_\lambda$ where $\lambda$ is the measure $\lambda = \mathcal{L}^d_{| \{  v \, \ge \,  0 \} } + \alpha(\tau) \mathcal{L}^d_{| \{  v \, < \,  0 \} }$ (being $\mathcal{L}^d$ the $d$-dimensional Lebesgue measure). Then we write 
\begin{align*}
      d \Psi_{\tau} ( v ) [ z ] & = 2 \int_{ \{ v \,\ge\, 0 \} } \! \! vz  \, dx  + 2 \alpha(\tau) \int_{ \{ v \,< \, 0 \} } \! \! vz \, dx = 2 \int_\Omega vz \, d \lambda = 2 \langle v , z \rangle_{L^2_\lambda} \le 2 \| v  \|_{L^2_\lambda} \, \| z  \|_{L^2_\lambda} \,.
\end{align*}
Note that $\| v \|_{L^2_\lambda} = | v |_{L^2_{\tau}}$ while in general $\| z \|_{L^2_\lambda} \neq | z |_{L^2_{\tau}}$; however, $\| z \|_{L^2_\lambda} = | z |_{L^2_{\tau}}$ if $z \in \hat{V}$, and thus the previous inequality yields
$$
	\sup \big\{  d \Psi_{\tau} ( v ) [ z ] : z \in L^2 ,\, | z |_{L^2_{\tau}} \le 1 ,\, z \in \hat{V} \big\}  \le 2 | v |_{L^2_{\tau}} .
$$
Clearly, $z = v / | v |_{L^2_{\tau}}$ is the maximizer, which gives \eqref{e.wlaf}. \qed 

Finally, let us introduce the slope 
$$
	| \partial \F |_{L^2_\tau} (t,u) = \limsup_{ v \to u} \frac{[ \F (t ,v) - \F ( t, u) ]_-}{| v - u |_{L^2_\tau}} = \sup \big\{ \!- \partial \F ( t ,u ) [ \phi] : | \phi |_{L^2_\tau} \le 1  \big\} ,
$$
where $v \to u$ in $L^2$ while the second identity follows by the same arguments as in the proof of Lemma \ref{l.slope-diff}. 
We remark that, if $| \phi |_{L^2_+} \le 1$, then $\phi \ge 0$ and $\| \phi \|_{L^2} \le 1$, thus $ | \phi |_{L^2_\tau} = \| \phi \|_{L^2} \le 1$, hence  $| \partial \F |_{L^2_\tau} (t,u) \ge| \partial \F |_{L^2_+} (t,u)$. 
%As a consequence, if $\tau_m \to 0$, $t_m \to t$ and $u_m \weakto u$ in $H^1_0(\Omega)$ then by Corollary  \label{c.lsc-slope}
%\begin{equation} \label{e.lsctau3}
%	| \partial \F |_{L^2_+} ( t,u) \le \liminf_{m \to +\infty} | \partial \F |_{L^2_+} ( t_m , u_m)  \le \liminf_{m \to +\infty} | \partial \F |_{L^2_{\tau_m}} ( t_m , u_m) .
%\end{equation}

\separe

\bigskip
Let $t_{n,k}$ and $\F_n$ be as in the previous sections. Define $u_{n,0} = u_0$ and by induction
\begin{gather} \label{e.scheme3}
	u_{n,k+1} \in \argmin \big\{ \F_n ( t_{n,k+1} , u  ) + \tfrac{1}{2\tau_n} \, | u - u_{n,k} |^2_{L^2_{\tau_n}} :  u \in H^1_0 \big\}  .
\end{gather}
Let $u_n : [0,T] \to L^2 (\Omega)$ and $u^\sharp_n : [0,T] \to L^2 (\Omega)$ denote respectively the piecewise affine and the piecewise constant interpolation of $u_{n,k}$ in $t_{n,k}$, as in the previous sections.  %Let us define also $f^\flat_n$ as the piecewise constant forward (right-continuous) interpolation of $f$ in the points $t_{n,k}$ and finally $t_n^\flat$ as the piecewise constant forward interpolation of $t$ in the points $t_{n,k}$. 

\begin{proposition} \label{p.mm3} Upon extracting a subsequence (non relabelled) $u_n \weakto u$ in $H^1( 0,T; L^2)$ where $u$ is unilateral gradient flow, in the sense of Definition \ref{d.curve}. \end{proposition}

\begin{lemma} \label{l.EL3} For every $t \in (t_{n,k} , t_{n,k+1})$
\begin{equation} \label{e.due3}
	 | \dot{u}_{n} (t) |^2_{L^2_{\tau_n}} =  - d \F_n ( t_{n,k+1} , u_{n,k+1})  )  [ \dot{u}_{n} (t) ]  = | \partial \F_n |_{L^2_{\tau_n}}^2 ( t , u^\sharp_{n} (t) )  \,.
\end{equation}
\end{lemma} 

\proof %Write for simplicity $u_{n,k+1} = u_n ( t_{n,k+1})$ and $\dot{u}_{n} = (u_{n,k+1} - u_{n,k} )/ \tau_n$  instead of $\dot{u}_{n} (t)$. 
By minimality $u_{n,k+1}$ satisfies the Euler-Lagrange equation
\begin{equation*} %\label{e.dvA}
	d \F_n ( t_{n,k+1}, u_{n,k+1} )  [ z ] + \tfrac{1}{2} d \Psi_{{\tau_n}} ( \dot{u}_{n} ) [ z ] = 0  \quad \text{for every $z \in H^1_0$.}
\end{equation*}
Hence,
\begin{align}
	| \partial \F_n |_{L^2_{\tau_n}} ( t_{n,k+1}, u_{n,k+1}) 
	& = \sup \big\{ \!- d \F_n  ( t_{n,k+1}, u_{n,k+1} ) [z] : z \in H^1_0 ,\, | z |_{L^2_{\tau_n}} \le 1 \big\} \nonumber \\
	& = \sup \big\{ \tfrac12 d \Psi_{{\tau_n}} ( \dot{u}_{n} ) [ z ] : z \in H^1_0 ,\, | z |_{L^2_{\tau_n}} \le 1  \big\} \nonumber
\end{align}
By Lemma \ref{l.Psi}, choosing $z = \dot{u}_n$ we obtain
$$
   - d \F_n  ( t_{n,k+1}, u_{n,k+1} ) [\dot{u}_n] = \tfrac12 d \Psi_{\tau_n} ( \dot{u}_n) [\dot{u}_n] = | \dot{u}_n |^2_{L^2_{\tau_n}} = | \partial \F_n |^2_{L^2_{\tau_n}} ( t_{n,k+1}, u_{n,k+1}) ,
$$
which concludes the proof. \qed

\begin{lemma} \label{p.enest3} For every $0 \le k \le n-1$, the following energy estimate holds
\begin{align*}
	\E (  u_n ( t_{n,k+1} ) )  & \le \E ( u_n ( t_{n,k} ) ) -  \int_{t_{n,k}}^{t_{n,k+1}} \tfrac12 | \dot{u}_{n} (t) |^2_{L^2_{\tau_n}}  + \tfrac12 | \partial \F_n |_{L^2_{\tau_n}}^2  ( t , u^\sharp_{n} (t) )  \, dt \ + \nonumber  \\ & 
+ \int_{t_{n,k}}^{t_{n,k+1}} \langle f_n (t) , \dot{u}_n (t) \rangle \, dt %- \tau_n \int_{t_{n,k+1}}^{t_{n,k+1}}  \| \dot{u}_n  (t) \|^2_{H^1_0} \, dt  \,. \label{e.enest-2}
\end{align*}
\end{lemma}

\proof It is enough to combine the proof of Lemma \ref{p.enest} with Lemma \ref{l.EL3}. \qed

\noindent {\bf Proof of Proposition \ref{p.mm3}.} As $| \dot{u}_n |_{L^2_{\tau_n}} \ge \| \dot{u}_n \|_{L^2}$, arguing as in the proof of Proposition \ref{p.mm}, it follows that $u_n \in H^1(0,T;L^2)$ and then $u_n \in L^\infty(0,T;H^1)$. Thus, upon extracting a subsequence (non relabelled) $u_n \weakto u$ in $H^1(0,T;L^2)$.

\separe

Given $t^* \in (0,T]$ let $k_n$ such that $t_{n,k_n} <  t^* \le t_{n,k_n+1}$; by the previous lemma we have
\begin{equation*} %\label{e.123}
	\E (  u_n ( t_{n,k_n+1} ) )   + \int_{0}^{t_{n,k_n+1}} \hspace{-5pt} \tfrac12 |\dot{u}_{n} (t) |^2_{L^2_{\tau_n}}  +  \tfrac12 | \partial \F_n |^2_{L^2_{\tau_n}}  ( t , u^\sharp_{n} (t) ) + \langle f_n(t) , \dot{u}_n (t)  \rangle  \, dt  \le \E (u_0 ) \,.
\end{equation*}
As in the proof of Proposition \ref{p.mm} we get that $u^\sharp_n (t) \weakto u(t)$ in $H^1_0$ and thus by convexity of the energy $\E (  u (t^*) ) \le \liminf_{n \to +\infty} \E (  u_n (  t_{n,k_n+1} ) ) $. 

Note that $| \dot{u}_n |_{L^2_{\tau_n}} \ge | \dot{u}_n |_{L^2_\tau}$ for every $\tau_n \le \tau$ and that $| \cdot |^2_{L^2_\tau}$ is positive and convex. Then, by weak convergence in $H^1(0,T; L^2)$ for every $\tau>0$ we get 
\begin{equation}  \label{e.wlag}
	\int_0^{t^*} |\dot{u} (t) |^2_{L^2_\tau} \, dt \le \liminf_{n \to +\infty} \int_0^{t^*} | \dot{u}_n (t) |^2_{L^2_\tau} \, dt \le \liminf_{n \to +\infty} \int_0^{t^*} | \dot{u}_n (t) |^2_{L^2_{\tau_n}} \, dt \le  \liminf_{n \to +\infty} \int_0^{t_{n,k_n+1}} | \dot{u}_n (t) |^2_{L^2_{\tau_n}} \, dt .
\end{equation}
 Note that 
$$	\sup_{\tau >0} | z |_{L^2_\tau} = | z |_{L^2_+} = \begin{cases} \| z \|_{L^2} & \text{if $z \ge 0$}  \\ +\infty & \text{otherwise.} \end{cases}  $$
Then, taking the supremum w.r.t. $\tau>0$ in \eqref{e.wlag} by monotone convergence we get %$\dot{u} \ge 0$ and 
$$
	\int_0^{t^*} | \dot{u} (t) |^2_{L^2_+} \, dt \le  \liminf_{n \to +\infty} \int_0^{t_{n,k_n+1}} | \dot{u}_n (t) |^2_{L^2_{\tau_n}} \, dt . 
$$
Since $u^\sharp_n (t) \weakto u(t)$ in $H^1_0$, by Corollary \ref{c.lsc-slope} together with $| \partial \F_n |_{L^2_+} (t, \cdot) \le | \partial \F_n |_{L^2_\tau} (t, \cdot )$ we obtain %that 
%As a consequence, if $\tau_m \to 0$, $t_m \to t$ and $u_m \weakto u$ in $H^1_0(\Omega)$ then by Corollary  \label{c.lsc-slope}
\begin{equation*}% \label{e.lsctau3}
	| \partial \F |_{L^2_+} ( t,u(t)) \le \liminf_{n \to +\infty} | \partial \F_n |_{L^2_+} ( t , u_n^\sharp(t) )  \le \liminf_{n \to +\infty} | \partial \F_n |_{L^2_{\tau_n}} ( t  , u_n^\sharp(t)) .
\end{equation*}
Passing to the limit in the discrete energy estimate, as in the proof of Proposition \ref{p.mm}, we conclude the proof. \qed

\subsection{Energy identity, uniqueness and strong convergence \label{s.uni}}

Before proving energy identity and uniqueness, we give a short proof of the following Lemma.

\begin{lemma} \label{l.int} Let $w \in H^1 ( 0, T ; L^2)$. Consider a sequence of finite subdivisions $t_{j,i} $ of $[0,T]$ with $0 = t_{j,0} < ... < t_{j,i} < t_{j,i+1} < ... < t_{j,I_j} = T$ and let $\tau_j =\max_i  (t_{j,i+1}  - t_{j,i} )$.  Let $w_j$ be the piecewise affine interpolant of $w$ in the points $t_{j,i}$. Then
\begin{gather}  \label{e.12norm}
	\int_0^T \| w_j - w \|^2_{L^2} \, dt \le 4 \tau_j^2 \int_0^T \| w' \|^2_{L^2} \, dt \,, \qquad 
	\int_0^T \| w'_j \|^2_{L^2} \, dt \le \int_0^T \| w' \|^2_{L^2} \, dt   \,.
\end{gather}
If $\tau_j \to 0 $ then $w_j \weakto w $ in $H^1 ( 0, T ; L^2)$ and $\| w'_j \|_{L^2} \weakto \| w' \|_{L^2}$ in $L^2(0,T)$.
\end{lemma}

\proof For $t \in (t_{j,i} , t_{j,i+1})$ we have $w'_j = \dashint_{t_{j,i}}^{t_{j,i+1}} w'(t) \, dt$. Hence, by Jenssen's inequality
\begin{equation} \label{e.that}
	\int_{t_{j,i}}^{t_{j,i+1}}  \! \!  \| w'_j (t) \|_{L^2}^2 \, dt = ( t_{j,i+1} - t_{j,i} ) \,  \| w'_j \|^2_{L^2} = ( t_{j,i+1} - t_{j,i} ) \,  \Big\|  \, \dashint_{t_{j,i}}^{t_{j,i+1}} \! \! \! w'(t) \, dt \,  \Big \|_{L^2}^2 \le  \int_{t_{j,i}}^{t_{j,i+1}}  \! \!  \| w'(t) \|_{L^2}^2 \, dt .
\end{equation}
Taking the sum for $i=0,...,I_j-1$ yields the second estimate in \eqref{e.12norm}. Using \eqref{e.that} for $t \in (t_{j,i} , t_{j,i+1})$ we can write 
\begin{align*}
	\| w_j (t) - w(t) \|^2_{L^2} & = \Big\|  \int^{t}_{t_{j,i}} \!\! w'_j (s) - w'(s)   \, ds \,  \Big\|^2_{L^2} \le 
\Big(  \int^{t_{j,i+1}}_{t_{j,i}} \! \| w'_j (s) \|_{L^2} + \| w'(s) \|_{L^2}  \, ds \Big)^2 \\ 
& \le 2 \, (t_{j,i+1} - t_{j,i})  \int_{t_{j,i}}^{t_{j,i+1}} \! \| w'_j (s) \|^2_{L^2} + \| w'(s) \|^2_{L^2} \, ds   \le  4 \, \tau_j    \int_{t_{j,i}}^{t_{j,i+1}} \! \!  \| w'(s) \|_{L^2}^2 \, ds .
\end{align*}
Hence
$$
	 \int^{t_{j,i+1}}_{t_{j,i}} \| w_j (t) - w(t) \|^2_{L^2} \, dt \le 4 \, \tau_j^2  \int_{t_{j,i}}^{t_{j,i+1}} \! \!  \| w'(t) \|_{L^2}^2 \, dt .
$$
Taking the sum for $i=0,...,I_j-1$ yields the first estimate in \eqref{e.12norm}. 

From \eqref{e.12norm} it is clear that $w_j \weakto w $ in $H^1 ( 0, T ; L^2)$ for $\tau_j \to 0$. Let us see that $\| w'_j \|_{L^2} \weakto \| w' \|_{L^2}$ in $L^2(0,T)$. First, we show that $\| w'_j \|_{L^2} \to \| w' \|_{L^2}$ a.e.~in $(0,T)$. For $t \in (t_{j,i} , t_{j,i+1})$ 
\begin{align*}
	\Big|  \| w'(t) \|_{L^2} - \| w'_j (t) \|_{L^2}  \Big|   & \le \|  w'(t) - w'_j(t) \|_{L^2}  \le \Big\| w'(t) - \dashint^{t_{j,i+1}}_{t_{j,i}} \!\!\! w'(s) \, ds \, \Big\|_{L^2} \\
	& \le \Big\|\dashint^{t_{j,i+1}}_{t_{j,i}}  w'(t) -  w'(s) \, ds \,\Big\|_{L^2} \le \dashint^{t_{j,i+1}}_{t_{j,i}}  \| w'(t) -  w'(s) \|_{L^2} \, ds \, \\
& \le 2 \, \dashint_{t -  | t_{j,i+1} - t_{j,i} |}^{t + | t_{j,i+1} - t_{j,i} |} \| w'(t) -  w'(s) \|_{L^2} \, ds . \,
\end{align*}
It is well known (see e.g.~\cite[Proposition 2.1.22]{GasinskiPapageorgiou06}) that as  $| t_{j,i+1} - t_{j,i} | \to 0$ the last term is infinitesimal for a.e.~$t\in (0,T)$.

Since $\| w'_j \|_{L^2} \to \| w' \|_{L^2}$ a.e.~in $(0,T)$ and $\| w'_j \|_{L^2}$ is bounded in $L^2 (0,T)$ we know, by classical results, that $\| w'_j \|_{L^2} \to \| w' \|_{L^2}$ strongly in $L^1(0,T)$ and thus weakly in $L^2(0,T)$. \qed

\begin{lemma} \label{l.kan} Let $u \in H^1(0,T ; L^2)$ with $ | \partial \F |_{L^2_+} ( \cdot , u (\cdot) ) \in L^2 (0,T)$. Then  for every $0 \le t^* \le T$ 
\begin{align*}
       \E ( u(t^*)) ) % & \ge \F (0, u_0) - \int_0^{t^*}  | \partial_+ \F| ( t, u ( t ) ) \| \dot{u} (t) \|_{L^2} \, dt  - \int_0^{t^*}  \langle f'(t) , u ( t ) \rangle \, dt \nonumber \\ &
 & \ge \E ( u_0 ) -  \int_0^{t^*}  | \partial \F|_{L^2_+} ( t, u ( t ) ) \, \|  \dot{u} (t) \|_{L^2_+} \, dt - \int_0^{t^*}  \langle f(t) ,\dot{u} ( t ) \rangle \, dt \\ 
& \ge \E ( u_0 ) -  \tfrac12 \int_0^{t^*}  | \partial \F|_{L^2_+}^2 ( t, u ( t ) ) + \| \dot{u} (t) \|^2_{L^2_+} \, dt - \int_0^{t^*}  \langle f(t) , \dot{u} ( t ) \rangle \, dt .
\end{align*}
\end{lemma}

\proof The lack of time regularity in $H^1_0$ prevents from employing the chain rule, we will use instead a Riemann sum argument, adapted from \cite[Lemma 4.12]{DalMFrancfToad05}, see also \cite[Proposition 3.8]{Negri_ACV19}.
%\cite[Proposition 3.8]{Negri_16}. 
%Note that by \eqref{e.enest} the slope $| \partial_+ \F | ( \cdot , u ( \cdot))$ is bounded in $L^2(0,T)$. 
Note that $ | \partial \F |_{L^2_+} ( \cdot , u (\cdot) ) < +\infty$ a.e.~in $(0,T)$.

Given $t^*$ let $0 < t_* < t^*$ with $| \partial \F |_{L^2_+} ( t_* , u (t_*) ) < +\infty$. We can find a sequence of finite subdivisions $\mathcal{T}_j = \{ t_{j,i} \}$ of $[t_*,t^*]$ with $t_* = t_{j,0} < ... < t_{j,i} < t_{j,i+1} < ... < t_{j,I_j} = t^*$, such that $\lim_{j \to +\infty}  \max_i \{ t_{j,i+1}  - t_{j,i} \} = 0 $ and\footnote{It is enough to apply the Riemann sum argument to $( | \partial \F|_{L^2_+} , f ) \in L^2(0,T; \mathbb{R} \otimes L^2)$}  
\begin{gather} \label{e.Riemann}
	S_j (\cdot) = \sum_{i=0}^{I_j-1}  | \partial \F|_{L^2_+} ( t_{j,i} , u (t_{j,i} )  )  \text{\large $\chi$}_{(t_{j,i} , t_{j,i+1} )} (\cdot) \to | \partial \F |_{L^2_+} ( \cdot , u (\cdot) ) 
	\quad \text{strongly in $L^2(t_*,t^*)$,}  \\ 
	F_j (\cdot) = \sum_{i=0}^{I_j-1}  f ( t_{j,i} )  \text{\large $\chi$}_{(t_{j,i} , t_{j,i+1} )} (\cdot) \to f 
	\quad \text{strongly in $L^2(t_*,t^*; L^2)$.}
 \label{e.fsum}
\end{gather}
(For sake of clarity, we remark that the points $\{ t_{j,i} \}$ do not coincide with the points $t_{n,k} = n \tau_n$ appearing in the discrete scheme).  By convexity of $\F ( t_{j,i+1} , \cdot)$ we write
\begin{align*}
          \F ( t_{j,i} , u ( t_{j,i+1} ) ) & \ge \F ( t_{j,i} , u ( t_{j,i} ) ) + d \F ( t_{j,i} , u ( t_{j,i} ) ) [ u ( t_{j,i+1} ) -  u ( t_{j,i} ) ] \\
	 & \ge \F ( t_{j,i} , u ( t_{j,i} ) ) - | \partial \F|_{L^2_+} ( t_{j,i} , u ( t_{j,i} ) ) \| u ( t_{j,i+1} ) -  u ( t_{j,i} ) \|_{L^2} .
%	& \ge  \F ( u ( t_{j,i} ) ) - | \partial_+ \F| ( u ( t_{j,i} ) ) \|  \tilde{u}'_j (t) \|_{L^2} (t_{j,i+1} - t_{j,i})  .
%	& \ge \F ( u ( t_{j,i} ) ) - \tfrac12 | \partial_+ \F|^2  ( u ( t_{j,i} ) )   (t_{j,i+1} - t_{j,i})   -  \tfrac12   \|  \tilde{u}'_n \|^2_{L^2} (t_{j,i+1} - t_{j,i}). 
\end{align*}
Denote by $u_j$ the piecewise affine interpolant of $u(t_{j,i})$.
%Let us introduce (by abuse of notation) $u^\flat_j$,  $u^\sharp_j$ and $u_j$ respectively as the forward piecewise constant, backward piecewise constant and piecewise affine interpolant of $u$ in the points $t_{j,i}$. Similarly, let $t^\flat_{j}$ be the forward piecewise affine interpolant of $t_{j,i}$ in the points $t_{j,i}$.  For every $t \in ( t_{j,i} , t_{j,i+1})$ we have 
Writing explicitely $\F ( t_{j,i} , u ( t_{j,i+1} ) )$ and $\F ( t_{j,i} , u ( t_{j,i} ) )$ we get 
$$
  \E ( u ( t_{j,i+1} )  )  \ge \, \E ( u ( t_{j,i} )  )  - \int_{t_{j,i}}^{t_{j,i+1}} | \partial \F|_{L^2_+} ( t_{j,i} , u ( t_{j,i}) ) ) \|\dot{u}_j (t) \|_{L^2} \, dt -  \int_{t_{j,i}}^{t_{j,i+1}}  \langle f(t_{j,i}) , \dot{u}_j (t) \rangle  \, dt .  
$$
Using the above estimate for $i=1,...,I_j$ we get, in terms of the functions $S_j$ and $F_j$, 
\begin{equation} \label{e.gigi}
\E (  u(t^*) ) \ge \E ( u (t_*) ) -  \int_{t_*}^{t^*}  S_j (t)  \|\dot{u}_j (t) \|_{L^2} \, dt
 - \int_{t_*}^{t^*} \langle F_j (t) , \dot{u}_j (t) \rangle \, dt  .
\end{equation} 
By \eqref{e.Riemann} we known that $S_j  ( \cdot ) \to | \partial\F |_{L^2_+} ( \cdot , u ( \cdot ))$ strongly in $L^2(t_*,t^*)$ and by \eqref{e.fsum} that $F_j  \to f$  strongly in $L^2(t_*,t^*; L^2)$. By Lemma \ref{l.int} we get  $\dot{u}_j \weakto \dot{u} $ in $L^2(t_*,t^* ; L^2 )$ and $\|\dot{u}_j \|_{L^2} \weakto \|\dot{u} \|_{L^2}$ in $L^2 (t_*,t^*)$.  In summary, we can pass to the limit in \eqref{e.gigi} and get, by Young's inequality,  
\begin{align*}
\E( u(t^*) ) & \ge \, \E ( u (t_*) ) - \int_{t_*}^{t^*}  | \partial \F|_{L^2_+} ( t, u ( t ) ) \| \dot{u} (t) \|_{L^2} \, dt  - \int_{t_*}^{t^*}  \langle f(t) , \dot{u} ( t ) \rangle \, dt 
\nonumber \\ & \ge  \, \E ( u (t_*) ) -  \tfrac12 \int_{t_*}^{t^*}  | \partial \F|^2_{L^2_+} ( t, u ( t ) ) + \| \dot{u} (t) \|^2_{L^2} \, dt - \int_{t_*}^{t^*}  \langle f(t) , \dot{u} ( t ) \rangle \, dt. % \label{e.dise}
\end{align*}
Taking the liminf of the right hand for $t_* \to 0^+$ we get 
$$
\E( u(t^*) )  \ge  \E ( u_0 ) -  \tfrac12 \int_{0}^{t^*}  | \partial \F|^2_{L^2_+} ( t, u ( t ) ) + \| \dot{u} (t) \|^2_{L^2} \, dt - \int_0^{t^*}  \langle f(t) , \dot{u} ( t ) \rangle \, dt ,
$$
which concludes the proof.  \qed  %is the opposite inequality of 

%\separe

\medskip
{\bf Energy identity.} Clearly, using \eqref{e.evol-ineq} and Lemma \ref{l.kan} it follows that for every $0 \le  t^*  \le T$ we get 
\begin{equation} \label{e.enide2}
\E (  u(t^*) )  %= - \tfrac12 \int_{t_1}^{t_2} \|\dot{u} (t) \|_{L^2} ^2 + | \partial_+ \F |^2 ( u (t) ) \, dt 
= \E ( u_0 ) - \tfrac12 \int_{0}^{t^*} \|\dot{u} (t) \|^2_{L^2_+} + | \partial \F |^2_{L^2_+} ( t, u (t) ) \, dt   - \int_{0}^{t^*}  \langle f(t) , \dot{u} ( t ) \rangle \, dt.
\end{equation}
As a consequence, the energy identity holds in every subinterval $[t_1 , t_2]$ of $[0,T]$.

%\subsection{Uniqueness}

%\separe

% --------------------------------------------------------------------------------------------------------------------------------------- ARGUMENT OF GIGLI 
%Since the energy $\F$ is convex by \cite[Corollary 2.4.10]{AmbrosioGigliSavare05} the slope is a strong upper gradient, thus for every curve $v \in AC ( 0, T ; L^2 )$ it holds 
%\begin{equation}  \label{e.FALSE!!}
%	| \F ( v (t_2) ) - \F ( v (t_1) ) | \le  \int_{t_1}^{t_2} | \partial_+ \F | ( v (t) ) \| v ' (t) \|_{L^2}  \, dt 
%\end{equation}
%for every $ 0 < t_1 < t_2 < T$; hence 
%$$
%	\F ( v (t_1) ) \le  \F ( v (t_2) )  + \tfrac12 \int_{t_1}^{t_2}  | \partial_+ \F |^2 ( v (t) )  + \| v' (t) \|^2_{L^2}  \, dt .
%$$
%If $v(0) = u_0$ taking the liminf for $t_1 \to 0^+$ in the left hand side and using the lsc of $\F$ we get
%\begin{equation}  \label{e.ine-v0}
%	\F ( u_0 ) \le  \F ( v (t_2) )  + \tfrac12 \int_{0}^{t_2}  | \partial_+ \F |^2 ( v (t) )  + \| v' (t) \|^2_{L^2}  \, dt .
%\end{equation}%
%\separe
%In particualr  curves $u$ of maximal slope are of class $AC ( 0,T ; L^2)$ with $u(0)=u_0$, hence for every $t^* >0$ %, again by \eqref{e.upgrad},t
%$$
%	\F ( u_0 ) \le  \F ( u (t^*) )  + \tfrac12 \int_{0}^{t^*}  | \partial_+ \F |^2 ( u (t) )  + \|\dot{u} (t) \|^2_{L^2}  \, dt ,
%$$
%which is the opposite of \eqref{e.evol-ineq}; thus  for every $t^* \in (0,T]$ we get the energy identity
%\begin{equation} \label{e.ide}
%     \F ( u_0 ) = \F ( u(t^*) ) +  \tfrac12 \int_{0}^{t^*} | \partial_+ \F |^2 ( u (t) )  + \|\dot{u} (t) \|^2_{L^2}  \, dt  .
%\end{equation}

\bigskip
{\bf Uniqueness.} Remember that in our weak setting solutions belongs only to $H^1(0,T;L^2) \cap L^\infty (0,T;H^1_0)$, therefore we are not in a position to employ any argument based on the chain rule for $\F$. Instead, we follow the contradiction argument of \cite[Theorem 15]{Gigli_CVPDE10}. Assume that $u_\I$ and $u_\II$ are different unilateral gradient flows with the same initial value $u_0$. Let $t^*$ such that $u_\I(t^*) \neq u_\II (t^*)$. Define $ u_\natural = \tfrac12 ( u_\I + u_\II)$. Writing the energy identity \eqref{e.enide2} for both $u_\I$ and $u_\II$ we get (for $i=\text{\small \rm I, II}$)
$$
   \tfrac12 \E ( u_0 ) = \tfrac12 \E ( u_i(t^*) ) + \tfrac12 \int_{0}^{t^*} \tfrac12 | \partial \F |^2_{L^2_+} ( t, u_i (t) )  + \tfrac12 \|\dot{u}_i (t) \|^2_{L^2}  \, dt
+ \tfrac12 \int_{0}^{t^*} \langle f(t) , \dot{u}_i(t) \rangle \, dt \, .
$$
Taking the sum for $i=\text{\small \rm I, II}$ and using the strict convexity of the energy $\E$, the  convexity of $\tfrac12 \| \cdot \|^2_{L^2}$, the convexity of $\tfrac12 | \partial \F |^2_{L^2_+} ( t, \cdot) $ (see Corollary \ref{c.lsc-slope}) and the linearity of $ \langle f(t) , \cdot \rangle$ we get 
$$
   \E (  u_0 )  >  \E ( u_\natural (t^*) ) + \tfrac12 \int_{0}^{t^*} | \partial \F |^2_{L^2_+} ( t, u_\natural (t) )  + \| \dot{u}_\natural (t) \|^2_{L^2}  \, dt  +  \int_{0}^{t^*} \langle f(t) , \dot{u}_\natural (t) \rangle \, dt \, .
$$
Hence $ | \partial \F |_{L^2_+} ( \cdot , u_\natural (\cdot) )$ belongs to $L^2(0,t^*)$ and clearly $u_\natural \in H^1(0,T; L^2)$. The previous inequality is a contradiction with Lemma \ref{l.kan}. 

Since the limit evolution is unique it is not necessary to extract any subsequence in Propositions \ref{p.mm}, \ref{p.mmbis} and \ref{p.mm3}.

\bigskip
{\bf Strong convergence.} To conclude, let us check that $u_n (t) \to u(t)$ in $H^1_0$ pointwise in $[0,T]$, where $u_n$ is the sequence provided by the discrete scheme of \S\,\ref{3.1}; the same property holds, with few changes, for the sequences $u_n$ of \S\,\ref{s.3-2} and \S\,\ref{s.pen}. 

Given $t^* \in [0,T]$ let us first prove that $u_n^\sharp (t^*) \to u(t^*)$ in $H^1_0$. Since $u^\sharp_n (t^*) \weakto u(t^*)$ in $H^1_0$ it is enough to show that $\E ( u^\sharp_n (t^*)) \to \E ( u(t^*) ) $, which implies that $u^\sharp_n (t^*) \to u(t^*)$ in $H^1_0$ endowed with the energy norm. Let $k_n$ s.t. $ t_{n,k_n} < t^* \le t_{n,k_n+1}$ and recall \eqref{e.123}, i.e.
$$
	\E ( u^\sharp _n (t^*) ) = \E (  u_n ( t_{n,k_n+1} ) ) \le \E ( u_0 ) - \int_{0}^{t_{n,k_n+1}} \hspace{-5pt} \tfrac12 |\dot{u}_{n} (t) |^2_{L^2_+}  +  \tfrac12 | \partial \F_n |^2_{L^2_+}  ( t , u^\sharp_{n} (t) ) + \langle f_n(t) , \dot{u}_n (t)  \rangle  \, dt  \,.
$$
Taking the limsup yields
\begin{align*}
	\E ( u(t^*)   ) & \le  \, \liminf_{n \to +\infty} \E ( u_n^\sharp (t^*) ) \le \limsup_{n \to +\infty} \E ( u_n^\sharp (t^*) ) \\
%	& \le \limsup_{n \to +\infty} \left(  \E ( u_0 ) - \int_{0}^{t_{n,k_n+1}} \hspace{-5pt} \tfrac12 |\dot{u}_{n} (t) |^2_{L^2_+}  +  \tfrac12 | \partial \F_n |^2_{L^2_+}  ( t , u^\sharp_{n} (t) ) + \langle f_n(t) , \dot{u}_n (t)  \rangle  \, dt \right) \\
	& \le \E ( u_0 ) - \liminf_{n \to +\infty } \tfrac12  \int_{0}^{t_{n,k_n+1}} \hspace{-7pt} |\dot{u}_{n} (t) |^2_{L^2_+}  +  | \partial \F_n |^2_{L^2_+}  ( t , u^\sharp_{n} (t) ) \, dt + \lim_{n \to +\infty } \int_{0}^{t_{n,k_n+1}}  \langle f_n(t) , \dot{u}_n (t)  \rangle  \, dt \\
	& \le \E ( u_0 ) - \tfrac12 \int_{0}^{t^*} |\dot{u} (t) |^2_{L^2_+} + | \partial \F |^2_{L^2_+} ( t, u (t) ) \, dt   - \int_{0}^{t^*}  \langle f(t) , \dot{u} ( t ) \rangle \, dt = \E (  u(t^*) ) , %= - \tfrac12 \int_{t_1}^{t_2} \|\dot{u} (t) \|_{L^2} ^2 + | \partial_+ \F |^2 ( u (t) ) \, dt 
\end{align*}
where, in the last line, we used (\ref{e.alprim}-\ref{e.alters}) from Proposition \ref{p.mm} together with the energy identity \eqref{e.enide2}. As a consequence, all inequalities above turn into equalities and $u^\sharp_n (t^*) \to u(t^*)$ in $H^1_0$; hence $ u_n ( t_{n,k_n+1} ) = u^\sharp _n (t^*) \to u ( t^*)$ in $H^1_0$. A similar argument shows that $u_n ( t_{n,k_n} )   \to u ( t^*)$ in $H^1_0$. Being $u_n (t^*)$ a convex combination of $u_n ( t_{n,k_n} ) $ and $u_n ( t_{n,k_n+1} )$ it converges strongly to $u (t^*)$ as well.

% !TeX root = article.tex

\section{Further properties of solutions}

\subsection{Comparison principle}

Since the unilateral gradient flow is unique it is enough to prove the maximum principle for the discrete solutions provided in \S\,\ref{s.3-2}. To this end, fix $\tau_n>0$ and assume that $u_{n,0} = u_0 \le v_0 = v_{n,0}$. We will show by induction that $u_{n,k} \le v_{n,k}$ for every index $k \ge 1$. We recall that $u_{n,k+1} = \max \{  \tilde u_{n,k+1} , u_{n,k}  \}$ and that 
$$
	\tilde u_{n,k+1}  \in \argmin \big\{ \F ( t_{n,k+1} , u  ) + \tfrac{1}{2\tau_n} \, \| u - u_{n,k} \|^2_{L^2} :  u \in H^1_0 \big\} . 
$$                               
Assume by induction that $u_{n,k} \le v_{n,k}$, we claim that $\tilde u_{n,k+1} \le \tilde v_{n,k+1}$ from which we get  $u_{n,k+1} \le v_{n,k+1}$. By minimality there exists $\xi_{n,k+1} \in \partial \F ( t_{n,k+1} , \tilde u_{n,k+1} ) \subset H^{-1}$ s.t.~$\xi_{n,k+1} + \frac{1}{\tau_n} ( \tilde u_{n,k+1} - u_{n,k}) = 0$ in $H^{-1}$. In a similar way, there exists $\zeta_{n,k+1} \in \partial \F ( t_{n,k+1} , \tilde v_{n,k+1} )$ such that  $\zeta_{n,k+1} + \frac{1}{\tau_n} ( \tilde v_{n,k+1} - v_{n,k}) = 0$ in $H^{-1}$.  Hence, using $[  \tilde u_{n,k+1} - \tilde v_{n,k+1} ]_+ \in H^1_0$ as a test function we get
$$ 
	\big(  \xi_{n,k+1}   -   \zeta_{n,k+1}  ,  [  \tilde u_{n,k+1} - \tilde v_{n,k+1} ]_+  \big)  +  \tfrac{1}{\tau_n} \big\langle  ( \tilde u_{n,k+1} - \tilde v_{n,k+1}) - (  u_{n,k} - v_{n,k} ) , [  \tilde u_{n,k+1} - \tilde v_{n,k+1} ]_+  \big\rangle  = 0.
$$
The first term is non-negative by $T$-monotonicity (see Remark \ref{r.en}) hence for the last term we can write 
$$
	0 \le \int_{\Omega} [  \tilde u_{n,k+1} - \tilde v_{n,k+1} ]_+^2 \, dx \le \int_\Omega (  u_{n,k} - v_{n,k} ) [  \tilde u_{n,k+1} - \tilde v_{n,k+1} ]_+ \, dx .
$$
Since $u_{n,k} \le  v_{n,k}$ the integrand in the right hand side is non-positive; it follows that $ [ \tilde u_{n,k+1} - \tilde v_{n,k+1} ]_+ = 0$ and thus $ \tilde u_{n,k+1}  \le \tilde v_{n,k+1}$.

\separe

If $u_{n,k} \le v_{n,k}$ for every $k \ge 0$ then $u_n \le v_n $ in $[0,T]$; passing to the limit weakly in $H^1(0,T; L^2)$ we get $u \le v $  in $[0,T]$.

\subsection{Continuous dependence \label{s.contdep}} 

In this section we will prove Proposition \ref{p.contdep}. We adopt the scheme of \cite{SandSerf04}. By definition we know that 
%%$\F ( u ( \cdot))$ is monotone non-increasing in $[0,T)$ and if 
for every $0 \le t^* \le T$ it holds 
\begin{equation}\label{e.evol-ineq-dc}
      \E ( u^m (t^*) ) \le \E ( u^m_0 )  -  \, \tfrac12 \int_{0}^{t^*} \| \dot{u}^m (t) \|_{L^2}^2 + \| [ A u^m(t) + f^m(t) ]_+   \|^2_{L^2} \, dt + \int_0^{t^*} \langle f^m (t) , \dot{u}^m (t) \rangle \, dt .
\end{equation}
Denote $\F^m ( t , u) = \E ( u) - \langle u , f^m (t) \rangle$.  We recall that by Proposition \eqref{p.pde} and Corollary \ref{c.B1} we have $\| \dot{u}^m (t) \|_{L^2} = \| [ A u^m(t) + f^m(t) ]_+ \|_{L^2} = | \partial \F^m |_{L^2_+} ( t, u^m(t)) $. Hence, choosing $t^* =T$ above we obtain 
$$
	\E ( u^m_0) \ge \int_{0}^{T} \| \dot{u}^m (t) \|_{L^2}^2  - \int_0^{T} \langle f^m (t) , \dot{u}^m (t) \rangle \, dt \ge \| \dot{u}^m \|^2_{L^2(0,T;L^2)} - \| f^m \|_{L^2(0,T;L^2)} \, \| \dot{u}^m \|_{L^2(0,T;L^2)} .
$$ 
Since $u^m_0 \to u_0$ in $H^1_0$ we have $\E (u^m_0) \to \E ( u_0)$; since $f^m \to f$ in $L^2(0,T;L^2)$ the above estimate implies that the sequence $\dot{u}^m$ is bounded in $L^2(0,T;L^2)$. since $u^m_0 \to u_0$ in $H^1_0$ it follows that $u_m$ is also bounded in $H^1(0,T;L^2)$. Moreover, by \eqref{e.evol-ineq-dc} and by coercivity of the stored energy, for every $0 \le t^* \le T$ we have 
$$
	c \| u^m (t^*) \|^2_{H^1_0} \le \E ( u^m (t^*) ) \le \E ( u^m_0 ) + \| f^m \|_{L^2(0,T;L^2)} \, \| \dot{u}^m \|_{L^2(0,T;L^2)} \le C .
$$
Hence, the sequence $u_m$ is bounded also in $L^\infty (0,T; H^1_0)$. In conclusion, there exists a subsequence (non relabelled) such that $u^m \weakto u$ in $H^1(0,T;L^2)$, as consequence $u_m$ is monotone non-decerasing. 
Moroever, arguing as in the proof of Proposition \ref{p.mm}, we get that $u_m (t) \weakto u(t)$ in $H^1_0$ for every $t \in [0,T]$. 

It remains to show that $u$ is the unilateral gradient flow for $\F$ with initial condition $u_0$. By \eqref{e.evol-ineq-dc} for every $t^* \in (0,T]$ we can write 
\begin{align*}
	\E ( u (t^*) ) & \le \liminf_{m \to +\infty} \E ( u^m (t^*)) \\
	& \le \limsup_{m \to +\infty} \bigg( \E ( u^m_0 )  -  \, \tfrac12 \int_{0}^{t^*} \| \dot{u}^m (t) \|_{L^2}^2 + \| [ A u^m(t) + f^m(t) ]_+   \|^2_{L^2} \, dt + \int_0^{t^*} \langle f^m (t) , \dot{u}^m (t) \rangle \, dt \bigg)  \\
	& \le \limsup_{m \to +\infty} \, \E ( u^m_0 ) - \tfrac12  \liminf_{m \to +\infty} \bigg( \int_{0}^{t^*} \| \dot{u}^m (t) \|_{L^2}^2 + | \partial \F^m |^2_{L^2_+}  ( t, u^m (t)) \, dt \bigg) \, + \\
	& \phantom{\le} + \limsup_{m \to +\infty} \int_0^{t^*} \langle f^m (t) , \dot{u}^m (t) \rangle \, dt .
\end{align*}
We know that $\E ( u^m_0 ) \to \E ( u_0 )$ because $u^m_0 \to u_0$ (strongly) in $H^1_0$. Since $\dot{u}^m \weakto \dot{u}$ in $L^2(0,T;L^2)$ we have $\| \dot{u} \|^2_{L^2(0,t^*;L^2)} \le \liminf_{m \to +\infty} \| \dot{u}^m \|^2_{L^2(0,t^*;L^2)}$. As $u^m \weakto u $ in $H^1_0$ a.e.~in $(0,T)$ we can apply Corollary \ref{c.lsc-slope} and then by Fatou's lemma we get
$$
	\int_{0}^{t^*} | \partial \F |^2_{L^2_+}  ( t, u(t)) \, dt \le \liminf_{m \to +\infty }
	\int_{0}^{t^*} | \partial \F^m |^2_{L^2_+}  ( t, u^m (t)) \, dt .
$$
Finally, $\int_0^{t^*} \langle f^m (t) , \dot{u}^m (t) \rangle \, dt  \to \int_0^{t^*} \langle f (t) , \dot{u} (t) \rangle \, dt $ by strong-weak convergence in $L^2(0,t^*;L^2)$. In conclusion, we get 
$$
      \E ( u (t^*) ) \le \E ( u (0) )  -  \, \tfrac12 \int_{0}^{t^*} | \dot{u} (t) |_{L^2_+}^2 + | \partial \F | _{L^2_+} ( t, u(t)) \, dt + \int_0^{t^*} \langle f (t) , \dot{u} (t) \rangle \, dt ,
$$
which is equivalent to \eqref{e.evol-ineq}.

Finally, in order to prove that $u_m (t) \to u(t)$ strongly in $H^1_0$ for every $t \in [0,T]$ it is enough to follow the proof of the strong convergence in \S\ref{s.uni}.

\subsection[{Parabolic equation and variational inequality in $L^2$}]{Parabolic equation and variational inequality in \boldmath{$L^2$} \label{PDE}}

\noindent {\bf Proof of (\ref{e.pde}).} Let $u$ be the unilateral gradient flow for $\F$ with initial datum $u_0$, in the sense of Definition \ref{d.curve}. By the energy identity we know that $\| \dot{u} (t) \|_{L^2} = | \partial \F |_{L^2_+} ( t, u (t))$ is a.e.~finite in $[0,T]$. Hence, by Corollary \ref{c.B1}, %we get that $- d \F( t, u(t)) = A u(t) + f(t)$ is a (locally finite) Radon measure and 
for a.e.~$t \in [0,T]$ we have 
\begin{equation} \label{e.etfir} \| \dot{u} (t) \|_{L^2} = | \partial \F |_{L^2_+} ( t, u (t)) = \| [ A u(t) + f(t)]_+ \|_{L^2} . \end{equation} 
Now, let us show that $ - d \F ( t , u (t) ) \le \dot{u}(t)$ in $H^{-1}$ for a.e.~$t$ in $[0,T]$. 
By uniqueness, we can rely on the discrete scheme of \S\,\ref{3.1}. By \eqref{e.dPDE} for every $t \in (t_{n,k} , t_{n,k+1})$ we have 
$$
    - d \F_n ( t , u^\sharp_{n} (t) )  [ \phi ] \le  \langle \dot{u}_{n} (t) , \phi \rangle  \quad \text{for every $\phi \in C^\infty_0$ with $\phi \ge 0$,}
$$
which reads, by symmetry of $a ( \cdot , \cdot)$, 
$$
	\langle u^\sharp_{n} (t) , A \phi \rangle + \langle f_n (t) , \phi \rangle  \le \langle \dot{u}_{n} (t) , \phi \rangle  . %= \int_\Omega \dot{u}_n (t) \phi \, dx .
$$
Thus, for every $0 \le t_1 < t_2 \le T$ we can write 
$$
	\int_{t_1}^{t_2} \langle u^\sharp_n(t)  , A \phi \rangle +  \langle f_n(t) , \phi \rangle \, dt \le \int_{t_1}^{t_2} \langle \dot{u}_{n} (t) , \phi \rangle \, dt \,.
$$
Passing to the limit, by the strong convergence of $u^\sharp_n (t)$ in $L^2(0,T;L^2)$  togheter with the strong convergence of $f_n$ in $L^2(0,T;L^2)$ and the weak convergence of $u_n$ in $H^1(0,T;L^2)$, we get
$$
	\int_{t_1}^{t_2} \langle u (t)  , A \phi \rangle +  \langle f (t) , \phi \rangle \, dt \le \int_{t_1}^{t_2} \langle \dot{u} (t) , \phi \rangle \, dt .
$$
Since the above inequality holds for any choice of $t_1 < t_2$, for a.e.~$t\in(0,T)$ we have
$$
     - d \F ( t , u(t) ) [ \phi ] =  \langle u (t)  , A \phi \rangle +  \langle f (t) , \phi \rangle \le \langle \dot{u} (t) , \phi \rangle ,
$$
which reads $( A u(t) + f(t) , \phi ) \le \langle \dot{u} (t) , \phi \rangle$. Therefore, applying Lemma \ref{l.ineqL2} we get $[ \Delta u(t) + f(t)]_+ \le \dot{u} (t) $. By \eqref{e.etfir}, we have $\| \dot{u} (t) \|_{L^2} = \| [ \Delta u(t) + f(t)]_+ \|_{L^2}$, it follows that $\dot{u} (t) =  [ \Delta u(t) + f(t)]_+$.

\bigskip

\noindent \textbf{Proof of (\ref{e.pvi}).} If $Au(t) + f(t)$ is a Radon measure and $\dot{u}(t) = [A u(t) + f(t) ]_+ \in L^2$ then for every $\phi \in C^\infty_0$ with $\phi \ge 0$ we have 
$$
	\langle \dot{u} (t) , \phi \rangle = \langle [A u(t) + f(t) ]_+ , \phi \rangle \ge ( A u(t) + f(t) , \phi ) \,.
$$

\subsection[Non-uniqueness for parabolic problems in $L^2$]{Non-uniqueness for parabolic problems in \boldmath{$L^2$} \label{non-uni}} 

%\subsection{Explicit solution of the unilateral gradient flow}

First of all, let us see that the set of solutions of the parobolic variational inequality \eqref{e.pvi} is larger than the set of solutions of the parabolic problem \eqref{e.pde}, when $f \in L^2$ is independent of time.  Let $u$ be a solution of \eqref{e.pde}, by the arguments of the previous sections we know that $u$ solves also \eqref{e.pvi}, i.e. 
$$
	( A u(t) + f , \phi ) \le \langle \dot{u} (t) , \phi \rangle 
	\qquad
	\text{for every $\phi \in H^1_0$ with $\phi \ge 0$.}
$$
Now, consider $u_\lambda (t) = u (\lambda t)$ for any $\lambda > 1$. Then 
$$
	( A u_\lambda (t) + f , \phi ) = ( A u (\lambda t) + f , \phi )  \le \langle \dot{u} (\lambda t) , \phi \rangle 
	= \lambda^{-1} \langle \dot{u}_\lambda (t) , \phi \rangle  
	\le \langle \dot{u}_\lambda (t) , \phi \rangle .
	%\text{for every $\phi \in H^1_0$ with $\phi \ge 0$.}
$$
Thus, any such $u_\lambda$ solves \eqref{e.pvi}.

\bigskip

Next, we provide an example in which the parabolic problem \eqref{e.pde} has many solutions, and thus it is not equivalent to Definition \ref{d.curve}. Let $u_0 \in H^1_0 (-1,1)$ be defined by $u_0(x) = 1 - | x|$.
%$$
%	 u_0 (x) = \begin{cases}
%	(x + L)/L & \text{if $x \in (-L,0)$} \\
%	(L - x) / L & \text{if $x \in (0,L)$.}
%	\end{cases}
%$$
We will denote $u'$ and $u''$ the first and second space derivatives, respectively. Clearly $u''_0 = - 2 \delta_0$,
%$$
%	u'_0 (x) = \begin{cases}
%	1/L & \text{if $x \in (-L,0)$,} \\
%	- 1 / L & \text{if $x \in (0,L)$,}
%	\end{cases}
%	\qquad
%	u''_0 = - (2 / L) \delta_0 \,,
%$$
where $\delta_0$ denotes Dirac's delta in the origin. Note that $u_0 >0$ in $(-1,1)$ and  $[u''_0]_+ = 0$. 

Let $f = u_0$ (independent of time) and consider the Dirichlet energy $\F : H^1_0(-1,1) \to \R$ given by
\begin{equation} \label{e.F1D}
	\F ( u ) = \tfrac12 \int_{(-1,1)} | u' |^2 \, dx - \int_{(-1,1)} f u \, dx = \E ( u ) - \langle f , u \rangle .
\end{equation}
Clearly $A u + f = u'' + f $. 
Let us define $u(t) = (1+t) u_0 $. Then, $\dot{u} (t) = u_0$ and $u'' (t) = (1+t) u''_0 = -2 (1+t)  \delta_0$. Since  $[ u''(t) + f ]_+ = [ - 2 (1+t) \delta_0 + u_0 ]_+ = u_0$, it turns out that $u$ solves
$$
\begin{cases}
		\dot{u} (t) = [ u'' (t) + f ]_+  & \text{in $(0,T)$} \\
		u (0) = u_0 .
\end{cases}
$$
On the other hand $u$ does not satisfy the energy identity, in the form \eqref{e.phys}, i.e.,
$$
	\E ( u (t) )  = \E ( u_0 ) - \int_0^t \| \dot{u} (t) \|^2_{L^2} \, dt  + \int_0^t  \langle f , \dot{u} (t) \rangle \, dt \,.
$$
Indeed, $\E ( u(t) ) = \tfrac12 (1+t)^2 \int_{(-1,1)} | u'_0 |^2 \, dx = (1 + t)^2$, $\E ( u_0 ) = 1$, $\| \dot{u} (t) \|^2_{L^2} =  \| u_0 \|^2_{L^2}$ while $\langle f , \dot{u} (t) \rangle = \| u_0 \|^2_{L^2}$.

%\begin{remark} \normalfont
%By Theorem \ref{t.exist} a curve of maximal unilateral slope is expected to satisfy 
%\begin{equation} \label{e.vs}
%	\F' (u(t)) = - \| u' (t) \|^2_{L^2} \quad \text{a.e.~in $(0,T)$.}
%\end{equation}
%A direct computation shows that 
%$$
%	\F ( u (t) ) = \tfrac12 \int_{(-L,L)} | \nabla u(t) |^2 - f u(t) \, dx  = \tfrac1L (1+t)^2 - \tfrac23 (1+t) L \,.
%$$
%Hence $\F' ( u(t)) = \tfrac2L (1+t) - \tfrac23 L$. Note that for $t \approx 1$ we have $\F' (u(t)) > 0$ if $L \ll 1$ and $\F' (u(t)) < 0$ if $L \gg 1$. Let us check that $u$ is not a curve of maximal slope, indeed, 
%%$$
%%	| \partial_+ \F |^2 ( u(t) ) = \| [ u'' + u_0 ]_+ \|^2_{L^2} = \| \dot{u} \|^2_{L^2}  = \| u_0 \|^2_{L^2} =  2L/3 
%%$$
%%and thus for $t >0$
%$$
%	\F' ( u(t))  = \tfrac2L (1+t) - \tfrac23 L > - \tfrac23 L  = - \| u_0 \|^2_{L^2} =  - \| u' (t) \|^2_{L^2} ,
%$$
%which contradicts \eqref{e.vs}. 
%It is also interesting to write the derivative of the energy, using the chain rule, in the following way:
%$$
%    \F' ( u(t)) = \int_{(-L,L)} D u(t)  Du' (t) -  f u' (t) \, dx = - ( \Delta u(t) + f , u' ) ,
%$$
%where $( \cdot , \cdot)$ denotes the duality between $H^{-1} (-L,L)$ and $H^1_0 (-L,L)$. Then, writing $\Delta u(t) + f  = -(1+t) \tfrac2L \delta_0 + u_0$ and $u'=u_0$  we can write
%$$
%	- ( \Delta u(t) + f , u') =  - ( u' , u' )  + ( (1+t) \tfrac2L \delta_0 , u_0 ) = - \| u ' \|_{L^2}^2 + (1+t) \tfrac2L  
%$$
%since $u_0 (0) = 1 $.
%\end{remark}

  \bigskip
\noindent{\bf Solution of the unilateral gradient flow.\label{explicit}} 
%The solution of the unilateral gradient flow is computed hereafter. %\S\,\ref{explicit}.
To better understand the behaviour of solutions with singularties, it is interesting to study the unilateral gradient flow for the functional \eqref{e.F1D} in more detail.
Let us start considering the sub-interval $(0,1)$ and the following parabolic problem 
\begin{equation} \label{e.num}
	\begin{cases}
		\dot{u} (t,x)  = u'' (t,x) + f (x) \quad \text{ in } (0,T) \times (0,1)  \\ 
		u(t,0) = 1 , \ u(t,1)=0 \quad \text{ for } t \in (0,T)  \\
		u(0,x) = u_0 \quad \text{ for }  x \in ( 0, 1) .
	\end{cases}
\end{equation}
%where, by abuse of notation, we still denote by $f$ its restriction to $(0,1)$. 
By classical results, see e.g.~\cite[Theorem 5 (ii) \S7.1]{Evans} there exists a unique solution $u_r$ which belongs to $ L^\infty (0,T; H^2(0,1)) \cap H^1 (0,T; H^1(0,1)) $.

\separe

\begin{lemma} Let $u_r$ be the solution of the above parabolic problem. Then $\dot{u}_r (t) \ge 0$, and thus %it solves also 
$$
	\begin{cases}
		\dot{u}_r (t,x)  = [ u''_r (t,x) + f (x) ]_+ \quad \text{ in } (0,T) \times (0,1)  \\ 
		u_r (t,0) = 1 , \ u_r (t,1)=0 \quad \text{ for } t \in (0,T)  \\
		u_r (0,x) = u_0 \quad \text{ for }  x \in ( 0, 1) .
	\end{cases}
$$
Moreover $u'_r (t, 0) < 0$ for a.e.~$t \in (0,T)$.
%and $u_r (t) \le u_{min}$ where $u_{min}$ is the minimizer of the energy 
%$$	\F_r ( u (t) ) = \tfrac12 \int_{(0,L)} | \nabla u(t) |^2 - f \,u(t) \, dx   $$
%in . Moreover, $Du_r (t, 0^+) < 0$.
\end{lemma}

\proof First, let us prove that $\dot{u}_r \ge 0$. For convenience, we introduce the set $\mathcal{U}_r = \{ u \in H^1(0,1) : \, u(0)=1 ,\, u(1)=0 \}$ and the energy $\F_r : \mathcal{U}_r \to \R$ given by 
$$
	\F_r ( u ) = \tfrac12 \int_{(0,1)} | u' |^2 \, dx - \int_{(0,1)} f u \, dx  .
$$
%i.e.~the ``restriction'' of the energy $\F$ to the right sub-interval $(0,1)$.
%It is well known that $u_r$ can be found by an implicit (uncostrained) Euler scheme. 
Consider again an implicit (uncostrained) Euler scheme. Let $\tau > 0$ such that $T/\tau$ is integer and let $t_{k} = k \tau$ for $k =0,...,T/\tau$. Given $u_0$ we define by induction
\begin{gather} \label{e.scheme1r}
	u_{k+1} \in \argmin \{ \F_r ( u  ) + \tfrac{1}{2\tau} \, \| u - u_{k} \|^2_{L^2} :  u \in \mathcal{U}_r \}  .
\end{gather}
We will prove, by induction, that $u_{k+1} \ge u_k$ for every index $k$. 
Let us see that $u_1 \ge u_0$. Clearly, being $u_0$ affine, 
$$
	\int_{(0,1)} u'_0 \, \phi' \, dx = 0 , \quad \text{for every $\phi \in H^1_0 (0,1)$.}
$$
 Consider the auxiliary function $u_* = u_0 + | u_1 - u_0|$
and note that 
$$
	u_* \ge \max \{ u_0 , u_1 \} , \qquad | u_* - u_0 | = | u_1 - u_0 | ,
	\qquad
	| | u_1 - u_0 |' | = | ( u_1 - u_0)'| .
$$
We claim that $\F_r(u_*) \le \F_r(u_1)$. Indeed, since $(u_1-u_0)$ and $| u_1 - u_0|$ belongs to $H^1_0(0,1)$ we can write 
\begin{align*}
	\F_r ( u_1) 
	& = \F_r ( u_0 + (u_1 - u_0)) =  \int_{(0,L)} \tfrac12 | u'_0 | ^2 + \tfrac12 | ( u_1 - u_0)'  |^2 + u'_0 ( u_1 - u_0)' - f u_1 \, dx \\
	& = \int_{(0,L)} \tfrac12 | u'_0 | ^2 + \tfrac12 | ( u_1 - u_0)'  |^2  - f  u_1 \, dx  
\end{align*}
and 
\begin{align*}
	\F_r ( u_*) 
	& = \F_r ( u_0 + |u_1 - u_0| ) = \int_{(0,L)} \tfrac12 | u'_0 | ^2 + \tfrac12 | | u_1 - u_0|'  |^2 + u'_0   | u_1 - u_0 |' - f u_* \, dx \\
	& = \int_{(0,L)} \tfrac12 |   u'_0 | ^2 + \tfrac12 |   ( u_1 - u_0)'  |^2  - f  u_* \, dx .
\end{align*}
Since $u_* \ge u_1$ and $f=u_0 > 0$ we have $\F_r (u_*) \le \F_r (u_1)$. From the latter inequality it follows that
$$ \F_r( u_* ) + \tfrac{1}{2\tau} \, \| u_* - u_0 \|^2_{L^2} \le \F_r ( u_1)  + \tfrac{1}{2\tau} \, \| u_1 - u_0 \|^2_{L^2} $$ 
and then, by uniqueness of the minimizer in \eqref{e.scheme1r}, that $u_1 = u_* \ge \max \{ u_0 , u_1 \} $, i.e. $u_1 \ge u_0$.

Next, let us see that $u_{k+1} \ge u_k$ for $k \ge 1$. In this case, the Euler-Lagrange equation for $u_{k}$ reads 
% $u_k \in \argmin  \{ \F_r ( u  ) + \tfrac{1}{2\tau} \, \| u - u_{k-1} \|^2_{L^2} :  u \in \mathcal{U}_r \} $ and thus 
$$
	\int_{(0,1)}   u'_k    \phi' - f \phi + \tfrac{1}{\tau}( u_k - u_{k-1} ) \phi \, dx = 0 
	\qquad \text{for every $\phi \in H^1_0 (0,1)$.}
$$
As in the case $k=0$, it is enough to check that $\F_r ( u_*) \le \F( u_{k+1})$ for $u_* = u_k + | u_{k+1} - u_k |$.
Choosing $\phi = u_{k+1} - u_k $ and $\phi = | u_{k+1} - u_k |$ in the Euler-Lagrange equation yields, respectively,  \begin{align*}
	\F_r ( u_{k+1}) 
	& = \F_r ( u_k + (u_{k+1} - u_k))  \\
	& = \int_{(0,1)} \tfrac12 |   u'_k | ^2 + \tfrac12 |   ( u_{k+1} - u_k)'  |^2 +   u'_k   ( u_{k+1} - u_k)' - f u_k - f ( u_{k+1} - u_k)  \, dx \\
	& = \int_{(0,1)} \tfrac12 |   u'_k | ^2 + \tfrac12 |   ( u_{k+1} - u_k)'  |^2  - f  u_{k} \, dx - \tfrac{1}{\tau}  \int_{(0,1)} (u_k - u_{k-1} ) ( u_{k+1} - u_k)  \, dx 
\end{align*}
and
\begin{align*}
	\F_r ( u_*) 
	& = \F_r ( u_k + |u_{k+1} - u_k| ) \\
	& = \int_{(0,1)} \tfrac12 |   u'_k | ^2 + \tfrac12 |   | u_{k+1} - u_k|'  |^2 +   u'_k   | u_{k+1} - u_k |' - f u_k - f | u_{k+1} - u_k | \, dx \\
	& = \int_{(0,1)} \tfrac12 |   u'_k | ^2 + \tfrac12 |   ( u_{k+1} - u_k)'  |^2  - f  u_{k} \, dx - \tfrac{1}{\tau}  \int_{(0,1)} (u_k - u_{k-1} ) | u_{k+1} - u_k | \, dx .
\end{align*}
By induction $u_k \ge u_{k-1}$ and thus $(u_k - u_{k-1} ) ( u_{k+1} - u_k) \le (u_k - u_{k-1} ) | u_{k+1} - u_k |$. Hence $\F_r ( u_*) \le \F( u_{k+1})$.

\separe

It is well known that up to subsequences the piecewise affine interpolant $u_\tau$ converges weakly in $H^1(0,T; L^2(0,1))$ to the unique solution $u_r$ of \eqref{e.num}. Therefore, $u_r$ is monotone non-decreasing in time.

\separe

It is simple to check that the minimizer of the energy $\F_r$ in $\mathcal{U}_r$ is the function
$$
	u_{min} (x) = - \tfrac16 (1 - x)^3 + \tfrac76 (L-x) .
$$
Since $u_0 < u_{min}$, the comparison principle for \eqref{e.num} implies $u_0 \le u(t) \le u_{min}$ in $(0,T)$. Moreover, $u_0 ( t, 0) = u_r (t , 0) = u_{min} (t,0) = 1$, hence $ u_r'(t, 0) \le u'_{min}(t,0) = - 2/3$. \qed 

\begin{proposition} The function $u$ defined by 
$$
	u(t,x) = \begin{cases} u_r (t, x)  & x \in (0,1) \\  u_r (t, - x)  & x \in (-1,0) \end{cases}
$$
is the unilateral gradient flow for the functional $\F$ defined in \eqref{e.F1D}.
\end{proposition}

\proof For convenience, denote $u_l (t,x) = u_r (t,-x)$ and note that, by previous lemma, it holds $\dot{u}_{l} (t) =  [ u''_{l} (t) + f ]_+$.
 % and accordingly denote
%$$
%	\F_l ( u ) = \tfrac12 \int_{(-1,0)} | u' |^2 \, dx - \int_{(-1,0)} f u \, dx   . 
%$$
Moreover, in terms of $u_l$ and $u_r$ the derivatives of $u$ reads %in the sense of distributions we have 
$$   \dot{u} (t) =  \dot{u}_l (t) \, \mathcal{L}_{|(-1,0)} + \dot{u}_r (t) \, \mathcal{L}_{|(0,1)} ,
	\qquad
	u'' (t) = u''_l (t) \, \mathcal{L}_{|(-1,0)} + u''_r (t) \, \mathcal{L}_{|(0,1)} + [ u'_r ( t, 0) - u'_l (t,0) ] \, \delta_0  ,   $$
where $\mathcal{L}$ and $\delta_0$ denote respectively Lebesgue measure and Dirac delta. 
We remark that $\dot{u} ( t , 0 ) = 0$. Note that $[ u'_r ( t, 0) - u'_l (t,0) ] = 2  u'_r (t,0) < 0$, hence 
\begin{align*}
     [u''(t) + f]_+ %& = [ u''_l (t) \, \mathcal{L}_{|(-1,0)} + u''_r (t) \, \mathcal{L}_{|(0,1)} + f \, \mathcal{L}_{|(-1,1)} ]_+ \\[1mm]
	& = [ u''_l (t) +  f ]_+ \, \mathcal{L}_{|(-1,0)}   + [ u''_r (t) +  f ]_+ \, \mathcal{L}_{|(0,1)} .
\end{align*}
In particular $\dot{u} (t) = [u''(t) + f]_+$. % in $L^2(-1,1)$.

By the regularity of $f$ we can characterize the unilateral gradient flow for $\F$ as in Proposition \ref{ACsol}, i.e.
$$
	\dot{\F}( u(t)) \le - \tfrac12 | \dot{u}(t) |^2_{L^2_+}  - \tfrac12 |\partial \F|^2_{L^2_+} ( u(t) ) = - \tfrac12 \| \dot{u}(t) \|^2_{L^2} - \tfrac12 \| [ u'' (t) + f]_+ \|^2_{L^2} .
$$
Being $u \in H^1(0,T; H^1_0(-1,1))$ the chain rule and the fact that $\dot{u}(t,0)=0$ yield
\begin{align*}
	\dot{\F}( u(t)) & = d\F (u(t) ) [\dot{u}(t)] %= d \F_ l ( u_l (t) ) [ \dot{u}_l (t)] + d \F_r ( u_r (t) ) [ \dot{u}_r (t)] \\ & 
= - (u''(t) + f ,\dot{u} (t)) \\ 
& = - \langle u''_l (t) + f  , \dot{u}_l (t)  \rangle - \langle u''_r (t) + f  , \dot{u}_r (t)  \rangle 
	- ( [ u'_r ( t, 0) - u'_l (t,0) ] \, \delta_0 ,  \dot{u} (t,0) ) \\
	& = - \| \dot{u}_l (t) \|^2_{L^2} - \| \dot{u}_r (t) \|^2_{L^2} = - \| \dot{u} (t) \|^2_{L^2} = - \| [u''(t) + f]_+ \|^2_{L^2} ,
\end{align*}
which concludes the proof. \qed

\separe

\section{Solutions for \boldmath{$f \in AC (0,T;L^2)$}}

%\note{It could be interesting to prove time regularity, for instance $u \in AC_{loc}(0,T; H^1_0)$, without using the discrete scheme ... maybe with some integral form of Gronwall lemma}

%Assume that the datum $f \in AC (0,T;L^2). % is of (pointwise) bounded variations, i.e.~that
%$$
%	\| f \|_{BV(0,T;L^2)} = \sup \Big\{   \sum_{i=0}^I  \| f( t_{i+1}) - f( t_i)  \|  :  0 \le t_0 \le ... \le t_{I-1} \le t_I \le T \Big\} < +\infty \,.
%$$

In this section we consider the case in which $f \in AC(0,T;L^2)$ and we will prove the assertions contained in Proposition \ref{ACsol}. 

\subsection{Characterization by power balance}

If $f \in AC(0,T;L^2)$ then the map $t \mapsto \langle f(t) , u(t) \rangle$ is absolutely continuous in $(0,T)$ and 
for a.e.~$t \in (0,T)$ it holds
$$
	\tfrac{d}{dt} \langle f(t) , u(t) \rangle =  \langle \dot{f}(t) , u(t) \rangle  +  \langle f(t) , \dot{u}(t) \rangle  .
$$
We already know, by Theorem \ref{t.exist}, that the stored energy $ t \mapsto \E ( u(t))$ is absolutely continuous and that for a.e.~$t \in (0,T)$ we have
$$
	\dot{\E} ( u (t) ) =  - \tfrac12 | \dot{u} (t) |^2_{L^2_+} - \tfrac12 | \partial \F |^2_{L^2_+} ( t, u (t) ) \, dt + \langle f (t) , \dot{u} (t) \rangle .
$$
Therefore $t \mapsto \F ( t, u(t) ) = \E (u(t)) - \langle f(t) , u(t) \rangle$ is absolutely continuous and for a.e.~$t \in (0,T)$
$$
	\dot{\F} ( t, u(t) ) = - \tfrac12 | \dot{u} (t) |^2_{L^2_+} - \tfrac12 | \partial \F |^2_{L^2_+} ( t, u (t) ) \, dt - \langle \dot{f}(t) , u(t) \rangle ,
$$
which gives \eqref{e.maxsl}.

Conversly, if $t \mapsto \F (t, u(t))$ is absolutely continuous and \eqref{e.maxsl} holds, integration in time easily leads to show that $u$ is a unilateral gradient flow in the sense of Definition \ref{d.curve}.

%\subsection{Compactness}

%By \eqref{e.evol-ineq} and Cauchy inequality for $0 <  t_2 \le T$ we get 
%\begin{equation}  \label{e.einz}
%       \F  ( u  (t_2) ) \le \F ( u_0 ) \, -  \, \tfrac12 \int_{0}^{t^*} \| \dot{u} (t) \|_{L^2}^2 + | \partial_+ \F |^2 ( u (t) ) \, dt \le 
% \F ( u_0 ) \, -  \, \int_{0}^{t^*}  | \partial_+ \F | ( u (t) ) \| u ' (t) \|_{L^2} \, dt .
%\end{equation}
%
%\separe
%
%By \eqref{e.upgrad}, applied to $u$ we get 
%$$
%	 \F ( u (t') ) \le \F ( u (t_2) ) + \int_{t'}^{t_2} | \partial_+ \F | ( u (t) ) \| u ' (t) \|_{L^2}  \, dt 
%$$
%for every $ 0 < t' < t_2 < T$;  taking the liminf for $t' \to 0^+$ we obtain
%\begin{equation}  \label{e.zwei}
%	\F ( u_0 ) \le  \F ( u (t_2) )  + \int_{0}^{t_2}  | \partial_+ \F | ( u (t) ) \| u ' (t) \|_{L^2}  \, dt .
%\end{equation}
%Joining  \eqref{e.einz} and \eqref{e.zwei} provides the energy identity
%$$
%          \F ( u_0 ) = \F ( u (t_2) )  + \int_{0}^{t_2}  | \partial_+ \F | ( u (t) ) \| u ' (t) \|_{L^2}  \, dt .
%$$

\subsection[{Characterization by differential inclusions in $H^1_0$}]{Characterization by differential inclusions in \boldmath{$H^1_0$}}

The results of this section are essentially an adaption of those contained in \cite{GianazzaSavare_RANSXLMM94}. For sake of completness, we provide some short alternative proofs.  

\begin{proposition} \label{p.Wloc} Let $u_0 \in H^1_0$ and $f \in AC (0,T ; L^2)$ then the unilateral gradient flow $u$ belongs to $W^{1,\infty}_{loc} ( 0 , T ; L^2 ) \cap W^{1,2}_{loc} (0,T ; H^1_0)$.
\end{proposition}

\proof We will employ the discrete evolutions $u_n$ obtained  by the implicit Euler scheme of \S\ref{3.1}. We will show that given $0 < T' < T$ the sequence $u_n$ is bounded in $W^{1,\infty} ( T' , T ; L^2 ) \cap W^{1,2} ( T' , T  ; H^1_0)$, from which it follows that $u \in W^{1,\infty}_{loc} ( 0 , T ; L^2 ) \cap W^{1,2}_{loc} (0,T ; H^1_0)$.

\separe 

Let us denote $\dot{u}_{n,k} = ( u_{n,k} - u_{n,k-1} ) / \tau_n$.  For $k\ge 0$ by \eqref{e.dvAA} we have
\begin{equation} \label{e.2k+1}
	 \| \dot{u}_{n,k+1} \|^2_{L^2} + d \F_n ( t_{n,k+1} , u_{n,k+1} )  [ \dot{u}_{n,k+1} ]  = 0  \,.
\end{equation}
For $k \ge 1$ %by minimality of $u_{n,k}$ we have 
%$$
%	d \F_n ( t_{n,k}, u_{n,k} )  [ v - u_{n,k} ] + \langle \dot{u}_{n,k} , v - u_{n,k} \rangle \ge 0  \quad \text{for every $v \in H^1_0$ with $v \ge u_{n,k-1}$.}
%$$
by \eqref{e.dPDE} we get 
$$
      d \F_n ( t_{n,k}, u_{n,k} )  [ \dot{u}_{n,k+1} ] + \langle \dot{u}_{n,k} , \dot{u}_{n,k+1} \rangle \ge 0 .
$$
Hence, for $k \ge 1$ we obtain
\begin{align}\label{e.three*}
	d \F_n ( t_{n,k}, u_{n,k} )  [ \dot{u}_{n,k+1} ] - d \F_n ( t_{n,k+1}, u_{n,k+1} )  [ \dot{u}_{n,k+1} ] & \ge -  \langle \dot{u}_{n,k} , \dot{u}_{n,k+1} \rangle +  \langle \dot{u}_{n,k+1} , \dot{u}_{n,k+1} \rangle . %\nonumber \\ & \ge \tfrac12  \| \dot{u}_{n,k+1} \|^2_{L^2}  - \tfrac12  \| \dot{u}_{n,k} \|^2_{L^2} 
\end{align}
Let us write explicitely the left hand side as 
\begin{align*}
	%d \F_n ( t_{n,k-1}, u_{n,k} )  [ \dot{u}_{n,k+1} ] - d \F_n ( t_{n,k-1}, u_{n,k} )  [ v - u_{n,k} ] = 
	a (  u_{n,k} , \dot{u}_{n,k+1} ) - \langle f_{n,k} , \dot{u}_{n,k+1} \rangle & - a ( u_{n,k+1} , \dot{u}_{n,k+1})  + \langle f_{n,k+1} , \dot{u}_{n,k+1} \rangle 
	= \\
	& = 
	- \tau_n a ( \dot{u}_{n,k+1} , \dot{u}_{n,k+1} ) +  \langle f_{n,k+1} - f_{n,k} , \dot{u}_{n,k+1} \rangle .
\end{align*}
A simple algebraic calculation gives $-  \langle \dot{u}_{n,k} , \dot{u}_{n,k+1} \rangle +  \langle \dot{u}_{n,k+1} , \dot{u}_{n,k+1} \rangle  \ge \tfrac12  \| \dot{u}_{n,k+1} \|^2_{L^2}  - \tfrac12  \| \dot{u}_{n,k} \|^2_{L^2}$. Hence by coercivity \eqref{e.three*} reads, for $k \ge 1$,
\begin{align}\label{e.three**}
   - c \tau_n \| \dot{u}_{n,k+1} \|^2_{H^1_0} + \| f_{n,k+1} - f_{n,k}\|_{L^2} \, \| \dot{u}_{n,k+1} \|_{L^2}  \ge \tfrac12  \| \dot{u}_{n,k+1} \|^2_{L^2}  - \tfrac12  \| \dot{u}_{n,k} \|^2_{L^2} .
\end{align}
Neglecting the $H^1_0$-norm and denoting $a_k = \| \dot{u}_{n,k} \|_{L^2}$ we obtain the discrete inequality $a^2_{k+1} - a^2_k \le b_k a_{k+1}$ where $b_k = 2 \| f_{n,k+1} - f_{n,k}\|_{L^2}$. By an elementary algebraic calculation we get $a_{k+1} \le a_k + b_k$. Then for every $k_0 < k$ we can write %$a_{k} \le a_{k_0} + \sum_{j=k_0}^{k-1} b_j$. Note that 
$$
     a_{k} \le a_{k_0} + \sum_{j=k_0}^{k-1} b_j \le a_{k_0} + 2 \,\| f \|_{AC (0,T; L^2)} \le a_{k_0} + C \,.
$$
Let $k'_n$ such that $ \tau_n (k'_n - 1) < T' \le \tau_n k'_n $. 
% Let $k_0$ and $k'_n$ such that $1 \le k_0 \le k'_n < k$. 
Given $k > k'_n$ the estimate $a_k \le a_{k_0} + C$ holds for every index $k_0$ such that $1 \le k_0 \le k'_n$ with $C=2 \,\| f \|_{AC (0,T; L^2)}$.
Taking the sum of $a_k \le a_{k_0} + C$  for $k_0=1,...,k'_n$  (and $k$ fixed) and deviding by $k'_n$ yields
%$$
%	\sum_{k_0=1}^{k'_n}  a_k \le \sum_{k_0=1}^{k'_n} a_{k_0} +  \sum_{k_0=1}^{k'_n} 2 \,\| f \|_{AC(0,T; L^2)} \,.
%$$
%which reads
$$
	a_k \le \frac{1}{k'_n} \sum_{k_0=1}^{k'_n} ( a_{k_0} + C ) \le \frac{1}{T'} \Big( \sum_{k_0=1}^{k'_n} \tau_n a_{k_0}  \Big) + C .
$$
Hence, for every $k > k'_n$ we have 
$$
	\| \dot{u}_{n,k} \|_{L^2} \le \frac{1}{T'}  \int_{0}^{T} \| \dot{u}_{n} (t) \|_{L^2} \, dt + C  %\le C \F_nrac{1}{T'^{1/2}} \| u_{n} \|_{H^1(0,T; L^2)}  + C 
\le C (T') ,
$$
where $C (T') >0$ is independent of $n$ because the sequence $u_n$ is bounded in $H^1 ( 0 , T ; L^2)$. (Note that $C(T')$ diverges as $T' \to 0^+$). Taking the supremum with respect to $k > k'_n$ it follows that the sequence $\dot{u}_n$ is bounded in $L^\infty( T'', T ; L^2)$ for every $T'' >T' > 0$.

\separe

Let us go back to \eqref{e.three**}. For every $k > k'_n$, now we can write %written in form 
$$
	  \tfrac12  \| \dot{u}_{n,k} \|^2_{L^2}  - \tfrac12  \| \dot{u}_{n,k+1} \|^2_{L^2} + C (T') \, \| f_{n,k} - f_{n,k-1} \|_{L^2} \ge c \tau_n \| \dot{u}_{n,k+1} \|^2_{H^1_0}  .
$$
Taking the sum for $k > k'_n$ we get
$$
	c \int_{T''}^{T} \| \dot{u}_n \|^2_{H^1_0} \, dt \le c \int_{\tau_n k'_n}^{T} \| \dot{u}_n \|^2_{H^1_0} \, dt \le C (T')  \| f \|_{AC(0,T;L^2)}  + \tfrac12 \| \dot{u}_{n,k'_n+1} \|^2_{L^2} \le C'(T') ,
$$
for every $T''> T' >0$, which concludes the proof. \qed

\separe

Let $\Phi: H^1_0 \to [0,+\infty]$ be defined by  $\Phi (u) = \tfrac12 \| u \|_{L^2}^2 + I_{ \{u \,\ge\, 0\} }$ (where $I$ denotes the indicator function) and let $\partial \Phi (u) \subset H^{-1}$ be its subdifferential. Let $\Phi^*$ be its Legendre transform, i.e. $  \Phi^*(\xi) = \sup \{ ( \xi , v ) - \Phi (v) : v \in H^1_0 \}$ for $\xi \in H^{-1}$. By Lemma \ref{l.measure} it is easy to check that the domain of $\Phi^*$ is $D \Phi^* = \{ \xi \in H^{-1} : \xi \in \mathcal{M}_{loc} \,,\, \xi_+ \in L^2 \}$ and then that $\Phi^* (\xi) = \tfrac12 \| \xi_+ \|_{L^2}$ for $\xi \in D\Phi^*$.

\begin{proposition} \label{p.=H10} Let $f \in AC (0,T;L^2)$ and $u \in W^{1,2}_{loc} ( 0,T ; H^1_0)$. Then $u$ is a unilateral gradient flow for $\F$ with intial condition $u_0$ if and only if it solves the differential inclusion \begin{equation*}%\label{e.pdeH1}
	\begin{cases}
		\partial \Phi ( \dot{u} (t) ) \ni  A u (t) + f (t)  & \text{ in $H^{-1}$  for a.e.~$t \in (0,T)$} \\
		u (0) = u_0 .
	\end{cases}
\end{equation*}
%where $\Phi: H^1_0 \to [0,+\infty]$ is defined by  $\Phi (u) = \tfrac12 \| u \|_{L^2}^2 + I_{ \{u \,\ge\, 0\} } $ (where $I$ denotes the indicator function) and $\partial \Phi (u) \subset H^{-1}$ denotes its subdifferential. 
\end{proposition}

\proof
Since $\dot{u} (t) \in H^1_0$,  for a.e.~$t \in (0,T)$, by classical results in convex analysis we have
\begin{eqnarray}
% & \dot{u} (t) \in \partial \Phi^* ( A u (t) + f (t) ) &  \nonumber \\
%	& \Updownarrow & \nonumber \\
 & A u (t) + f (t) \in \partial \Phi ( \dot{u}(t) ) & \nonumber \\
	& \Updownarrow & \nonumber \\
 & \Phi (\dot{u}(t) ) + \Phi^* ( A u (t) + f (t) ) = ( A u (t) + f (t) , \dot{u}(t) ) & \nonumber \\
    & \Updownarrow & \nonumber \\
& \tfrac12 \| \dot{u}(t) \|^2_{L^2} + \tfrac12 \| [ A u(t) + f(t) ]_+ \|^2_{L^2} = - d \F ( t, u(t)) [ \dot{u}(t) ] & \nonumber \\
 & \Updownarrow & \nonumber \\
 &   \tfrac12 | \dot{u}(t) |^2_{L^2_+} + \tfrac12 | \partial \F|^2_{L^2_+} ( t, u(t) ) = - d \F ( t, u(t)) [ \dot{u}(t) ] . &  \label{e.2-ei}
\end{eqnarray}
The fact that $\partial \Phi ( \dot{u} (t) ) \neq \emptyset$ implicitely says that $ \dot{u}(t) \in D \Phi$ and hence
$ \dot{u} (t) \ge 0$. 
Moreover, if $u \in W^{1,2}_{loc} (0,T; H^1_0)$ with $\dot{u} \ge 0$ then by \eqref{e.2-ei} for a.e.~$t \in (0,T)$ we can write 
$$
	\dot{\F} ( t, u(t) ) = d\F ( t, u(t) ) [ \dot{u} (t) ] + \partial_t \F ( t, u(t) ) \le	
	- \tfrac12 | \dot{u} |^2_{L^2_+}  - \tfrac12 | \partial \F |^2_{L^2_+}  ( t, u(t)) - \langle \dot{f} (t) , u(t) \rangle ,
$$	
which is indeed the characterization \eqref{e.maxsl} of unilateral gradient flows. Conversly, if the energy identity \eqref{e.maxslen} holds then for a.e.~$t \in (0,T)$ we have 
$$
	\dot{\F} ( t, u(t) ) = d\F ( t, u(t) ) [ \dot{u} (t) ] + \partial_t \F ( t, u(t) ) = 	
	- \tfrac12 | \dot{u}(t) |^2_{L^2_+} - \tfrac12 | \partial \F|^2_{L^2_+} ( t, u(t) )  - \langle \dot{f} (t) , u(t) \rangle 
$$	
and then \eqref{e.2-ei} holds.

\section{A characterization for \boldmath{$f$} independent of time \label{s.auto}}

In this section we will prove Proposition \ref{p.AE}.

\medskip \noindent {\bf Step I.} %We employ the discrete scheme of \S\ref{3.1} and 
We claim that for every $v \in H^1_0$ with $u_0 \le v \le u_{n,k+1}$
we have
\begin{equation} \label{e.backintime}
	\F ( u_{n,k+1} ) + \tfrac{1}{2 \tau_n} \| u_{n,k+1} - u_{n,k} \|^2_{L^2} \le \F ( v ) + \tfrac{1}{2 \tau_n} \| v - u_{n,k} \|^2_{L^2} .
\end{equation}
Given $0 \le m \le k+1$ let $\Omega_{n,m} = \{ v \ge u_{n,m} \}$. Since $u_0 \le v \le u_{n,m}$ we have $\Omega_{n,0} = \Omega$ and $\Omega_{n,m+1} = \emptyset$, while the monotonicity of $u_{n,m}$ w.r.t.~$m$ implies that the sets $\Omega_{n,m}$ are monotone non-increasing w.r.t.~$m$ and $\{ \Omega_{n,m} \setminus \Omega_{n,m+1} : \text{ for } m =0 , ...,k \}$ is a disjoint partition of $\Omega$. 
For any measurable subset $O$ of $\Omega$ or we employ the notation 
$$
	\F_{O} ( w ) = \int_{O} \nabla w \cdot B \,\nabla w + b \, w^2 - f w \,  dx .
$$
For each index $0 \le m \le k$ let us define $v_{n,m} = \min \{ \max \{  v , u_{n,m} \} , u_{n,m+1} \} $, i.e., 
$$
	v_{n,m} = 
	\begin{cases}  
		u_{n,m+1} & \text{in $\Omega_{n,m+1}$} \\
		v & \text{in $\Omega_{n,m} \setminus \Omega_{n,m+1}$} \\
		u_{n,m} & \text{in $\Omega_{n,m}^c$} .	
	\end{cases}
$$
By minimality we have 
$$
	\F ( u_{n,m+1} ) + \tfrac{1}{2 \tau_n} \| u_{n,m+1} - u_{n,m} \|^2_{L^2} \le \F ( v_{n,m} ) + \tfrac{1}{2 \tau_n} \| v_{n,m} - u_{n,m} \|^2_{L^2} .
$$
Since $u_{n,m} \le v_{n,m} \le u_{n,m+1}$ we have $\| u_{n,m+1} - u_{n,m} \|^2_{L^2} \ge \| v_{n,m} - u_{n,m} \|^2_{L^2}$ and thus $\F ( u_{n,m+1} ) \le \F ( v_{n,m} ) $; the latter inequality can be written as 
$$
   \F_{\Omega_{n,m+1}} ( u_{n,m+1} ) + \F_{\Omega_{n,m+1}^c} ( u_{n,m+1} ) \le \F_{\Omega_{n,m+1}} ( v_{n,m} ) + \F_{\Omega_{n,m} \setminus \Omega_{n,m+1}} ( v_{n,m} ) + \F_{\Omega_{n,m}^c} ( v_{n,m} ) 
$$
which, by definition of $v_{n,m}$, yields 
\begin{equation} \label{e.indu}
		\F_{\Omega_{n,m+1}^c} ( u_{n,m+1} ) \le \F_{\Omega_{n,m} \setminus \Omega_{n,m+1}} ( v ) + \F_{\Omega_{n,m}^c} ( u_{n,m} ) .
\end{equation}
Next, being $\Omega_{n,k+1} = \emptyset$ by minimality we can write 
\begin{align*}
	\F ( u_{n,k+1} ) + \tfrac{1}{2 \tau_n} \| u_{n,k+1} - u_{n,k} \|^2_{L^2} 
		& \le \F ( v_{n,k} ) + \tfrac{1}{2 \tau_n} \| v_{n,k} - u_{n,k} \|^2_{L^2} \\
		& \le 	\F_{\Omega_{n,k}} ( v ) + \F_{\Omega_{n,k}^c} ( u_{n,k} ) + 
				\tfrac{1}{2 \tau_n} \int_{\Omega_{n,k}} | v - u_{n,k} |^2 \, dx .
\end{align*}
Remember that $\Omega_{n,0}^c = \emptyset$ and that $\{ \Omega_{n,m} \setminus \Omega_{n,m+1} : \text{ for } m =0 , ...,k \}$ is a disjoint partition of $\Omega$. Hence, using iteratively \eqref{e.indu} we get 
\begin{align*}
	\F ( u_{n,k+1} ) + \tfrac{1}{2 \tau_n} \| u_{n,k+1} - u_{n,k} \|^2_{L^2} 
		& \le \bigg(  \sum_{k=0}^m \F_{\Omega_{n,k} \setminus \Omega_{n,k+1}} ( v ) \bigg) + 
				\tfrac{1}{2 \tau_n} \int_{\Omega_{n,k}} | v - u_{n,k} |^2 \, dx \\
		& \le \F ( v ) + \tfrac{1}{2 \tau_n} \| v - u_{n,k} \|^2_{L^2} ,
\end{align*}
which proves the claim \eqref{e.backintime}. 

\medskip \noindent {\bf Step II.} Thanks to \eqref{e.backintime}, the incremental problem \eqref{e.scheme1} can be replaced by the fixed obstacle problem
\begin{equation} \label{e.minu_0}
	u_{n,k+1} \in \argmin \big\{ \F ( u  ) + I_+ (u-u_0) + \tfrac{1}{2\tau_n} \, \| u - u_{n,k} \|^2_{L^2} :  u \in H^1_0 \big\}  .
\end{equation}
Indeed, let $v \ge u_0$. Denote $v^+ = \max \{ v , u_{n,k+1} \}$ and $v^- = \min \{ v , u_{n,k+1} \}$, accordingly let $\Omega^+ = \{ v \ge u_{n,k+1} \}$ and $\Omega^- = \Omega \setminus \Omega^+$. By minimality of $u_{n,k+1}$ we can write 
$$
	\F ( u_{n,k+1} ) + \tfrac{1}{2 \tau_n} \| u_{n,k+1} - u_{n,k} \|^2_{L^2} 
		\le \F ( v^+ ) + \tfrac{1}{2 \tau_n} \| v^+ - u_{n,k} \|^2_{L^2} 
$$
and then, 
$$
	\F_{\Omega^+} ( u_{n,k+1} ) + \tfrac{1}{2 \tau_n} \| u_{n,k+1} - u_{n,k} \|^2_{L^2(\Omega^+)} 
		\le \F_{\Omega^+} ( v ) + \tfrac{1}{2 \tau_n} \| v- u_{n,k} \|^2_{L^2(\Omega^+)} . 
$$
Since $u_0 \le v^- \le u_{n,k+1}$ we can invoke \eqref{e.backintime} which gives
$$
	\F ( u_{n,k+1} ) + \tfrac{1}{2 \tau_n} \| u_{n,k+1} - u_{n,k} \|^2_{L^2} 
		\le \F ( v^- ) + \tfrac{1}{2 \tau_n} \| v^- - u_{n,k} \|^2_{L^2} 
$$
and then, 
$$
	\F_{\Omega^-} ( u_{n,k+1} ) + \tfrac{1}{2 \tau_n} \| u_{n,k+1} - u_{n,k} \|^2_{L^2(\Omega^-)} 
		\le \F_{\Omega^-} ( v ) + \tfrac{1}{2 \tau_n} \| v - u_{n,k} \|^2_{L^2(\Omega^-)} . 
$$
Taking the sum of estimates in $\Omega^\pm$ shows that \eqref{e.minu_0} holds. As a consequence, the limit evolution $u$ turns out to be the $L^2$-gradient flow for the functional $\tilde{\F} ( u ) = \F (u) + I_+ (u - u_0)$.

\medskip \noindent {\bf Step III.} 
Now, let us show that the (unilateral) gradient flow is a solution of the parabolic obstacle problem \eqref{e.AE}, i.e.,
\begin{equation} \label{e.AE*}
	\begin{cases} 
		\dot{u}(t) - A u(t) - f \ge 0 & \text{in $H^{-1}$ for a.e.~$t \in (0,T)$} \\
		( u(t) - u_0 ,  \dot{u}(t) - A u(t) - f ) = 0 & \text{for a.e.~$t \in (0,T)$} \\
		u(0) = u_0 , \quad u(t) \ge u_0 & \text{for a.e.~$t \in (0,T)$.} 
	\end{cases}
\end{equation}
%(In the next step we will prove that the solution of \eqref{e.AE*} is unique). 
By Proposition \ref{p.pde} we know that $\dot{u} (t) = [ A u(t) + f ]_+$ hence
$$
	\dot{u} (t) - (A u(t) + f) = \dot{u} (t) - [ A u(t) + f ]_+ + [ A u(t) + f ]_- = [ A u(t) + f ]_- \ge 0.
$$
In particular $( u(t) - u_0 ,  \dot{u}(t) - A u(t) - f ) \ge 0$. Therefore, re-writing the second line in \eqref{e.AE*}, it remains to show that for a.e.~$t \in (0,T)$ we have 
\begin{equation} \label{e.AEaux}
	\langle u(t) - u_0 ,  \dot{u}(t) \rangle + d \F ( u(t) ) [ u(t) - u_0 ] \le 0 . 
\end{equation}
The functional $\mathcal{J} ( w) = \F (w) + \tfrac{1}{2 \tau_n} \| w - u_{n,k} \|^2_{L^2}$ is convex, thus \eqref{e.backintime} implies that $d \mathcal{J} (u_{n,k+1}) [ u_0 - u_{n,k+1} ] \ge 0$, i.e.
$$
	d\F ( u_{n,k+1}) [u_0 - u_{n,k+1} ] + \tfrac{1}{\tau_n} \langle u_{n,k+1} - u_{n,k} , u_0 - u_{n,k+1} \rangle \ge 0 .
$$
In terms of the piece-wise affine interpolant $u_n$ and the piece-wise constant interpolant $u_n^\sharp$ (see \S\ref{3.1}), for every $t \in (t_{n,k} , t_{n,k+1})$ the previous inequality reads
$$
	d\F ( u_{n}^\sharp(t) ) [u_0 - u_{n}^\sharp(t) ] + \langle \dot{u}_{n} (t) , u_0 - u_{n}^\sharp(t) 
\rangle \ge 0 .
$$
Given $0< t^- < t^+ < T$, %let $k^\pm_n$ such that $t^\pm \in [t_{n,k^\pm_n} , t_{n,k^\pm_n+1} )$. 
we obtain
$$
	\int_{t^-}^{t^+} d\F ( u_{n}^\sharp(t) ) [u_0 - u_{n}^\sharp(t) ] + \langle \dot{u}_{n} (t) , u_0 - u_{n}^\sharp(t) \rangle \, dt \ge 0 .
$$
By Remark \ref{r.uniqueness} we know that $u_n \weakto u$ in $H^1(0,T;L^2)$, that $u_n^\sharp$ is bo that $u_n^\sharp$ is bounded in $L^\infty(0,T;H^1_0)$ and that $u_n^\sharp (t) \to u(t)$ strongly in $H^1_0$ for a.e.~$t$ in $(0,T)$unded in $L^\infty(0,T;H^1_0)$ and that $u_n^\sharp (t) \to u(t)$ strongly in $H^1_0$ for a.e.~$t$ in $(0,T)$. Therefore, we can pass to the limit in the previous inequality and get 
$$
	\int_{t^-}^{t^+} d\F ( u (t) ) [u_0 - u(t) ] + \langle \dot{u} (t) , u_0 - u(t) \rangle \, dt \ge 0 .
$$
By arbitrariness of $t^-$ and $t^+$ we conclude that \eqref{e.AEaux} holds a.e.~in $[0,T]$. 

\medskip \noindent {\bf Step IV.} We prove that there exists a unique solution of \eqref{e.AE*}. Assume, by contradiction, that $u_\I$ and $u_\II$ are solutions of \eqref{e.AE*} with $u_\I(t^*) \neq u_\II (t^*)$. We define $u_\natural = \tfrac12 ( u_\I + u_\II )$. By linearity, for a.e.~$t \in (0,T)$ we have 
$\dot{u}_\natural (t) - A u_\natural (t) - f \ge 0$ and 
%		( u(t) - u_0 ,  \dot{u}(t) + A u(t) + f ) = 0 & \text{} \\
$	u_\natural (t) \ge u_0$, 
and thus 
$$
	\int_0^{t^*}  ( u_\natural (t) - u_0 ,  \dot{u}_\natural (t) - A u_\natural (t) - f ) \, dt \ge 0 .
$$
On the other hand, by the strict convexity of the $L^2$-norm we can write % taking the integral in $(0,t^*)$ yields
\begin{align*}
	\int_0^{t^*} ( u_\natural (t) - u_0 ,  \dot{u}_\natural (t) ) \, dt & = \tfrac12 \| u_{\natural} (t^*) - u_0 \|^2_{L^2}  < \tfrac12 \big(  \tfrac12  \| u_\I (t^*) - u_0 \|^2_{L^2} + \tfrac12 \| u_\II (t^*) - u_0 \|^2_{L^2} \big) \\
	& = \tfrac12 \int_0^{t^*} ( u_\I (t) - u_0 ,  \dot{u}_\I (t) ) + ( u_\II (t) - u_0 ,  \dot{u}_\II (t) ) \, dt .
\end{align*}
Moreover, by convexity of the stored energy $ \E ( \cdot) = a ( \cdot, \cdot) $ we have 
\begin{align*}
	( u_\natural (t) - u_0 , - A u_\natural (t) - f )  & = ( u_\natural (t) , - A u_\natural (t) ) - (u_\natural (t) , f  ) + ( u_0 , A u _\natural (t) + f ) \\
	& = a ( u_\natural (t) , u_\natural (t) ) - (u_\natural (t) , f  ) + ( u_0 , A u _\natural (t) + f ) \\
	& \le \tfrac12 a ( u_\I (t) , u_\I (t) ) + \tfrac12 a ( u_\II (t) , u_\II (t) ) - \tfrac12  (u_\I (t) , f  ) - \tfrac12  (u_\II (t) , f  ) \,  + \\ 
& \phantom{\le.} \tfrac12 ( u_0 , A u_\I (t) + f ) + \tfrac12 ( u_0 , A u_\II (t) + f ) \\
	& = \tfrac12 ( u_\I (t) - u_0 , - A u_\I (t) - f )   + \tfrac12 ( u_\II (t) - u_0 , - A u_\II (t) - f )  
\end{align*}
Taking the integral in $(0,t^*)$ we obtain the contradiction
\begin{align*}
	\int_0^{t^*}  ( u_\natural (t) - u_0 , &  \dot{u}_\natural (t) - A u_\natural (t) - f ) \, dt < \\ & < \tfrac12 \int_0^{t^*}  ( u_\I (t) - u_0 ,  \dot{u}_\I (t) - A u_\I (t) - f ) +  ( u_\II (t) - u_0 ,  \dot{u}_\II (t) - A u_\II (t) - f ) \, dt = 0 ,
\end{align*}
which concludes the proof. \qed

\begin{remark} \normalfont \label{r.auto} If $f$ depends on time then, in general, unilateral gradient flows do not enjoy \eqref{e.AE}, as the following example shows. Let $\Omega= (0,1)$ and $u_0 (x) = x ( x-1)$. Define
$$
	u(t) = \begin{cases} 
		(1-t) u_0 & \text{if $t \in [0,1]$} \\ 
		0	& \text{if $t >1$,}
	\end{cases}
	\qquad \qquad 
	f (t ) = \begin{cases} 
		\dot{u} (t) - u'' (t) & \text{if $t \in [0,1]$} \\ 
		u_0	& \text{if $t >1$.}
	\end{cases}
$$
\end{remark}
It is easy to check that $u$ is monotone non-decreasing and that $\dot{u} (t) = u'' (t) + f(t)$ for $t \in (0,1)$. In particular $\| \dot{u} (t) \|_{L^2} = \| \dot{u} (t) \|_{L^2_+}$ and  $| \partial \F |_{L^2_+} (t, u(t)) = \| u'' (t) + f(t) \|_{L^2}$.  Thus, by the chain rule
$$
	\dot{\F} ( t, u(t)) =  d \F ( u(t) ) [ \dot{u} (t) ] - \partial_t \F ( t , u(t) ) = d \F ( u(t) ) [ \dot{u} (t) ] - \langle \dot{f} (t) , u(t) \rangle ,
$$
and 
$$
	d \F ( u(t) ) [ \dot{u} (t) ] = - \int_\Omega ( u''(t) + f(t) ) \, \dot{u}(t) \, dx = - \tfrac12 \| \dot{u} (t) \|^2_{L^2} - \tfrac12 \| u''(t) + f(t) \|^2_{L^2} .
$$
Hence $u$ is the unilateral gradient flow for $t \in [0,1]$. For $t > 1$ we have $\dot{\F} ( t , u(t)) = 0$ and $u (t) =0$, hence $\| \dot{u} (t) \|_{L^2} = 0$. In this case we have $| \partial \F |_{L^2_+} (t, u(t)) = \| [ u'' (t) + f(t) ]_+ \|_{L^2} = 0$. Hence $u$ is the unilateral gradient flow also for $t > 1$. On the other hand for $t >1$ we have 
$$
	( u (t) - u_0 , \dot{u} ( t) - u'' (t) - f ) = \int_{(0,1)} u^2_0 \, dx  \neq 0 ,
$$
thus \eqref{e.AE} does not hold. \qed

%\input{pdeLp}
%\include{Rn}
% !TeX root = article.tex

\appendix

\section{Metric settings \label{app.A}}

Having in mind the modern theory of gradient flows \cite{AmbrosioGigliSavare05} it is interesting to see if and how our framework fits into some sort of metric setting. For sake of simplicity, we assume that the force $f$ is constant. 

\medskip
\noindent {\bf A singular metric for \boldmath $L^2_+$.} Consider the ``singular'' metric
\begin{equation*}
    d_+ ( v , u) = \begin{cases}
	\| v - u  \|_{L^2} &  \text{ if $v \ge u$} \\
	+ \infty & \text{ otherwise.}
	\end{cases}
\end{equation*}
Accordingly, we will say that $u_m \to u$ when $d_+ ( u_m , u) \to 0$, i.e.~when $u_m \ge u$ and $u_m \to u$ in $L^2$. In this way, we can consider $L^2_+$ as the space $L^2$ ``endowed'' with $d_+$.

It is clear that if $u \in AC_{loc} (0,T; L^2_+)$ then $u$ is monotone non-decreasing and the metric derivative (cf.~\cite[Theorem 1.1.2]{AmbrosioGigliSavare05}) exists a.e.~in $(0,T)$ and reads
\begin{equation*}
	| \dot{u} |_{L^2_+} (t) = \lim_{h \to 0^+}  \frac{d_+ ( u(t+h) , u (t) ) }{h} = \| \dot{u} (t) \|_{L^2} = | \dot{ u} (t) |_{L^2_+} .
\end{equation*}
Moreover, the unilateral slope reads (cf.~\cite[Definition 1.2.4]{AmbrosioGigliSavare05})
\begin{equation*}
 | \partial \F |_{L^2_+} (u) = \limsup_{v \to u} \frac{[\F(v) - \F(u)]_-}{d_+ (v, u)} .
\end{equation*}

A key point in \cite{AmbrosioGigliSavare05} is the fact that slopes are upper gradients 
(cf.~\cite[Definition 1.2.1 and 1.2.2]{AmbrosioGigliSavare05}), however in our case we have the following (negative) result.

\begin{remark} The unilateral slope $| \partial \F |_{L^2_+}$ is neither a strong nor a weak upper gradient for $\F$. %as the following %example shows. %see \cite[Definition 1.2.1 and Definition 1.2.2]{AmbrosioGigliSavare05}. 
\end{remark}

The counter-example of \S\,\ref{non-uni} applies also here. Let $u_0 (x) = 1 - |x|$, $u(t) = (1+t) u_0$ and $f = u_0$. Clearly $\dot{u} = u_0$ and $u \in AC (0,T; L^2 (-1,1))$. Moreover, 
$$
	\F(u(t)) = \tfrac12 \| u (t) \|^2_{H^1_0} - \langle f , u(t) \rangle = \tfrac12 (1+t)^2 \| u_0 \|^2_{H^1_0} - (1+t) \| u_0 \|^2_{L^2} = \F (u_0) + (\tfrac12 t^2 + t) \| u_0 \|^2_{H^1_0} - t\| u_0 \|^2_{L^2} .
$$
Hence, 
$$
	\dot{\F} ( u(t)) = (t+1) \| u_0 \|^2_{H^1_0} - \| u _0 \|^2_{L^2} .
$$

Now, let us  heck that the inequality $| \dot{\F} ( u(t)) |  \le  | \partial \F |_{L^2_+} (u(t)) \, | \dot{u}(t) |_{L^2_+}$ fails for $t$ sufficiently large. For every $z \in H^1_0(-1,1)$ with $z \ge 0$ and $\| z \|_{L^2} \le 1$ we have 
$$- d\F(u(t)) [z] = ( u''(t) + f , z ) = ( - 2 (1+t) \delta_0 + u_0 , z ) \le \langle u_0 , z \rangle \le \| u_0  \|_{L^2} ,  $$
then $
	| \partial \F |_{L^2_+}  ( u(t) )  
		=  \sup  \big\{ \! - d\F(u(t)) [z] :  z \in H^1_0 , \, | z |_{L^2_+} \le 1 \big\}   \le \| u_0 \|_{L^2}$ and thus 
$$ | \partial \F |_{L^2_+} (u(t)) \, | \dot{u}(t) |_{L^2_+} \le \| u_0  \|^2_{L^2} . $$

\bigskip
\noindent {\bf A quasi-metric for \boldmath $L^2_\tau$.} In the setting of  \S\ref{s.pen} it is natural to introduce 
$$d_\tau ( u, v) = | u - v |_{L^2_\tau} = \left( \int_\Omega \psi_\tau (u-v) \, dx \right)^{1/2} 
	\quad \text{where } \quad 	\psi_\tau ( u-v ) = 
		\begin{cases}
		(u-v)^2 & \text{if $u \ge v$} \\
		\alpha(\tau) (u-v)^2 & \text{if $u < v$.}
		\end{cases} 
$$
It is not difficult to check that $d_\tau$ is a quasi-metric in $L^2(-1,1)$, i.e., that
$$
	d_\tau (u,v) \ge 0 , \quad d_\tau( u,v) = 0 \Leftrightarrow  u=v , \quad d_\tau(u,v) \le d_\tau(u,z) + d_\tau(z,v) ;
$$
however $d_\tau$ is not a metric since it is not symmetric. First, note that the metric derivative coincides with $ | \dot{v} |_{L^2_\tau}$; indeed, if $v \in AC ( 0,T ; L^2 )$ then, in every differentiability point $t$, we can write
\begin{align*}
	\lim_{h \to 0} \frac{d_\tau ( v(t+h) , v(t) ) }{h} & %= \lim_{h \to 0} \frac{\Psi^{1/2}_\tau ( v(t+h) - v(t))}{h}
	= \lim_{h \to 0} \left(  \int_\Omega \frac{\psi_\tau (v(t+h) -v(t)) }{h^2} \, dx \right)^{1/2} \\
	& = \lim_{h \to 0} \left( \int_\Omega \psi_\tau \left( \frac{v(t+h) - v(t)}{h}  \right) dx \right)^{1/2} 
	= \left( \int_\Omega \psi_\tau \left( \dot{v} (t) \right) dx \right)^{1/2} = | \dot{v} |_{L^2_\tau} .
\end{align*}
Therefore the slope of $\F ( t, \cdot)$ with respect to the quasi-metric $d_\tau$ actually coincides with the slope $| \partial \F |_{L^2_\tau} ( t, \cdot)$ defined in \S\ref{s.pen}. In this case, following the arguments of \cite[Theorem 1.2.5]{AmbrosioGigliSavare05}, the slope $|\partial \F|_{L^2_\tau} ( t, \cdot)$ turns out to be a strong upper gradient for $\F ( t, \cdot)$.

\section{A unilateral \boldmath{$L^2$}-subdifferential \label{app.B}}

In this section we propose a notion of unilateral subdifferential in $L^2$ and show its connection with the unilateral slope \eqref{e.slope}. For sake of simplicity we consider again an autonomous functional $\F$. 

\begin{definition} For $u \in H^1_0$ define
\begin{equation}\label{e.sub+}
	\partial_+ \F (u) := \{ \xi \in L^2 : \F (v) \ge \F (u) + \langle \xi  , v - u \rangle  \text{ for $v \in L^2$, $ v \ge u$} \} .
\end{equation}
% hence $| \partial_+ \F | (u)  =  \min   \{  \| \xi \|_{L^2} :  \xi \in \partial_+ \F (u)   \}    $
\end{definition}

As usual, $\partial_+ \F ( u) = \emptyset$ if $\F(u) = + \infty$ (i.e.~if $u \not\in H^1_0$). It is easy to check that $\partial_+ \F (u)$ is convex and closed. Note also that $\partial_+ \F ( u)$ is much larger than the (single valued) $\partial \F (u)$, indeed if $\zeta \in \partial \F (u)$ then every $\xi \ge \zeta$ belongs to $ \partial_+ \F (u)$. However, in analogy with the unconstrained setting, we can define the ``minimal selection''
$$ \partial_+^\circ \F ( u) = \argmin \{ \| \xi \| _{L^2} : \xi \in \partial_+ \F ( u) \} . $$ 
Note that, being $\partial_+ \F ( u)$ closed in $L^2$, the minimum is attained. Next lemma shows the natural relationship between unilateral slope and unilateral subdifferential.

\begin{lemma}\label{l.slope-sub} If $u \in H^1_0$ then $| \partial \F |_{L^2_+} ( u)  =  \min   \{  \| \xi \|_{L^2} :  \xi \in \partial_+ \F (u)   \}$. Moreover, if $\partial_+ \F (u) \neq \emptyset$ then $\partial^\circ_+ \F (u) = \{ - [ A u + f ]_+ \}$.
%\begin{equation}\label{e.slope-sub}
%\end{equation}
\end{lemma}

\proof  In the definition of slope \eqref{e.slope} it is not restrictive to consider $v \ge u$ such that $\F (v) \le \F(u)$. In this case, for $\xi \in \partial_+ \F (u)$ we have 
$$
	[\F(v) - \F(u)]_- = \F(u) - \F(v) \le - \langle \xi , v- u \rangle \le \| \xi \|_{L^2} \, \| v -u \|_{L^2} .
$$
Hence $| \partial \F |_{L^2_+} (u)  \le \|  \xi \|_{L^2}$ for every $\xi \in \partial_+ \F (u)$ from which 
$$  | \partial \F |_{L^2_+} (u)  \le  \min   \{  \|  \xi \|_{L^2} :  \xi \in \partial_+ \F (u)   \}   .   $$
The previous inequality holds also when $\partial_+ \F (u) = \emptyset$, in which case the right-hand side is infinite. 

Let us prove that $ | \partial \F |_{L^2_+} (u)  \ge  \min   \{  \|  \xi \|_{L^2} :  \xi \in \partial_+ \F (u)   \}$. 
If $ | \partial \F |_{L^2_+} (u)$ is infinite there is nothing to prove. Otherwise, $u \in H^1_0$ and by Lemma \ref{l.slope-diff}
$$
     | \partial \F |_{L^2_+} (u)  =  \sup  \left\{ - d\F( u) [z] :  z \in H^1_0 , \, z \ge 0 , \, \| z \|_{L^2} \le 1 \right\}   < +\infty .
$$
By Lemma \ref{l.measure} it follows that $- d \F (u) = A u + f $ is a Radon measure $\mu$ with positive part $\mu_+$ in $L^2$ and 
$ | \partial \F |_{L^2_+} (u)  = \| \mu_+ \|_{L^2}$. To conclude, it is enough to show that $- [Au + f]_+ = - \mu_+ \in \partial_+ \F ( u)$:  
by convexity, for every $v \in H^1_0$ with $v \ge u$ we have
\begin{align*}
	\F (u) - \F (v) \le - d \F (u) [ v -u] =  ( \mu ,  v - u ) \le \langle \mu_+ , v -u \rangle ,
\end{align*}
which concludes the proof.  \qed

%Invoking Lemma \ref{l.measure} together with Lemma \ref{l.slope-diff} and \eqref{e.dF} we get the following Corollary ...
%
%\begin{corollary} \label{c.B1} If  $| \partial_+ \F | (t, u) < +\infty$ then $\Delta u + f (t)$ is a Radon measure with $[ \Delta u + f (t)]_+ \in L^2$ and $\partial^\circ_+ \F ( t, u) = \{ - [\Delta u + f(t) ]_+\}$. % in particular $| \partial_+ \F | (t, u) = \| [ \Delta u + f (t)]_+ \|_{L^2}$. 
%\end{corollary}

\separe

Within this setting, the parabolic problem \eqref{e.pde} would read
$$
	\begin{cases}
		\dot{u} (t) = - \partial_+^\circ \F ( u (t) ) & \text{in $L^2$ for a.e.~$t \in (0,T)$} \\
		u (0) = u_0 .
	\end{cases}
$$
However, besides the non-uniqueness issues for \eqref{e.pde}, we have the following remark. %$\partial_+^\circ \F (u)$ is not a monotone operator in $L^2$. 

\begin{remark} $\partial_+^\circ \F (u)$ and $\partial_+ \F (u)$ are not monotone in $L^2$. %s the following example shows. 
\end{remark}

For $f = -1$ let 
$ \F (u) = \tfrac12 \int_{(-1,1)} | u' |^2 + u \, dx$ and $u_0=0$. Clearly $\F(u) \ge \F (u_0)$ for every $u \ge u_0$, hence  $0 = \xi_0 \in \partial_+ \F (u_0)$. Clearly $\xi_0 \in \partial_+^\circ \F (u_0)$.

On other hand, let $u \in H^1_0$ with $u > 0$ in $(-1,1)$ and $[u''-1]_+ \not\equiv 0$. Then, by Lemma \ref{l.slope-sub} we have $\partial_+^\circ \F (u) = \{ - [u'' - 1 ]_+ \}$.
Denoting $\xi_u = - [u'' - 1 ]_+$ we have
$$
	\langle \xi_u - \xi_0 , u - u_0 \rangle = - \langle u'' -1 , u \rangle < 0  
$$
and thus $\partial_+^\circ \F (u)$ is not monotone.
As a consequnce, $\partial_+ \F (u)$ is not monotone, since it is larger than $\partial_+^\circ \F (u)$.

%\separe
%For $m \in \mathbb{N}$ (to be choosen later) let 
%$$
%	u (x) = 
%	\begin{cases}
%		n^2 x^2  & \text{if $|x|<\tfrac1n$} \\
%		( | x| -1) \frac{n}{1-n} & \text{otherwise,}
%	\end{cases} 
%	\qquad
%	\bar{u} (x) = 
%	\begin{cases}
%		1  & \text{if $|x|<\tfrac1n$} \\
%		( | x| -1) \frac{n}{1-n} & \text{otherwise.}
%	\end{cases} 
%$$
%In a similar way, given $\tau>0$ is easy to see that for $n \gg 1$ we have 
%\begin{gather*}
%	v \in \argmin \{ \F (t ,w) + \tfrac1\tau \| w - v \|^2_{L^2} : w \in H^1_0 (-1,1) , \, w \ge v  \} ,  \\
%	u \not \in \argmin \{ \F (t ,w) + \tfrac1\tau \| w - u \|^2_{L^2} : w \in H^1_0 (-1,1) , \, w \ge u  \}   ,
%\end{gather*}
%because $\F (t , \bar{u} ) + \frac1\tau \| \bar{u}  - u \|^2_{L^2} < \F ( t , u)$. As a consequence the resolvent operator $J_\tau [\cdot]$ can be expansive, indeed, by the above argument we have $J_\tau [v] = v = 0$ while  $J_\tau [u] \ge u \ge 0$ and $J_\tau [u] \neq u$,  hence 
%$$
%	\| J_\tau [ u ] - J_\tau [ v ] \|_{L^2}  > \| u -v \|_{L^2} \,.
%$$
%\end{example}
%
%\separe
%
%

\separe

% Indeed, dividing by $t_2 - t_1$ and pasing to the limit we would get $ | \dot{\F} (u(t) ) | \le   | \partial_+ \F| (u(t) ) \, \| u_0 \|_{L^2} $. In the counterexample we have $| \partial_+ \F| (u(t) ) = \| u_0 \|_{L^2} = (2L/3)^{1/2}$ while $ \dot{\F} (u(t) ) = (2 t /L) -  (2L/3)$. For $t \gg1$ the inequality $| (2 t /L) -  (2L/3) | \le (2L/3)$ is false.

% In general we have only the one sided upper gradient inequality ...

%\section{Representation}

%\section{Interpolation}

%\section{Riemann sums}
%
%The following lemma follows easily from the proof of the similar \cite[Lemma 4.12]{DalMasoFrancfortToader}.
%
%\begin{lemma} \label{l.DFT} ...
%\end{lemma}

\bibliographystyle{plain}
\bibliography{references} % use symbolic links

\end{document}